\documentclass{amsart}

\usepackage{latexsym,amsmath,amstext,mathrsfs,amsthm,amscd,amsopn,verbatim,amsrefs}

\usepackage{graphicx}
\usepackage[dvips]{color}
\usepackage{epsfig}
\usepackage{psfrag}
\usepackage[all]{xy}
\usepackage{amscd}
\usepackage{amssymb}

\newtheorem*{MainThmA}{Theorem A}
\newtheorem{Thm}{Theorem}[section]
\newtheorem{Prop}[Thm]{Proposition}
\newtheorem{Lem}[Thm]{Lemma}
\newtheorem{Cor}[Thm]{Corollary}

\theoremstyle{remark}
\newtheorem{Rem}[Thm]{Remark}
\newtheorem*{Rem*}{Remark}

\theoremstyle{definition}

\newtheorem{Def}[Thm]{Definition}
\newtheorem{Exa}[Thm]{Example}

\newcommand{\mg}{\mathfrak{g}}
\renewcommand{\bul}{{\bullet}}
\renewcommand{\d}{{\rm d}}
\renewcommand{\H}{{\rm H}}
\newcommand{\C}{{\rm C}}
\renewcommand{\S}{{\rm S}}
 \newcommand{\U}{{\rm U}}
\newcommand{\norm}[1]{|\!|#1|\!|}

\title[Compatibility with cap-products in Tsygan's formality]
{Compatibility with cap-products in Tsygan's formality and homological Duflo isomorphism}
\author{Damien Calaque}
\address{Universit\'e de Lyon, Universit\'e Lyon 1, CNRS, UMR5208, Institut Camille Jordan,
43 blvd du 11 novembre 1918, F-69622 Villeurbanne-Cedex, France}
\email{calaque@math.univ-lyon1.fr}
\author{Carlo A. Rossi}
\address{Department of mathematics, ETH Zurich, 8092 Zurich, Switzerland}
\email{carlo.rossi@math.ethz.ch}
\thanks{The result of this paper were mainly obtained when D.C.~was working in ETH (on leave 
of absence from Universit\'e Lyon 1). His research was fully supported by the European Union 
thanks to a Marie Curie Intra-European Fellowship (contract number MEIF-CT-2007-042212). }

\begin{document}

\begin{abstract}
In this paper we prove, with details and in full generality, that the isomorphism induced on tangent 
homology by the Shoikhet-Tsygan formality $L_\infty$-quasi-isomorphism for Hochschild chains is compatible 
with cap-products. 
This is a homological analog of the compatibility with cup-products of the isomorphism induced on tangent 
cohomology by Kontsevich formality $L_\infty$-quasi-isomorphism for Hochschild cochains. 

As in the cohomological situation our proof relies on a homotopy argument involving 
a variant of {\bf Kontsevich eye}. In particular we clarify the r\^ole played by the 
{\bf I-cube} introduced in \cite{CR1}.  

Since we treat here the case of a most possibly general Maurer-Cartan element, not forced to be a bidifferential 
operator, then we take this opportunity to recall the natural algebraic structures on the pair of Hochschild 
cochain and chain complexes of an $A_\infty$-algebra. In particular we prove that they naturally inherit the 
structure of an $A_\infty$-algebra with an $A_\infty$-(bi)module. 
\end{abstract}

\setcounter{tocdepth}{1}

\maketitle

\tableofcontents

\section*{Introduction}

Given a (possibly curved) $A_\infty$-algebra $(A,\gamma_0,\gamma_1,\gamma_2,\dots)$, it is known that its Hochschild 
cochain complex $C^\bullet(A,A)$ is naturally a (non curved) $A_\infty$-algebra, with structure maps ${\rm d}_{\gamma,k}$ 
being defined thanks to the famous brace operations \cite{GV} introduced by Gerstenhaber and Voronov: 
\begin{equation}\label{eq-Ainfty1}
{\rm d}_{\gamma,1}(P):=\sum_i(\gamma_i\{P\}\mp P\{\gamma_i\})\qquad\textrm{and}\qquad
{\rm d}_{\gamma,k}(P_1,\dots,P_k):=\sum_{i}\gamma_i\{P_1,\dots,P_k\}\quad(k\geq2)\,.
\end{equation}
This statement can be reformulated and proved using $B_\infty$-algebras \cite{GJ} and twisting procedure 
for them with respect to Maurer-Cartan elements, following Getzler and Jones. 
Namely, given a $B_\infty$-algebra $(B,{\rm d},\mathrm m)$ and a 
{\bf Maurer-Cartan element} (shortly MCE) $\gamma$, i.e.~a degree $1$ element in $B$ satisfying the 
{\bf Maurer-Cartan equation}
$$
{\rm d}(\gamma)+\mathrm m(\gamma,\gamma)={\rm d}_1(\gamma)+\mathrm m_{1,1}(\gamma,\gamma)=0\,,
$$
then there is a new $B_\infty$-algebra $(B,{\rm d}_\gamma,\mathrm m)$ with 
$$
{\rm d}_\gamma:={\rm d}+\mathrm m(\gamma\otimes\bullet)-\mathrm m(\bullet\otimes\gamma)\,.
$$
For any graded vector space $A$ the brace operations on $B=\underline{{\rm End}}(A)$, with 
$$
\underline{{\rm End}}(A):=\bigoplus_{n\geq0}{\rm Hom}(A^{\otimes n},A)[1-n]\,,
$$
define a $B_\infty$-algebra structure on $B$ such that a MCE $\gamma$ tantamounts to a curved $A_\infty$-algebra structure on $A$, and the structure maps ${\rm d}_{\gamma,k}$ of ${\rm d}_\gamma$ 
are precisely given by \eqref{eq-Ainfty1}. 
From this formalism it is clear that two homotopy equivalent $A_\infty$-algebra structures 
on $A$ induce homotopy equivalent $A_\infty$-algebra structures on $B=\underline{{\rm End}}(A)$. 
All this is recalled in the first Section of the paper. 

\medskip

The first aim of the present paper is to develop a similar machinery for Hochschild chains of a curved 
$A_\infty$-algebra $(A,\gamma_0,\gamma_1,\gamma_2,\dots)$. Unfortunately, things do not appear to go as easily 
as in the case of cochains. We can nevertheless prove that there is an $A_\infty$-bimodule structure on the Hochschild 
chain complex $C_{-\bullet}(A,A)$ (with reversed grading), over the $A_\infty$-algebra $C^\bullet(A,A)$. 
To do so, we prove that there are two distinct left $B_\infty$-actions of $\underline{{\rm End}}(A)$ on 
$$
A\underline{\otimes}A:=\bigoplus_{n\geq0}A\otimes A^{\otimes n}[n]\,.
$$
Then, to any (curved) $A_\infty$-algebra structure on $A$, we define the $A_\infty$-bimodule structure 
on $A\underline{\otimes}A=C_{-\bullet}(A,A)$, as usual, as the adjoint action of the MCE $\gamma$. 
Being easier to say than to do, the above claim requires some work, and 
to introduce new notions such as $B_\infty$-(bi)modules. 
This is the subject of Section 2, which can be viewed as the explanation of the sketch of a construction of Tamarkin--Tsygan~\cite{TT} regarding dualities between Hochschild cochain and chain complex.

\bigskip

The above constructions extend to the following setting: a smooth real manifold $X$ and a commutative 
and unital differential graded algebra (shortly, DGA) $(\mathfrak m,{\rm d}_{\mathfrak m})$ splitting as 
$\mathfrak m=\mathfrak n\oplus R$, with $\mathfrak n$ a (pro)nilpotent ideal and $R$ a unital subalgebra 
concentrated in degree $0$. Then the complex of $\mathfrak{m}$-valued polydifferential operators 
$D_{\rm poly}^{\mathfrak m}(X)$, is naturally a $B_\infty$-algebra in which MCEs are deformations of 
the DGA $C^\infty(X,\mathfrak m)$ as an $A_\infty$-algebra over $(\mathfrak{m},{\rm d}_{\mathfrak m})$. 
\footnote{In the following, we assume that $\mathfrak m$ is bounded below as a graded vector space: 
$\mathfrak m^k=\{0\}$ for $k<\!\!<0$. 
Moreover, tensor products with $\mathfrak m$ have to be understood as completed tensor 
products with respect to the $\mathfrak n$-adic topology. }

Similarly, the complex of $\mathfrak{m}$-valued Hochschild chains $C^{{\rm poly},\mathfrak{m}}(X)$ 
with reversed grading (see e.g.~\cite{Dol} for a precise definition) naturally carries two distinct left 
$B_\infty$-actions (of $D_{\rm poly}^{\mathfrak m}(X)$), for which any MCE in $D_{\rm poly}^{\mathfrak m}(X)$ 
induces an $A_\infty$-bimodule structure on $C^{{\rm poly},\mathfrak{m}}(X)$. 

\medskip

Then we recall from \cite{K} that there exists an $L_\infty$-quasi-isomorphism $\mathcal{U}$ from 
the differential graded Lie algebra (shortly, DGLA) $T_{\rm poly}^{\mathfrak m}(X)$ of $\mathfrak{m}$-valued polyvector fields to the 
DGLA $D_{\rm poly}^{\mathfrak m}(X)$. Therefore, given a MCE $\gamma$ in 
$T_{\rm poly}^{\mathfrak n}(X)$ one obtains a chain map 
$$
\mathcal U_{\gamma,1}\,:\,\Big(T_{\rm poly}^{\mathfrak m}(X),{\rm d}_{\mathfrak m}+[\gamma,\ ]\Big)\,
\longrightarrow\,\Big(D_{\rm poly}^{\mathfrak m}(X),{\rm d}_{\mathfrak m}+[\widetilde\gamma,\ ]\Big)\,,
$$
where (below $\mu$ denotes the standard commutative product on $C^\infty(X,\mathfrak{m})$)
$$
\widetilde{\gamma}:=\mu+\sum_{n\geq1}\frac1{n!}\mathcal U_n(\underbrace{\gamma,\dots,\gamma}_{n\textrm{ times}})\,.
$$
Moreover, Kontsevich claimed and sketchily proved in \cite[Section 8]{K} (see \cite{MT,CR} for detailed 
proofs in particular cases) that $\mathcal U_{\gamma,1}$ is {\bf compatible with cup-products} in the sense that it 
induces an {\bf algebra} isomorphism 
$$
\mathrm H^\bullet\Big(T_{\rm poly}^{\mathfrak m}(X),{\rm d}_{\mathfrak m}+[\gamma,\ ]\Big)\,
\tilde\longrightarrow\,\mathrm H^\bullet\Big(D_{\rm poly}^{\mathfrak m}(X),{\rm d}_{\mathfrak m}+[\widetilde\gamma,\ ]\Big)\,.
$$

Analogously, we recall from \cite{Sh} that there exists an $L_\infty$-quasi-isomorphism $\mathcal{S}$ from the 
differential graded Lie module (shortly, DGLM) $\mathcal A^{\mathfrak{m}}(X)$ of $\mathfrak{m}$-valued differential forms (with reversed grading) to the 
DGLM $C^{{\rm poly},\mathfrak{m}}(X)$. Therefore, given a MCE $\gamma$ as above one obtains a quasi-isomorphism 
$$
\mathcal S_{\gamma,0}\,:\,\Big(C^{{\rm poly},\mathfrak m}(X),{\rm d}_{\mathfrak m}+{\rm L}_{\widetilde\gamma}\Big)\,
\longrightarrow\,\Big(\mathcal A^{\mathfrak m}(X),{\rm d}_{\mathfrak m}+{\rm L}_\gamma\Big)\,.
$$

The second aim of the paper is to prove the following Theorem, which is a very natural generalization of \cite{CR1}: 
\begin{MainThmA}\label{t-capprod}
The quasi-isomorphism $S_{\gamma,0}$ is {\bf compatible with cap-products} in the sense that it induces an isomorphism 
of $\mathrm H^\bullet\Big(T_{\rm poly}^{\mathfrak m}(X),{\rm d}_{\mathfrak m}+[\gamma,-]\Big)$-modules 
\[
\mathrm H^\bullet\Big(C^{{\rm poly},\mathfrak m}(X),{\rm d}_{\mathfrak m}+{\rm L}_{\widetilde\gamma}\Big)\,
\tilde\longrightarrow\,\mathrm H^\bullet\Big(\mathcal A^{\mathfrak m}(X),{\rm d}_{\mathfrak m}+{\rm L}_\gamma\Big)\,.
\]
\end{MainThmA}
\begin{Rem*}
A.~Cattaneo pointed our attention on this possible generalization of our previous work, in which we made 
use of the geometry of the I-cube. 
In \cite{CR1} and previous versions of the present paper we made use of the geometry of the I-cube to prove 
that a homotopy exists for the cap-products (following Kontsevich's original idea \cite{K} for cup-products), 
but in the simpler cases (namely, when $\gamma$ is at most a bivector) we have considered, three boundary faces did not 
contribute. They actually do in the general context we consider in the present paper; nevertheless, we give here 
a cleaner proof, in which the I-cube is, at the end, not strictly needed. 
\end{Rem*}

\medskip

The proof of Theorem A requires several steps. As it is now usual in deformation quantization we first 
prove the desired result in the local situation $X=\mathbb{R}^d$. For this purpose we recall, in Section 3, 
the construction of Kontsevich's \cite{K} and Shoikhet's \cite{Sh} local formality maps. 
The main ingredients of both constructions are appropriate compactified configuration spaces and integrals 
of angle forms over them. We detail in particular two remarkable compactified configuration spaces which 
are of some use for the compatibility between cup and cap products: {\bf Kontsevich's eye} and the {\bf I-cube}. 

\medskip

We then review quickly in Section 4 the proof of the compatibility between cup products in our 
very general framework for $X=\mathbb{R}^d$. The main argument, of homotopical nature, was sketched 
by Kontsevich in \cite{K}, later clarified by Manchon and Torossian in \cite{MT} in the framework of 
deformation quantization, and finally adapted to the case of Q-manifolds in \cite{CR}. 
The globalisation of the compatibility between cup products was first seriously considered in \cite{CVdB}, 
and is addressed in Section 8. 

\medskip

The proof of the compatibility between cap products for $X=\mathbb{R}^d$ occupies Section 5, and is based 
on a homotopy argument very similar to the one of Section 4. Contrarily to what we first guessed in \cite{CR1}, 
the I-cube is not strictly needed for the proof, but definitely gives insight to understand how things work. 
Again, the question of globalisation is pushed-forward to the final Section of the paper. 

\bigskip

Before going through the globalisation of the previous results, we discuss three special cases of interest 
and an application. Namely, the cases of interest, detailed in Section 6, are the following ones. 
The first case is when $\mathfrak m=\mathbb{R}[\![\hbar]\!]$; it corresponds to Shoikhet's conjecture \cite{Sh}, 
which is originally motivated by deformation quantization, and is proved in \cite{CR1}. The second one 
is when the MCE $\gamma$ is of polyvector degree at most 1; then one can prove that so is its image 
$\mathcal U(\gamma)$, which can be interpreted in terms of a Fedosov connection and its Weyl curvature 
on a deformed algebra, following the terminology of \cite{CFT}. Finally, the third case of interest 
is when the MCE $\gamma$ is precisely a vector field; we are able to compute explicitly the quasi-isomorphisms 
$\mathcal U_{\gamma,1}$ and $\mathcal S_{\gamma,0}$ by means of a rooted Todd class $j(\gamma)$, following 
\cite{CVdB} (see also \cite{CR}). 

\medskip

In Section 7 we present an application of the third case described in the preceding Section. Shortly speaking, 
we prove a (co)homological {\em analogon} of the so-called Duflo isomorphism. Here the rooted Todd class $j(\gamma)$ 
is precisely the Duflo element that is used to modify the Poincar\'e-Birkhoff-Witt isomorphism. 
We remind the reader that the use of the compatibility between cup products to prove the Duflo isomorphism goes back 
to Kontsevich's seminal paper \cite{K}, where its cohomological extension was claimed (and of which one can find 
a complete proof in \cite{PT}). In \cite{CR1} we proved a version of the Duflo isomorphism on coinvariants, 
and the result presented in this Section ends the story by extending it to homology. 

\bigskip

The final Section of the paper is devoted to the proof of the main Theorem A. It is basically obtained 
by means of now standard globalisation methods. These methods were introduced by Fedosov \cite{Fedo} for the 
deformation quantization of symplectic manifolds, generalized by Cattaneo-Felder-Tomassini \cite{CFT} to the 
case of Poisson manifolds, and finally adapted (and popularized) by Dolgushev \cite{Dol2,Dol} to the context 
of the formality (both for cochains and chains). Our presentation follows closely \cite{Dol} and is quite sketchy, 
focusing essentially on the main specific points for the compatibility. We end the Section in explaining 
how the approach of Cattaneo-Felder-Tomassini is contained in this description. 

\begin{Rem*}
We finally mention that our main result can be obtained as a consequence of a very recent preprint 
\cite{DTT} of Dolgushev-Tamarkin-Tsygan, where they prove the formality of the homotopy calculus algebra 
of Hochschild (co)chains. 
Their proof is more conceptual and does not require to check compatibility with cup and cap products, 
as both are part of the generating operations of the coloured operad {\bf calc} of calculus algebras. 
Nevertheless, our approach seems to have the advantage of being, to some extent, computable 
(see e.g.~Subsection \ref{ss-6-3}). 
\end{Rem*}

\subsection*{Acknowledgements}

We thank Alberto Cattaneo, who raised the question of generalizing our results in a context 
where all the boundary faces of the I-cube could contribute, for his interest in our work. 
We also thank Giovanni Felder for many interesting discussions on this project. 
We finally thank Vasiliy Dolgushev for a particularly enlighting remark about the globalisation 
of the compatibility between cup products, and Dmitry Tamarkin for his kind help on a tricky 
point in Section 2. 

\section*{Notation}

Unless otherwise specified, we work over a field $k$ of characteristic zero: 
algebras, modules etc ... are over $k$. 

\medskip

A graded vector space means a $\mathbb{Z}$-graded $k$-vector space. 
The category of graded vector spaces is symmetric monoidal with non-trivial commutativity 
isomorphism $\sigma$ being given by the Koszul sign rule: 
$$
\sigma_{VW}:V\otimes W\,\longrightarrow\,W\otimes V\,,\,v\otimes w\,\longmapsto\,(-1)^{|v||w|}w\otimes v\,.
$$
Dealing with graded vector spaces and morphisms between them, this rule is always tacitly assumed. 

\medskip

If $(\mathcal O,{\rm d}_{\mathcal O},{\rm m}_{\mathcal O})$ is a (possibly colored) DG operad, then a 
{\it $\mathcal O$-algebra up to homotopy} is the data of a DG vector space $(V,{\rm d})$ together with a 
DG linear map $\rho:(\mathcal O,{\rm d}_{\mathcal O})\to\big(\underline{{\rm End}}(V),[{\rm d}_V,-]\big)$ 
such that $\rho\circ{\rm m}_{\mathcal O}$ is homotopic to 
${\rm m}_{\underline{{\rm End}}(V)}\circ(\rho\underline{\otimes}\rho)$. 
In particular, any $\mathcal O_\infty$-algebra is a $\mathcal O$-algebra up to homotopy, while the converse is 
false. 

\medskip

Exemplarily, an associative algebra up to homotopy is a DG vector space $(A,{\rm d}_A)$ together with a product 
${\rm m}_A:A\otimes A\to A$ such that 
$[{\rm d}_A,{\rm m}_A]=0$, and 
$({\rm m}_A\otimes{\rm id})\circ{\rm m}_A$ is homotopic to $({\rm id}\otimes{\rm m}_A)\circ{\rm m}_A$. 

\medskip

The homotopy is not considered as a part of the structure. 
The obvious notion of a morphism of $\mathcal O$-algebra up to homotopy can then be guessed by the reader. 
Exemplarily, a morphism of associative algebras up to homotopy is a graded linear map 
$f:A\to B$ such that $f\circ{\rm d}_A={\rm d}_B\circ f$, and $f\circ{\rm m}_A$ is homotopic to 
${\rm m}_B\circ (f\otimes f)$. 

\section{$B_\infty$-structure on Hochschild cochains}\label{s-1}

In this Section, we discuss in its generality the $B_\infty$-algebra structure on $\underline{\mathrm{End}}(E)$, 
for a graded vector space $E$, and the twisting procedure that allows one to deduce from this the $B_\infty$-structure 
on the Hochschild cochain complex of an ($A_\infty$-)algebra $A$. 
It has been first exploited by Getzler--Jones~\cite{GJ} and Gerstenhaber--Voronov~\cite{GV}, to which we refer for 
more details and for complete proofs; we will nonetheless write down explicitly certain formul\ae\ and some 
arguments, which will be helpful for upcoming computations. 

\subsection{$B_\infty$-algebras and twistings}\label{ss-1-1}

We consider a graded vector space $V$: for a homogeneous element $v$ in $V$, we denote by $|v|$ its degree.
The {\bf cofree, coassociative coalgebra with counit cogenerated by $V$} is the tensor coalgebra 
$T(V)=\bigoplus_{n\geq n} V^{\otimes n},\ V^{\otimes 0}=k$, with the natural coproduct, resp.\ counit, 
\begin{equation*}
\begin{aligned}
\Delta(v_1\otimes \cdots\otimes v_p)
&=\sum_{i=0}^p\left(v_1\otimes\cdots\otimes v_i\right)\otimes\left(v_{i+1}\otimes\cdots\otimes v_n\right),\ \text{resp.}\\
\varepsilon(v_1\otimes\cdots\otimes v_n)&=\begin{cases}
1, & n=0\\
0, & n\geq 1,
\end{cases}
\end{aligned}
\end{equation*}
where $v_0=v_{n+1}=1$, $1$ being the unit of the ground field $k$.
\begin{Def}\label{d-B_inf}
$V$ is called a {\bf $B_\infty$-algebra}, if there exist linear maps 
$$
\mathrm d:T(V)\to T(V)\quad\textrm{and}\quad\mathrm m:T(V)\otimes T(V)\to T(V)\,,
$$
such that the $6$-tuple $(T(V),\mathrm m,\Delta,\eta,\varepsilon,\mathrm d)$, where $\eta$ 
is the natural unit, is a DG bialgebra. 
\end{Def}
In other words, $\mathrm m$ is an associative product of degree $0$ on $T(V)$ and a morphism of coalgebras, 
i.e.\ the following identities hold true
\[
\mathrm m\circ (\mathrm m\otimes 1)=\mathrm m\circ (1\otimes \mathrm m),\ \Delta\circ \mathrm m=(\mathrm m\otimes \mathrm m)\circ (1\otimes \tau\otimes 1)\circ (\Delta\otimes \Delta),
\] 
where $\tau$ denotes the standard braiding in the category of graded vector spaces.

Further, $\mathrm d$ is a linear operator on $T(V)$ of degree $1$, which squares to $0$, and which is simultaneously 
a derivation w.r.t.\ $\mathrm m$ and a coderivation w.r.t.\ $\Delta$: more explicitly,
\[
\mathrm d^2=0,\ \mathrm d\circ \mathrm m=\mathrm m\circ (\mathrm d\otimes 1+1\otimes\mathrm d),\ \Delta\circ \mathrm d=(\mathrm d\otimes 1+1\otimes\mathrm d)\circ \Delta,
\]
tacitly assuming Koszul's sign rule.

The fact that $\mathrm m$ is a morphism of coalgebras, resp.\ $\mathrm d$ is a coderivation of degree $1$, 
implies that $\mathrm m$, resp.\ $\mathrm d$, is uniquely specified by its components 
\[
\mathrm m_{p,q}:V^{\otimes p}\otimes V^{\otimes q}\to V,\ \text{resp.}\ \mathrm d_p:V^{\otimes p}\to V.
\]
For the sake of clarity of upcoming computations, we need to write down explicitly the product $\mathrm m$ and the 
differential $\mathrm d$ in terms of its components: namely, $\mathrm m$, resp.\ $\mathrm d$, is determined via the 
formul\ae
\begin{align}
\label{eq-B_inf-prod}\mathrm m(v_1\otimes\cdots\otimes v_p\otimes \widetilde v_1\otimes\cdots\otimes\widetilde v_q)&=\sum_{l=1}^{p+q}\sum_{(\mu,\nu)\in \mathcal P_l(p)\times\mathcal P_l(q)\atop |\mu\vee\nu|=l} \sigma(\mu,\nu)\mathrm m_{\mu,\nu}(v_1\otimes\cdots\otimes v_p\otimes \widetilde v_1\otimes\cdots\otimes\widetilde v_q),\ \text{resp.}\\
\label{eq-B_inf-diff}\mathrm d(v_1\otimes\cdots\otimes v_p)&=\sum_{i=1}^p\sum_{j=0}^{p-i}(-1)^{\sum_{k=1}^j|v_k|} v_1\otimes\cdots\otimes\mathrm d_i(v_{j+1}\otimes\cdots\otimes v_{j+i})\otimes\cdots\otimes v_p.
\end{align}
Formula (\ref{eq-B_inf-prod}) needs some explanations, also for later computations.

For positive integers $l$, $p$, we define 
\[
\mathcal P_l(p)=\left\{(\mu_1,\dots,\mu_l)\in \mathbb Z^l:\ \mu_i\geq 0,\ |\mu|=\sum_{i=1}^l\mu_i=p\right\},
\]
the set of (generalized) partitions of $p$ into $l$ subsets; we observe that the entries of a generalized 
partition $\mu$ are {\bf not} ordered.
Furthermore, for positive integers $l$, $p$ and $q$, the pairing $\vee$ between $\mathcal P_l(p)$ and 
$\mathcal P_l(q)$ is defined via 
\[
(\mu\vee\nu)_i=\begin{cases}
1,\ &\mu_i,\nu_i\geq 1\\
\mu_i,& \mu_i\geq 1,\ \nu_i=0\\
\nu_i,& \mu_i=0,\ \nu_i\geq 1\\
0,& \mu_i=\nu_i=0.
\end{cases}
\]
For non-negative integers $l$, $p$ and $q$, such that $1\leq l\leq p+q$, a pair $(\mu,\nu)$ in 
$\mathcal P_l(p)\times \mathcal P_l(q)$, such that $|\mu\vee\nu|=l$, determines a linear map 
$\mathrm m_{\mu,\nu}$ from $V^{\otimes p}\otimes V^{\otimes q}$ to $V^{\otimes l}$, {\em via}
\[
\mathrm m_{\mu,\nu}(v_1\otimes\cdots\otimes \overline v_1\otimes\cdots)=\bigotimes_{i=1}^l \mathrm m_{\mu_i,\nu_i}\left((v_{\sum_{j=1}^{i-1} \mu_j+1}\otimes\cdots\otimes v_{\sum_{j=1}^{i} \mu_j})\otimes (\widetilde v_{\sum_{j=1}^{i-1} \nu_j+1}\otimes\cdots\otimes \widetilde v_{\sum_{j=1}^{i} \nu_j})\right),
\]
where the factors in the tensor product are ordered from $1$ to $l$ from the left to the right; 
if either $\mu_i=0$ or $\nu_i=0$, for some $i=1,\dots,l$, then we set $\mathrm m_{\mu_i,\nu_i}=\mathrm{id}$; 
if both indices are $0$, we set $\mathrm m_{0,0}=0$.
Finally, the sign $\sigma(\mu,\nu)$ is determined by Koszul's sign rule.

For later purposes, it is useful to write down explicitly the associativity condition for the product $\mathrm m$ in terms of its components:
\begin{equation}\label{eq-B_inf-ass}
\begin{aligned}
&\sum_{l=1}^{p+q}\sum_{(\mu,\nu)\in \mathcal P_l(p)\times\mathcal P_l(q)\atop |\mu\vee\nu|=l}\sigma(\mu,\nu)\mathrm m_{l,r}(\mathrm m_{\mu,\nu}(v_1\otimes\cdots\otimes v_p\otimes\overline v_1\otimes\cdots\otimes \overline v_q)\otimes \widetilde v_1\otimes\cdots\otimes \widetilde v_r)=\\
&=\sum_{l=1}^{q+r}\sum_{(\nu,\pi)\in \mathcal P_l(q)\times\mathcal P_l(r)\atop |\nu\vee\pi|=l}\sigma(\nu,\pi)\mathrm m_{p,l}(v_1\otimes\cdots\otimes v_p\otimes (\mathrm m_{\nu,\pi}(\overline v_1\otimes\cdots\otimes \overline v_q\otimes\widetilde v_1\otimes\cdots\otimes \widetilde v_r)).
\end{aligned}
\end{equation}
\begin{Rem}\label{r-gerst}
Writing down explicitly the previous families of identities for a few simple cases, we find that, 
if $V$ is a $B_\infty$-algebra, the binary operations 
\[
[v_1,v_2]:=\mathrm m_{1,1}(v_1\otimes v_2)-(-1)^{|v_1||v_2|}\mathrm m_{1,1}(v_2\otimes v_1)\quad\textrm{and}
\quad v_1\cup v_2:=\mathrm d_2(v_1\otimes v_2)\,,
\]
together with the differential $\mathrm d_1$, endow $V[1]$ with the structure of a 
Gerstenhaber algebra up to homotopy, with homotopies expressible via $\mathrm d_3$ 
(for the associativity of $\cup$), $\mathrm m_{1,2}$ (for the Leibniz rule between $\cup$ and $[,]$), 
$\mathrm m_{2,1}$ (for the Jacobi identity of $[,]$), and ${\rm m}_{1,1}$ (for the commutativity of $\cup$). 
\end{Rem}

Let now $V$ be a $B_\infty$-algebra in the sense of Definition~\ref{d-B_inf}.
\begin{Def}\label{d-MC-B_inf}
A {\bf Maurer--Cartan element} $\gamma$ for the $B_\infty$-algebra $V$ is an element of $V$ of degree $1$, 
which obeys the Maurer--Cartan equation
\begin{equation}\label{eq-MC-B_inf}
\mathrm d\gamma+\mathrm m(\gamma,\gamma)=0.
\end{equation}
\end{Def}
Since $\gamma$ belongs to $V$, it is obviously primitive in the bialgebra $T(V)$; further, 
(\ref{eq-MC-B_inf}) simplifies to 
\[
\mathrm d_1\gamma+\mathrm m_{1,1}(\gamma\otimes\gamma)=0.
\]
The MC equation (\ref{eq-MC-B_inf}) for the MC element $\gamma$, together with the primitivity of $\gamma$, 
implies that the map
\[
\mathrm d_\gamma=\mathrm d+\mathrm m(\gamma\otimes \bul)-\mathrm m(\bul\otimes \gamma),
\]
tacitly assuming Koszul's sign rule, defines a twisted $B_\infty$-structure on $V$, 
i.e.\ $(T(V),{\rm m},\Delta,\eta,\varepsilon,\mathrm d_\gamma)$ is a DG bialgebra. 

\subsection{$B_\infty$-algebra structure on the Hochschild cochain complex}\label{ss-1-2}

We consider a graded vector space $E$: to it, we associate 
\[
V=\underline{\mathrm{End}}(E):=\bigoplus_{n\geq0}\mathrm{Hom}(E^{\otimes n},E)[1-n]\,.
\]
\begin{Rem}
Regarding the grading on $V$, we use the following notation: since $E$ is graded, then any tensor power of 
$E$ is also naturally graded, as well as $\mathrm{Hom}(E^{\otimes n},E)$.
The degree referring to this grading will be denoted by $|\cdot|$.
Further, if $P$ is an element of $\mathrm{Hom}(E^{\otimes p},E)$ of degree $|P|$, then its degree in $V$ 
is called its {\bf (shifted) total degree} and is denoted by $\norm P$. We thus have 
\begin{equation}\label{eq-totdeg}
\norm P=p-1+|P|\,.
\end{equation}
\end{Rem}
A $B_\infty$-algebra structure on $V$ has been constructed explicitly by Getzler--Jones~\cite{GJ} and 
Gerstenhaber--Voronov~\cite{GV}: we review here its construction and some of its main features.

The differential $\mathrm d$ on $T(V)$ is the trivial one; the multiplication $\mathrm m$ on $T(V)$ is defined by components $\mathrm m_{p,q}$, which are non-trivial precisely when $p\leq 1$, with no restrictions on $q$: the unit axiom for $\mathrm m$ forces $\mathrm m_{0,q}$ to be equal to the identity map, while $\mathrm m_{1,q}$ is defined via
\begin{equation}\label{eq-brace}
\begin{aligned}
&\mathrm m_{1,q}(P,Q_1,\dots,Q_q)(e_1,\dots,e_n)=P\{Q_1,\dots,Q_q\}(e_1,\dots,e_n)=\\
&=\sum_{1\leq i_1\leq\cdots\leq i_q\leq n}(-1)^{\sum_{k=1}^q |\!|Q_k|\!|\left(i_k-1+\sum_{j=1}^{i_k-1}|e_j|\right)} P(e_1,\dots,Q_1(e_{i_1},\dots),\dots,Q_q(e_{i_q},\dots),\dots,e_n),
\end{aligned}
\end{equation}
with the previous grading conventions; in (\ref{eq-brace}), $n=p+\sum_{a=1}^q(q_a-1)$, $1\leq i_1$, $i_k+q_k\leq i_{k+1}$, $k=1,\dots,q-1$, $i_q+q_q-1\leq n$. 
It is not difficult to prove that the brace operations~\eqref{eq-brace} have (total) degree $0$.


It is useful to have a pictorial representation of certain operations: we depict an operator $P$ of with $p$ inputs 
(and one output) as a corolla with $p$ leaves, going from the bottom to the top.
Thus, the component $\mathrm m_{1,q}$ has the following graphical representation:
\bigskip
\begin{center}
\resizebox{0.25 \textwidth}{!}{\input{brace.pstex_t}}\\
\text{Figure 1 - Pictorial representation of the component $\mathrm m_{1,q}$} \\
\end{center}
\bigskip
The conditions for $V$ to be a $B_\infty$-algebra reduce to the associativity condition~\eqref{eq-B_inf-ass}, 
which simplifies in the present situation to
\begin{equation}\label{eq-brass-1}
\begin{aligned}
&(P\{Q_1,\dots,Q_q\})\{R_1,\dots,R_r\}=\\
&=\sum_{1\leq i_1\leq \cdots\leq i_q\leq r}(-1)^{\sum_{a=1}^q \norm{Q_a}\left(\sum_{b=1}^{i_a-1}\norm{R_b}\right)}
P\{R_1,\dots,R_{i_1-1},Q_1\{R_{i_1},\dots\},\dots,R_{i_q-1},Q_q\{R_{i_q},\dots\},\dots,R_r\},
\end{aligned}
\end{equation}



We recall from Remark \ref{r-gerst} that we have a bracket of total degree $0$ on $V$:
$$
[P_1,P_2]=P_1\{P_2\}-(-1)^{\norm{P_1}\norm{P_2}}P_2\{P_1\}\,.
$$
For $q=r=1$, Condition (\ref{eq-brass-1}) simplifies to
\[
(P_1\{P_2\})\{P_3\}=P_1\{P_2\{P_3\}\}+P_1\{P_2,P_3\}+(-1)^{\norm{P_2}\norm{P_3}}P_1\{P_3,P_2\}\,.
\]
Whence the bracket satisfies the Jacobi identity\footnote{This can be recovered from the fact 
(see Remark \ref{r-gerst}) that ${\rm m}_{2,1}$, which is the homotopy for the Jacobi identity in a general 
$B_\infty$-algebra, vanishes here. }, and thus $V$ is a DGLA with trivial differential.

Since the differential $\mathrm d$ is trivial, a MCE $\gamma$ satisfies the identity
\[
\gamma\{\gamma\}=\frac12[\gamma,\gamma]=0.
\]
\begin{Exa}\label{ex-Ainf}
Now we consider a (possibly curved) $A_\infty$-algebra $(A,\gamma_0,\gamma_1,\gamma_2,\dots)$. 
By abuse of terminology, we may say that 
$$
\gamma:=\gamma_0+\gamma_1+\gamma_2+\cdots\in\prod_{n\geq0}{\rm Hom}(A^{\otimes n},A)^{1-n}
$$
is a MCE (in the present situation the Maurer-Cartan equation \eqref{eq-MC-B_inf} makes sense since 
each of its homogeneous component is a finite sum). 

Therefore the twisting procedure of Subsection~\ref{ss-1-1} applies and one obtains a new $B_\infty$-structure 
on $\underline{{\rm End}}(A)$, such that 
$$
{\rm d}_{\gamma,1}(P)=[\gamma,P]
$$
is (up to a sign) the standard Hochschild coboundary operator for cochains of an $A_\infty$-algebra, and 
$$
P_1\cup_\gamma P_2:={\rm d}_{\gamma,2}(P_1\otimes P_2)=\gamma\{P_1,P_2\}
$$
defines a product which is associative up to homotopy\footnote{Recall that it is actually commutative up to homotopy. }. 
More precisely, if $p_i$ is the number of entries of $P_i$, then, 
for $p_1+p_2\leq p$, we have, by construction, 
\[
\begin{aligned}
P_1\cup_\gamma P_2(a_1,\dots,a_p)=
&\sum_{1\leq j_1,\ j_1+p_1\leq j_2\atop j_2+p_2-1\leq p}
(-1)^{\sum_{i=1}^2\norm{P_i}\left(j_i-1+\sum_{k=1}^{j_i-1}|a_k|\right)}\gamma_{p-p_1-p_2}\!
\left(a_1,\dots,P_1(a_{j_1},\dots),\dots\right.\\
&\phantom{\sum_{1\leq j_1\leq j_2\leq p}(-1)^{\sum_{i=1}^2\norm{P_i}\left(j_i-1+\sum_{k=1}^{j_i-1}|a_k|\right)}}
\left.\dots,P_2(a_{j_2},\dots),\dots,a_p\right)\,.
\end{aligned}
\]
\end{Exa}
\begin{Exa}
As a special case, if $A$ has the structure of a graded algebra, the associative product 
$\mu$ on $A$ is a MCE and the differential $\mathrm d_\mu$ has only two components: 
$-(-1)^{p-1}\mathrm d_{\mu,1}$, resp.\ $(-1)^{p_1(p_2-1)}\mathrm d_{\mu,2}$, is the standard 
Hochschild differential, resp.\ the standard product, on the Hochschild cochain complex of the 
algebra $A$. 
\end{Exa}

\subsection{A more general example of a $B_\infty$-algebra structure}\label{ss-1-3}

We may consider more generally a commutative DG algebra $(\mathfrak{m},{\rm d}_{\mathfrak{m}})$ as 
in the introduction. Its differential extends to a differential 
$\mathrm d_{\mathfrak{m}}:={\rm id}\otimes\mathrm d_{\mathfrak{m}}$ 
on $V:=\underline{{\rm End}}(E)\otimes\mathfrak m$ of total degree $1$, 
which further extends to a (co)differential on $T(V)$, which we denote, 
by abuse of notations, by the same symbol.

The brace operations defined on $\underline{{\rm End}}(E)$ naturally extend to $V$ in the following way: 
$$
P\otimes m\{Q_1\otimes n_1,\dots,Q_r\otimes n_r\}:=(-1)^\epsilon P\{Q_1,\dots,Q_r\}\otimes mn_1\cdots n_r\,,
$$
where the sign $(-1)^\epsilon$ is determined by the appropriate braiding (with corresponding Koszul's sign rule). 

Then, the construction of the previous Subsection can be repeated {\em verbatim}, except that we have the additional 
non zero structure map ${\rm d}_1:=\mathrm d_{\mathfrak{m}}$. 

Therefore the Maurer-Cartan equation reads 
$$
{\rm d}_{\mathfrak{m}}(\gamma)+\gamma\{\gamma\}=0
$$
and makes sense for a generalized element
\[
\gamma=\gamma_0+\gamma_1+\gamma_2+\cdots\in\prod_{n\geq0}\big({\rm Hom}(E^{\otimes n},E)\otimes\mathfrak{m}\big)^{1-n}\,.
\] 
Such MCEs are in bijection with $(\mathfrak{m},{\rm d}_\mathfrak{m})$-$A_\infty$-algebra structures 
on $E\otimes\mathfrak{m}$. 

We implicitly make use of this $B_\infty$-structures in Sections \ref{s-4} and \ref{s-5} below 
(see also the Introduction above). 

\section{$B_\infty$-structures on Hochschild chains}\label{s-2}

In this Section, we discuss two $B_\infty$-module structures on $E\underline\otimes E$, 
for a graded vector space $E$, and the twisting procedure that allows one to deduce from these two distinct left 
$B_\infty$-module structures on the Hochschild chain complex of an ($A_\infty$-)algebra $A$. 
We believe this clarifies and makes more explicit a construction roughly sketched by Tamarkin--Tsygan in~\cite{TT}. 

\subsection{$B_\infty$-bimodules}\label{ss-2-1}

We assume the graded vector space $V$ to be a $B_\infty$-algebra; we borrow the main notations 
from Subsection~\ref{ss-1-1}. Let $W$ be another graded vector space. 
\begin{Def}\label{d-B_inf-mod}
A {\bf $B_\infty$-bimodule} structure on $W$ (over $V$) is a $B_\infty$-algebra structure on 
$V\oplus W[-1]$ such that 
\begin{itemize}
\item $V$ is a $B_\infty$-subalgebra, 
\item all components of structure maps involving $W$ more than once are zero. 
\end{itemize}
\end{Def}
\begin{Rem}
To $W$, we may associate the bi-comodule $\overline W$ cogenerated by $W$, namely, 
$\overline W=T(V)\otimes W\otimes T(V)$, with left-, resp.\ right-, coaction $\Delta_L$, 
resp.\ $\Delta_R$, defined via
\[
\Delta_L=\Delta\otimes 1,\ \text{resp.}\ \Delta_R=1\otimes \Delta.
\]
A $B_\infty$-bimodule structure on $W$ w.r.t.\ the $B_\infty$-algebra structure on $V$ tantamounts to the data of linear maps 
\[
\mathrm b:\overline W\to \overline W\,,\
\mathrm m_L:T(V)\otimes \overline W\to \overline W\,,\ \mathrm m_R:\overline W\otimes T(V)\to \overline W\,,
\]
such that the $6$-tuple $(\overline W,\mathrm b,\mathrm m_L,\Delta_L,\mathrm m_R,\Delta_R)$ is a DG bi-(co)module over the DG bialgebra $(T(V),\rm m,\Delta,\eta,\varepsilon,\rm d)$. 

More precisely, $\mathrm m_L$ (resp.~$\mathrm m_R$) defines a left (resp.~right) action of $T(V)$ on $\overline W$, 
which is required to be a morphism of bi-comodules; moreover, the left and right actions are required to commute, 
$\mathrm b$ is required to square to $0$, and to be a bi-(co)derivation of the bi-(co)module structure on $\overline W$.

From this, we see that there is an obvious notion of left (resp.~right) $B_\infty$-module. 
\end{Rem}
Compatibility of $\mathrm b$, $\mathrm m_L$ and $\mathrm m_R$ with coalgebra and comodule structures implies 
that they are uniquely determined by their structure maps (i.e.\ their evaluation on homogeneous components, composed 
with the standard projection $\overline{W}\twoheadrightarrow W$): 
\begin{align*}
\mathrm b_{p,q}&:V^{\otimes p}\otimes W\otimes V^{\otimes q}\to W\,,\\
\mathrm m_L^{p,q,r}&:V^{\otimes p}\otimes V^{\otimes q}\otimes W\otimes V^{\otimes r}\to W\,,\\
\mathrm m_R^{p,q,r}&:V^{\otimes p}\otimes W\otimes V^{\otimes q}\otimes V^{\otimes r}\to W\,.
\end{align*}
Exemplarily, we write down the condition for $\mathrm m_L$ to be a left action w.r.t.\ $\mathrm m$ in terms of their respective components:
\begin{equation}\label{eq-L-act}
\begin{aligned}
&\sum_{l=1}^{p+q}\sum_{(\mu,\nu)\in\mathcal P_l(p)\times\mathcal P_l(q)\atop |\mu\vee\nu|=l}\sigma(\mu,\nu)\mathrm m_L^{l,r,s}(\mathrm m_{\mu,\nu}(v_1\otimes\cdots\otimes v_p\otimes \overline v_1\otimes \cdots\otimes\overline v_q)\otimes \widetilde v_1\otimes\cdots\otimes \widetilde v_r\otimes w\otimes \widehat v_1\otimes\cdots\otimes \widehat v_s)=\\
&=\sum_{l=1}^{q+r+s+1}\sum_{i=1}^l\sum_{(\mu_i,\nu_i)\in\mathcal P_l(q)\times\mathcal P_l(r+s+1)\atop |\mu_i\vee\nu_i|=l,\ w\in \nu_{ii}}\sigma(\mu_i,\nu_i)\mathrm m_L^{p,i-1,l-i}(v_1\otimes\cdots\otimes \mathrm m_L^{\mu_i,\nu_i}(\overline v_1\otimes\cdots\otimes\widetilde v_1\otimes\cdots\otimes w\otimes\widehat v_1\otimes\cdots)),
\end{aligned}
\end{equation}
the notations are obvious generalizations of those introduced in Subsection~\ref{ss-1-1}.
Finally, we consider a MC element $\gamma$ for the $B_\infty$-algebra $V$ as in Definition~\ref{d-MC-B_inf}, 
Subsection~\ref{ss-1-1}: if $W$ is a $B_\infty$-bimodule as in Definition~\ref{d-B_inf-mod}, then $\gamma$ 
determines a twisted differential $\mathrm b_\gamma$ on $\overline W$ via
\[
\mathrm b_\gamma=\mathrm b+\mathrm m_L(\gamma\otimes\bul)-\mathrm m_R(\bul\otimes\gamma),
\]
tacitly using Koszul's sign rule, and the $6$-tuple 
$(\overline W,\mathrm b_\gamma,\mathrm m_L,\Delta_L,\mathrm m_R,\Delta_R)$ again defines a $B_\infty$-bimodule 
structure on $W$.
\begin{Exa}\label{ex-binf}
Let $V=\oplus_{n\in\mathbb{Z}}V_n$ be a $\mathbb{Z}$-graded $B_\infty$-algebra, whose 
structure maps are degree preserving. 
Then $V_0$ is obviously a $B_\infty$-algebra w.r.t.\ the restriction of $\mathrm m$. 
If we assume that $V_k=\{0\}$ when $k<-1$, then $V_{-1}[-1]$ is a $B_\infty$-bimodule 
over $V_0$, with left, resp.\ right, action $\mathrm m_L$, resp.\ $\mathrm m_R$, whose components are given by
\[
\mathrm m_L^{p,q,r}=\mathrm m_{p,q+1+r},\ \mathrm m_R^{p,q,r}=\mathrm m_{p+1+q,r}.
\] 
It is clear e.g.\ that~\eqref{eq-L-act} follows immediately from~\eqref{eq-B_inf-ass}, and analogous arguments imply the claim. 
Moreover, if $\gamma$ is a MCE in $V_0$, then the $B_\infty$-bimodule structure induced by the twisted ($\mathbb{Z}$-graded) $B_\infty$-algebra $(V,{\rm d}_\gamma,\mathrm m)$ on $V_{-1}[-1]$ 
obviously coincides with the twisted $B_\infty$-bimodule structure $({\rm b}_\gamma,{\rm m}_L,{\rm m}_R)$ 
on it. 
\end{Exa}

\subsection{Left $B_\infty$-module structures on the Hochschild chain complex}\label{ss-2-2}

Let $E$ be a graded vector space, to which we associate 
\[
E\underline{\otimes}E:=\bigoplus_{n\geq 0} E\otimes E^{\otimes n}[n]\,.
\]
We then define $F:=E\oplus E^*$ with the following additional $\mathbb{Z}$-grading: $E$, resp.~$E^*$, has $\mathbb{Z}$-degree $0$, resp.~$-1$. 
Therefore, $V:=\underline{\rm End}(F)$, becomes a $\mathbb{Z}$-graded $B_\infty$-algebra that satisfies the condition 
of Example \ref{ex-binf}. 
Explicitly,  
$$
V_0=\underline{\rm End}(E)\oplus\bigoplus_{p,q\geq0}{\rm Hom}(E^{\otimes p}\otimes E^*\otimes E^{\otimes q},E^*)[-p-q]
\qquad
\textrm{and}
\qquad
V_{-1}=\bigoplus_{n\geq0}{\rm Hom}(E^{\otimes n},E^*)[1-n]\,.
$$
In particular, $V_{-1}[-1]$ is canonically isomorphic to $(E\underline\otimes E)^*$: explicitly, the identification is given by
\[
\langle P(e_1,\dots,e_n),e_0\rangle=(-1)^{|e_0|\left(\sum_{i=1}^n|e_i|\right)}\langle \widetilde P,(e_0|\cdots|e_n)\rangle,\ P\in V_{-1},\ (e_0|\cdots|e_n)\in E\underline\otimes E.
\]
Moreover, there is an inclusion $P\mapsto\overline P$
$$
{\rm Hom}(E^{\otimes n+1},E)\hookrightarrow\bigoplus_{p+q=n}{\rm Hom}(E^{\otimes p}\otimes E\otimes E^{\otimes q},E),
$$
explicitly given by the formula
\begin{equation}\label{eq-incl}
\langle\overline P(e_1,\dots,e_i,\xi,e_{i+1},\dots,e_p),e_0\rangle
:=(-1)^{|\xi|\left(\norm P+\sum_{j=1}|e_j|\right)+\left(i+1+\sum_{j=0}^i|e_j|\right)\left(p-i+\sum_{j=i+1}^p|e_j|\right)}\langle\xi,P(e_{i+1},\dots,e_p,e_0,e_1,\dots,e_i)\rangle\,.
\end{equation}
We observe that cyclic permutations enter into the game explicitly at this step.
In turn, Formula~\eqref{eq-incl} induces an inclusion 
\begin{equation}\label{eq-inc}
\underline{\rm End}(E)\hookrightarrow
\underline{\rm End}(E)\oplus\bigoplus_{p,q\geq0}{\rm Hom}(E^{\otimes p}\otimes E^*\otimes E^{\otimes q},E^*)[-p-q]
\subset V_0,\ P\mapsto P+\overline P.
\end{equation}
Obviously, the identity morphism preserves the $B_\infty$-algebra structure.
On the other hand, we may compute, using Formula~\eqref{eq-incl}, the inclusion $\overline{P\{Q_1,\dots,Q_q\}}$, for $P$, $Q_i$, $i=1,\dots,q$, general elements of $\underline{\mathrm{End}}(E)$, and we get
\begin{equation}\label{eq-notcomp}
\begin{aligned}
\overline{P\{Q_1,\dots,Q_q\}}&=\sum_{i=1}^q \overline{P}\{Q_i,\dots,Q_q,Q_1,\dots,Q_{i-1}\}+\sum_{i=1}^q \overline{Q_i}\{\overline P\{Q_{i+1},\dots,Q_q,Q_1,\dots,Q_{i-1}\}.
\end{aligned}
\end{equation}
The two terms on the right hand-side of Identity~\eqref{eq-notcomp} need some explanations.
The cyclic permutations of the elements $Q_i$, $i=1,\dots,q$, appear evidently because of Formula~\eqref{eq-incl}: $\overline P\{Q_i,\dots,Q_q,\dots\}$, resp.\ $\overline{Q_i}\{\overline{P}\{Q_{i+1},\dots,Q_q,Q_1,\dots,Q_{i-1}\}\}$, acts non-trivially precisely on those terms, where the argument labelled by $e_0$ is placed between $Q_{i-1}$ and $Q_i$, resp.\ as an argument of $Q_i$. 
Finally, we observe that we have omitted the signs in Identity~\eqref{eq-notcomp}: these are easily obtained by Koszul's sign rule w.r.t.\ total degree.

In the special case where we consider only $P$ and $Q$, the defect of Inclusion~\eqref{eq-inc} to be a $B_\infty$-algebra morphism can be characterized in a nice way, namely
\[
\iota(P)\{\iota(Q)\}-\iota(P\{Q\})=[\overline P,\overline Q].
\]
Since the inclusion \eqref{eq-inc} is NOT a $B_\infty$-algebra morphism, then we do NOT obtain a 
$B_\infty$-bimodule structure on $E\underline{\otimes}E$ over $\underline{\rm End}(E)$. 
Nevertheless, as we will now explain, we will get two distinct left $B_\infty$-module structure, which we now explicitly describe. 

\medskip

The only non-trivial structure maps of the right $B_\infty$-module structure ${\rm m}_R$ on 
$V_{-1}[-1]$ over $V_0$ are 
$$
{\rm m}_R^{0,0,q}(P,Q_1,\dots,Q_q)=P\{Q_1,\dots,Q_q\}\,,
$$
for $P\in V_{-1}$ and $Q_1,\dots,Q_q\in\underline{\rm End}(E)\subset V_0$. In particular, we can 
see that the induced left $B_\infty$-bimodule structure ${\rm m}_{L,2}$ on $E\underline\otimes E$, over 
$\underline{\rm End}(E)$,\footnote{We have an induced left $B_\infty$-bimodule structure since that, for $P,Q_1,\dots,Q_q$ 
as above, $P\{\iota(Q_1),\dots,\iota(Q_q)\}=P\{Q_1,\dots,Q_q\}$. } has only non-trivial structure maps 
${\rm m}_{L,2}^{p,0,0}$ given as follows: for any $c=(e_0|\cdots|e_m)\in E\underline\otimes E$, 
\begin{equation}\label{eq-L-2}
\mathrm m_{L,2}^{p,0,0}(P_1,\dots,P_p,c)=\sum_{1\leq i_1\leq \cdots\leq i_p\leq m}
(-1)^{\sum_{a=1}^p\norm{P_a}\left(i_a-1+\sum_{d=0}^{i_a-1}|e_d|\right)}\left(e_0|\cdots|P_1(e_{i_1},\dots)|\cdots|P_p(e_{i_r},\dots)|\cdots|e_m\right)\,,
\end{equation}
where the summation is over indices $i_1,\dots,i_p$, such that $1\leq i_1$, $i_k+p_k\leq i_{k+1}$, 
$k=1,\dots,r-1$, $i_p+p_p-1\leq m$ (if $p=0$, $\mathrm m_{L,2}^{0,0,0}$ is the identity map).
The fact that summation is over all indices $1\leq i_1$ is a consequence of the duality between $V_{-1}$ and $(E\underline \otimes E)^*$, which highlights the special element $e_0$ in a given chain.

As for the brace operations, we have a pictorial representation for those structure maps: 
\begin{center}
\resizebox{0.32 \textwidth}{!}{\input{m_L-2.pstex_t}}\\
\text{Figure 3 - Pictorial representation of a component of $\mathrm m_{L,2}^{p,0,0}$} \\
\end{center}

\medskip

The only non-trivial structure maps of the left $B_\infty$-module structure ${\rm m}_L$ on $V_{-1}[-1]$ are 
$$
{\rm m}_L^{1,q,r}(P,Q_1,\dots,Q_q,S,R_1,\dots,R_r)=P\{Q_1,\dots,Q_q,S,R_1,\dots,R_r\}\,,
$$
where $P\in {\rm Hom}(E^{\otimes k}\otimes E^*\otimes E^{\otimes l},E^*)$, 
$Q_a$ and $R_b$, $a=1,\dots,q$, $b=1,\dots,r$ elements of $\underline{\rm End}(E)$ and $S\in V_{-1}$. 
In particular, we get an induced left $B_\infty$-bimodule structure ${\rm m}_{L,1}$ on 
$E\underline\otimes E$, over $\underline{\rm End}(E)$, with only non-trivial 
structure maps ${\rm m}_{L,1}^{1,q,r}$, via
\[
\left\langle {\overline P\{R_1,\dots,R_r,S,Q_1,\dots,Q_q\}}\widetilde\ ,c\right\rangle=\left\langle \widetilde S,\mathrm m_{L,1}^{1,q,r}(P,Q_1,\dots,Q_q,c,R_1,\dots,R_r)\right\rangle,
\]
where $S$, resp.\ $c$, is a general element of $V_{-1}$, resp.\ $E\underline \otimes E$, such that the previous expression makes sense. 

The brace identities~\eqref{eq-brass-1}, Subsection~\ref{ss-1-2}, together with Identity~\eqref{eq-notcomp}, imply that the previous formula yields a left $B_\infty$-action: still, we observe that two dualizations are hidden in the previous formula, the first one in the inclusion $P\mapsto \overline P$, the second one between $V_{-1}$ and $(E\underline\otimes E)^*$. 

\begin{center}
\resizebox{0.55 \textwidth}{!}{\input{m_L-1.pstex_t}}
\text{Figure 4 - Pictorial representation of the component $\mathrm m_{L,1}^{1,q,r}$} \\
\end{center}
For $P$, $Q_a$, $R_b$ as before in $\underline{\rm End}(E)$, and $c=(e_0|\cdots|e_m)$ in $E\underline\otimes E$, we have the explicit form
\begin{equation}\label{eq-L-sg}
\begin{aligned}
&\mathrm m_{L,1}^{1,q,r}(P,Q_1,\dots,Q_q,c,R_1,\dots,R_r))=\\
&=\sum_{l\leq j_1\leq\cdots\leq j_q\leq m\atop 1\leq k_1\leq\cdots\leq k_r\leq n-1}(-1)^{\left(l+\sum_{j=0}^{l-1}|e_j|\right)\left(m-l+1+\sum_{j=l}^m|e_i|\right)+\sum_{a=1}^q |\!|Q_a|\!|\left(j_a-l+\sum_{d=l}^{j_a-1}|e_d|\right)+\sum_{b=1}^r|\!|R_b|\!|\left(k_b-l-1+\sum_{f=k_b}^{l-1}|e_f|\right)}\\
&\left(P\left(e_l,\cdots,Q_1(e_{j_1},\dots),\dots,Q_q(e_{j_q},\dots),\dots,e_0,\dots,R_1(e_{k_1},\dots),\dots,R_r(e_{k_r},\dots),\dots\right)|e_n|\cdots|e_{l-1}\right),
\end{aligned}
\end{equation}
and the indices in the summation satisfy $j_i+q_i\leq j_{i+1}$, $i=1,\dots,q-1$, $j_q+q_q-1\leq m$, and $k_i+r_i\leq k_{i+1}$, $i=1,\dots,r-1$, $k_r+r_r-1\leq n-1$.

Now we consider a (possibly curved) $A_\infty$-algebra on $(A,\gamma_0,\gamma_1,\gamma_2,\dots)$. 
We allow ourselves the same abuse of language as in Example \ref{ex-Ainf} and consider the formal sum 
$\gamma=\gamma_0+\gamma_1+\gamma_2+\cdots$ as a MCE of the $B_\infty$-algebra $\underline{\rm End}(A)$. 
Then we observe that, even if $\iota$ is NOT a $B_\infty$-algebra morphism,  
$$
\iota(\gamma)=\sum_{n\geq0}\gamma_n+\sum_{p,q\geq0}\gamma_{p,q}\,,
$$ 
defines a MCE (again, by abuse of notation) in $\underline{\rm End}(A\oplus A^*)$. 
Namely, $\gamma_{p,q}$ ($p,q\geq0$) are the structure maps of the natural $A_\infty$-bimodule structure on 
$A^*$.\footnote{To prove that it is truly a MCE, we simply recall that, if one has an $A_\infty$-algebra $B$ 
together with an $A_\infty$-bimodule $M$, then by definition $B\oplus M$ is an $A_\infty$-algebra.}

We may then apply the twisting procedure sketched at the end of Subsection~\ref{ss-1-1} to $\underline{\mathrm{End}}(A\oplus A^*)$ w.r.t.\ the MCE $\iota(\gamma)$.
Following the same lines of reasoning as above, we get the following
\begin{Thm}\label{t-MC-mod}
$A\underline{\otimes}A$ has two distinct left $B_\infty$-module structures over 
$(\underline{{\rm End}}(A),{\rm d_\gamma},{\rm m})$, given by ${\rm m}_{L,i}$, $i=1,2$, and the formula  
$$
{\rm b_\gamma}:={\rm m}_{L,1}(\gamma\otimes\bullet)-{\rm m}_{L,2}(\gamma\otimes\bullet),
$$
specifies a degree $1$ operator, which squares to $0$ (i.e.\ an $A_\infty$-module structure on $A\underline\otimes A$ over $\underline{\mathrm{End}}(A)$). 
\end{Thm}
We only observe that the twisting procedure, as in Subsection~\ref{ss-1-1}, cannot be applied {\em verbatim} in the present situation because of Identity~\eqref{eq-notcomp}.
Still, the same identity implies that $\mathrm b_\gamma$ squares to $0$, as can be verified by a direct computation.
\begin{Rem}
It follows directly from Subsection \ref{ss-1-3} that this construction generalizes to the situation 
where we tensorize $\underline{\rm End}(E)$ and $E\underline{\otimes}E$ by a commutative DG algebra 
$(\mathfrak{m},{\rm d}_{\mathfrak{m}})$ as in the introduction. 
\end{Rem}

\subsection{$T$-algebra structure on Hochschild (co)homology}\label{ss-2-3}

First, we observe that from the very definition of a $B_\infty$-bimodule, and according to the fact that 
any $B_\infty$-algebra is a Gerstenhaber algebra up to homotopy (see Remark \ref{r-gerst}), we have the following: 
\begin{quotation}
{\it For any $B_\infty$-algebra $V$ together with a $B_\infty$-bimodule $W$, the pair $(V[1],W)$ 
naturally inherits the structure of Gerstenhaber algebra and module up to homotopy. }
\end{quotation}

The crucial point is that, on the Hochschild chain complex of an $A_\infty$-algebra we only have two left 
$B_\infty$-module structures, having the same differential, but NOT satisfying the axioms of a 
$B_\infty$-bimodule (see the previous subsection). We therefore do NOT have the structure of a 
Gerstenhaber module up to homotopy. 
Nevertheless, we prove below that we have something close to a Gerstenhaber module; namely, a {\it $T$-algebra}. 

\medskip

We first recall the definition of a $T$-algebra (or {\it precalculus} following the 
recent terminology of \cite{DTT}) from \cite{TT}. 
\begin{Def}\label{def-Talg}
A {\bf $T$-algebra} is a pair $(V,W)$ of a Gerstenhaber algebra $(V,\cup,[,])$ and a 
graded vector space $W$ together with 
\begin{itemize}
\item an action $\cap$ of the GA $(V,\cup)$, turning $(W,\cap)$ into a GM, 
\item an action ${\rm L}$ of the GLA $(V[-1],[,])$, turning $(W,{\rm L})$ into 
a GLM, 
\end{itemize}
such that the following identities hold true for any $v_1,v_2\in V$ and $w\in W$: 
\begin{align}
\label{eq-T-1}\mathrm L_{v_1}(v_2\cap w)&=[v_1,v_2]\cap w+(-1)^{|v_1|(|v_2|-1)} v_2\cap({\rm L}_{v_1}w),\\
\label{eq-T-2}{\rm L}_{v_1\cup v_2}w&={\rm L}_{v_1}(v_2\cap w)+(-1)^{|v_1|}v_1\cap(L_{v_2}w)\,.
\end{align}
We may say, by abuse of language, that $W$ is a {\bf $T$-module} over the Gerstenhaber algebra $V$. 
\end{Def}
\begin{Rem}
It is worth mentioning that Identity~\eqref{eq-T-2} in Definition~\ref{def-Talg} can be re-written as 
\begin{equation}\label{eq-T-mod}
{\rm L}_{v_1\cup v_2}w=(-1)^{|v_2|(|v_1|-1)}v_2\cap\left({\rm L}_{v_1}w\right)+(-1)^{|v_1|}v_1\cap({\rm L}_{v_2}w)
-(-1)^{|v_1|}[v_1,v_2]\cap w\,,
\end{equation}
and we therefore see that a $T$-module is almost a Gerstenhaber module, the default being given by the last 
term in the r.h.s.~of \eqref{eq-T-mod}. 
The modified Identity~\eqref{eq-T-mod} will be particularly useful in the upcoming computations. 
\end{Rem}

We now prove that $A\underline{\otimes}A$ is a $T$-module up to homotopy over $\underline{{\rm End}}(A)$. 

We first recall that $\underline{{\rm End}}(A)$ is a graded Lie algebra with bracket 
$[P,Q]:=P\{Q\}-(-1)^{|\!|P|\!||\!|Q|\!|}Q\{P\}$. 
We observe that the same is true for 
$\underline{{\rm End}}(A\oplus A^*)$. 
This allows us to define a right Lie action of $\underline{{\rm End}}(A)$ onto $V_{-1}$. Namely, for 
$P\in\underline{{\rm End}}(A)$ and $Q\in V_{-1}$, we set
$$
{\rm R}_P(Q):=[Q,\iota(P)]=Q\{P\}-(-1)^{|\!|P|\!||\!|Q|\!|}\overline{P}\{Q\}\,.
$$
We now prove that the previous formula defines a right Lie action.
First of all, we evaluate explicitly $\rm R_{P_1}(\rm R_{P_2}Q)$, using the brace relations~\eqref{eq-brass-1}, Subsection~\ref{ss-1-2}:
\[
\begin{aligned}
\rm R_{P_1}(\rm R_{P_2}Q)&=Q\{P_2,P_1\}+Q\{P_2\{P_1\}\}+(-1)^{\norm{P_1}\norm{P_2}}Q\{P_1,P_2\}-\\
&\phantom{=}-(-1)^{\norm Q\norm{P_2}}\overline{P_2}\{Q,P_1\}-(-1)^{\norm Q\norm{P_2}}\overline{P_2}\{Q\{P_1\}\}-(-1)^{\norm Q\left(\norm{P_1}+\norm{P_2}\right)}\overline{P_2}\{P_1,Q\}-\\
&\phantom{=}-(-1)^{\norm{P_1}\left(\norm Q+\norm{P_2}\right)}\overline{P_1}\{Q\{P_2\}\}-(-1)^{\norm{P_1}\left(\norm Q+\norm{P_2}\right)}\overline{P_1}\{\overline{P_2}\{Q\}\}.
\end{aligned}
\]
A similar expression is obtained evaluating $\rm R_{P_2}(\rm R_{P_1}Q)$: summing up the two terms with the correct signs, we find
\[
\begin{aligned}
\left[\rm R_{P_1},\rm R_{P_2}\right](Q)&=Q\{P_2\{P_1\}\}-(-1)^{\norm{P_1}\norm{P_2}}Q\{P1\{P_2\}\}-\\
&\phantom{=}-(-1)^{\norm Q\norm{P_2}}\overline{P_2}\{Q,P_1\}-(-1)^{\norm Q\left(\norm{P_1}+\norm{P_2}\right)}\overline{P_2}\{P_1,Q\}-(-1)^{\norm Q\left(\norm{P_1}+\norm{P_2}\right)+\norm{P_1}\norm{P_2}}\overline{P_1}\{\overline{P_2}\{Q\}\}+\\
&\phantom{=}+(-1)^{\norm{P_1}\left(\norm Q+\norm{P_2}\right)}\overline{P_1}\{Q,P_2\}+(-1)^{\norm Q\left(\norm{P_1}+\norm{P_2}\right)+\norm{P_1}\norm{P_2}}\overline{P_1}\{P_2,Q\}+(-1)^{\norm Q\left(\norm{P_1}+\norm{P_2}\right)}\overline{P_2}\{\overline{P_1}\{Q\}\}.
\end{aligned}
\]
Obviously, the first two terms on the right hand-side sum up to $Q\{[P_2,P_1]\}$.
On the other hand, we consider $\overline{[P_2,P_1]}\{Q\}$: explicitly, in virtue of Formula~\eqref{eq-incl}, Subsection~\ref{ss-2-2},
\[
\overline{P_2\{P_1\}}\{Q\}=\overline{P_2}\{P_1,Q\}+(-1)^{\norm Q\norm{P_1}}\overline{P_2}\{Q,P_1\}+(-1)^{\norm{P_1}\norm{P_2}}\overline{P_1}\{\overline{P_2}\{Q\}\},
\]
and a similar formula holds true for $\overline{P_2\{P_1\}}\{Q\}$ with obvious due changes.
This yields $\left[\rm R_{P_1},\rm R_{P_2}\right](Q)=R_{[P_2,P_1]}(Q)$.
Dually, $\rm R$ defines a (left) graded Lie module structure on $A\underline{\otimes}A$, which we denote by ${\rm L}$. 

We have not yet considered the $A_\infty$-algebra structure on $A$, i.e. the MCE $\gamma$. 
Since we have a graded Lie algebra together with a graded 
Lie module, then we can twist them by the MCE $\gamma$ and obtain a DGLA 
$\big(\underline{{\rm End}}(A),[\gamma,\ ],[\ ,\ ]\big)$ together with a DGLM 
$\big(A\underline{\otimes}A,{\rm L}_\gamma,{\rm L}\big)$. 

Further, we recall that we have a product $\cup_\gamma:={\rm d}_{\gamma,2}$ which makes $\underline{{\rm End}}(A)[1]$ into 
a Gerstenhaber algebra up to homotopy. Analogously, we define an action (from the right) of 
$\underline{{\rm End}}(A)[1]$ onto $V_{-1}$ as follows: for $P\in\underline{{\rm End}}(A)$ and $Q\in V_{-1}$, we set
$$
Q\cap_\gamma P:=Q\cup_{\iota(\gamma)}\iota(P)=\iota(\gamma)\{Q,P\}=\bar{\gamma}\{Q,P\}\,.
$$
Since $\cup_{\iota(\gamma)}$ is associative up to homotopy, with homotopy being given by ${\rm d}_{\iota(\gamma),3}$, 
then $\cap_\gamma$ defines the structure of a right module up to homotopy on $V_{-1}$, over $\underline{{\rm End}}(A)[1]$.
Dually, we have a left module structure up to homotopy on $A\underline{\otimes}A$, 
which we again denote by $\cap_\gamma$. 

Moreover, the differential ${\rm d}_{\gamma,1}=[\gamma,\ ]$, resp.~${\rm b}_{\gamma,0,0}={\rm L}_\gamma$, 
is by definition compatible with the product $\cup_\gamma$, resp.~the action $\cap_\gamma$. 
\begin{Rem}\label{r-comm-actions}
A similar formula, where we switch $P$ and $Q$, defines accordingly a right action of $\underline{\mathrm{End}}(A)$ on $A\underline\otimes A$, which we will also denote by $\cap_\gamma$: one may think that this would lead to some confusion, but that both actions commute in the graded sense, since $\cup$ commutes up to homotopy, thus, later on, we will not distinguish between left and right action from the notational point of view. 
\end{Rem} 
\begin{Prop}
$\big(A\underline{\otimes}A,{\rm L}_\gamma,{\rm L},\cap_\gamma\big)$ is a $T$-module up to homotopy over the 
Gerstenhaber algebra up to homotopy 
$\big(\underline{{\rm End}}(A)[1],\rm d_{\gamma,1},[\ ,\ ],\cup_\gamma\big)$. 
\end{Prop}
\begin{proof}
By the above arguments and computations, it remains to prove the homotopical versions of Identities~\eqref{eq-T-1} and~\eqref{eq-T-2} in Definition~\ref{def-Talg}: in particular, we observe that we will prove the homotopical version of the modified Identity~\eqref{eq-T-mod}.

We will only write down the explicit homotopy formul\ae\ with signs: the computations leading to their proof make use of the brace identities~\eqref{eq-brass-1}, Subsection~\ref{ss-1-2}, and of Identity~\eqref{eq-notcomp}, Subsection~\ref{ss-2-2}, since $\mathrm L$ and $\cap$ can be described explicitly in terms of the brace operations on $V_{-1}$, and Identity~\eqref{eq-notcomp}, Subsection~\ref{ss-2-2} measures the failure of the two left $B_\infty$-actions $\rm m_{L,i}$, $i=1,2$, of being compatible.

Explicitly, we have the homotopy formul\ae, 
\[
\begin{aligned}
\mathrm L_{P_1}(P_2\cap_\gamma c)&=[P_1,P_2]\cap_\gamma c+(-1)^{\norm{P_1}(\norm{P_2}-1)}P_2\cap_\gamma \mathrm L_{P_1}c+(-1)^{\norm{P_1}}\left(\mathrm L_\gamma(\mathrm m_{L,1}^{1,1,0}(P_1,P_2,c))-\right.\\
&\phantom{=}\left.-\mathrm m_{L,1}^{1,1,0}(\mathrm d_{\gamma,1} P_1,P_2,c)-(-1)^{\norm{P_1}}\mathrm m_{L,1}^{1,1,0}(P_1,\mathrm d_{\gamma,1} P_2,c)-(-1)^{\norm{P_1}+\norm{P_2}}\mathrm m_{L,1}^{1,1,0}(P_1,P_2,\mathrm L_\gamma c)\right),
\end{aligned}
\]
and 
\[
\begin{aligned}
&\mathrm L_{P_1\cup_\gamma P_2}c+(-1)^{(\norm{P_1}-1)(\norm{P_2}-1)}\mathrm L_{P_2\cup_\gamma P_1}c=\\
&=(-1)^{\norm{P_1}}\left(P_1\cap_\gamma \mathrm L_{P_2}c+(-1)^{(\norm{P_1}-1)(\norm{P_2}+\norm{c}-1)}\mathrm L_{P_2}c\cap_\gamma P_1\right)+\\
&\phantom{=}+(-1)^{(\norm{c}-1)(\norm{P_2}-1)}\left(\mathrm L_{P_1}c\cap_\gamma P_2+(-1)^{(\norm{P_1}+\norm{c}-1)(\norm{P_2}-1)}P_2\cap_\gamma \mathrm L_{P_1}c\right)+\\
&\phantom{=}+(-1)^{(\norm{P_2}-1)}\left([P_1,P_2]\cap_\gamma Q+(-1)^{(\norm{c}-1)(\norm{P_1}+\norm{P_2}-1)}c\cap_\gamma [P_1,P_2]\right)+\\
&\phantom{=}+(-1)^{\norm{P_2}}\left(\mathrm L_\gamma(\mathrm m_{L,2}^{2,0,0}(P_1,P_2,c))-\mathrm m_{L,2}^{2,0,0}( \rm d_{\gamma,1}P_1,P_2,c)-(-1)^{\norm{P_1}}\mathrm m_{L,2}^{2,0,0}(P_1,\mathrm d_{\gamma,1} P_2,c)-\right.\\
&\phantom{=}\left.-(-1)^{\norm{P_1}+\norm{P_2}}\mathrm m_{L,2}^{2,0,0}(P_1,P_2,\rm L_\gamma c)\right)+(-1)^{\norm{P_1}}\left(\mathrm L_\gamma (\mathrm m_{L,2}^{2,0,0}(P_2,P_1,c))-\mathrm m_{L,2}^{2,0,0}(\rm d_{\gamma,1} P_2,P_1,c)-\right.\\
&\phantom{=}\left.-(-1)^{\norm{P_2}}\mathrm m_{L,2}^{2,0,0}(P_2,\mathrm d_{\gamma,1} P_1,c)-(-1)^{\norm{P_1}+\norm{P_2}}\mathrm m_{L,2}^{2,0,0}(P_2,P_1,\rm L_\gamma c)\right),
\end{aligned}
\]
for $P_i$, $i=1,2$, in $\underline{\mathrm{End}}(A)$ and $c$ in $A\underline{\otimes} A$.

We observe that there is a homotopy formula, similar to the first one we have written down, for the right action $\cap_\gamma$: this explains the appearance of many terms in the second homotopy formula.
We observe that the graded anti-commutators in the second homotopy formula sum up in the corresponding cohomology, whence the second homotopy formula restricts on cohomology to~\eqref{eq-T-mod}.
\end{proof}

\section{Configuration spaces and integral weights}\label{s-3}

In this Section we discuss in some details compactifications of configuration spaces of $i)$ points in the complex 
upper-half plane $\mathcal H$ and on the real axis $\mathbb R$, and $ii)$ points in the interior of the punctured 
unit disk $D^\times$ and on the unit circle $S^1$.

We will focus our attention on $\mathcal C_{2,0}\cong \mathcal D_{1,1}$ and on its boundary stratification: it will 
play a central r\^ole in the proof of both compatibilities with cup and cap products.
We will also take a better look at the compactified configuration space $\mathcal C_{2,1}\cong \mathcal D_{1,2}^+$: 
though it is not crucial in the forthcoming proofs, its boundary stratification leads to a better understanding of 
the homotopy formula for the compatibility between cap products, see e.g.~\cite{CR1}.

\subsection{Configuration spaces and their compactifications}\label{ss-3-1}

In this Subsection we recall compactifications of configuration spaces of points in the complex 
upper-half plane $\mathcal H$ and on the real axis $\mathbb R$, and of points in the interior 
of the punctured unit disk $D^\times$ and on the unit circle $S^1$. 

\subsubsection{Configuration spaces $C_{A,B}^+$ and $C_A$}\label{sss-3-1-1}

We consider a finite set $A$ and a finite (totally) ordered set $B$.
We define the open configuration space $C_{A,B}^+$ as 
$$
C_{A,B}^+:=\left\{(p,q)\in \mathcal H^A\times \mathbb R^B\,|\,p(a)\neq p(a')\textrm{ if }a\neq a'\,,\,
q(b)< q(b')\textrm{ if }b<b'\right\}/G_2,
$$
where $G_2$ is the semidirect product $\mathbb R^+\ltimes \mathbb R$, which acts diagonally on $\mathcal H^A\times \mathbb R^B$ via
\[
(\lambda,\mu)(p,q)=(\lambda p+\mu,\lambda q+\mu)\qquad(\lambda\in\mathbb R^+,\ \mu\in\mathbb R)\,.
\]
The action of the $2$-dimensional Lie group $G_2$ on such $n+m$-tuples is free, precisely when $2|A|+|B|-2\geq 0$: 
in this case, $C_{A,B}^+$ is a smooth real manifold of dimension $2|A|+|B|-2$. 

The configuration space $C_A$ is defined as
$$
C_A:=\left\{p\in\mathbb C^A\,|\,p(a)\neq p(a')\textrm{ if }a\neq a'\right\}/G_3,
$$
where $G_3$ is the semidirect product $\mathbb R^+\ltimes \mathbb C$, which acts diagonally on $\mathbb C^A$ via 
\[
(\lambda,\mu)p=\lambda p+\mu\qquad(\lambda\in\mathbb R^+,\ \mu\in\mathbb C)\,.
\] 
The action of $G_3$, which is a real Lie group of dimension $3$, is free precisely when $2|A|-3\geq 0$, 
in which case $C_A$ is a smooth real manifold of dimension $2|A|-3$. 

Finally, we observe that the spaces $C_{A,B}^+$ and $C_A$ are orientable, see e.g.~\cite{AMM} for a complete 
discussion of orientations of such configuration spaces. 

The configuration spaces $C_{A,B}^+$, resp.\ $C_A$, admit compactifications \`a la Fulton--MacPherson, 
obtained by successive real blow-ups: we will not discuss here the construction of their compactifications 
$\mathcal C_{A,B}^+$, $\mathcal C_A$, which are smooth manifolds with corners, referring to~\cite{K},~\cite{MT},
~\cite{CR} for more details, but we focus mainly on their stratification, in particular on the boundary strata 
of codimension $1$ of $\mathcal C_{A,B}^+$. 

Namely, the compactified configuration space $\mathcal C_{A,B}^+$ is a stratified space, and its boundary 
strata of codimension $1$ look like as follows:
\begin{enumerate}
\item[$i)$] there is a subset $A_1$ of $A$, resp.\ an ordered subset $B_1$ of successive elements of $B$, such that 
\begin{equation}\label{eq-upbound1}
\partial_{A_1,B_1}\mathcal C_{A,B}^+\cong
\mathcal C_{A_1,B_1}^+\times \mathcal C_{A\smallsetminus A_1,B\smallsetminus B_1\sqcup \{*\}}^+:
\end{equation}
intuitively, this corresponds to the situation, where points in $\mathcal H$, labelled by $A_1$, and successive points in $\mathbb R$ 
labelled by $B_1$, collapse to a single point labelled by $*$ in $\mathbb R$.
Obviously, we must have $2|A_1|+|B_1|-2\geq 0$ and $2(|A|-|A_1|)+(|B|-|B_1|+1)-2\geq 0$.
\item[$ii)$] there is a subset $A_1$ of $A$, such that
\begin{equation}\label{eq-upbound2}
\partial_{A_1}\mathcal C_{A,B}^+\cong \mathcal C_{A_1}\times \mathcal C_{A\smallsetminus A_1\sqcup \{*\},B}^+:
\end{equation}
this corresponds to the situation, where points in $\mathcal H$, labelled by $A_1$, collapse together to a single point $*$ in $\mathcal H$, labelled by $*$.
Again, we must have $2|A_1|-3\geq 0$ and $2(|A|-|A_1|+1)+|B|-2\geq 0$.
\end{enumerate}

\subsubsection{Configuration spaces $D_{A,B}^+$ and $D_A$}\label{sss-3-1-2}
We consider a finite set $A$ and a finite, cyclically ordered set $B$. 
We define the open configuration space $D_{A,B}^+$ as
$$
\left\{(p,q)\in(D^\times)^A\times (S^1)^B\,|\,p(a)\neq p(a'),\,a\neq a',\ q(b_1)< q(b_2)<\cdots<q(b_1),\,b_1< b_2<\cdots<b_1\right\}/S^1\,,
$$
where $D^\times$ denotes the punctured unit disk. Here the group $S^1$ acts on $D_{A,B}^+$ by rotations: 
the action is free, precisely when $2|A|+|B|-1\geq 0$, in which case $D_{A,B}^+$ is a smooth 
real manifold of dimension $2|A|+|B|-1$.

We also consider the configuration space
$$
D_A=\left\{p\in (\mathbb C^\times)^A\,|\,p(a)\neq p(a'),\,a\neq a'\right\}/\mathbb R^+\,,
$$
where $\mathbb R^+$ acts by rescaling.
It is obviously a smooth real manifold of dimension $2|A|-1$, when $2|A|-1\geq 0$. 

We also observe that, analogously to the configuration spaces $C_{A,B}^+$ and $C_A$ of the previous Subsection, 
$D_{A,B}^+$ and $D_A$ are orientable: it follows along the same lines of the discussion in~\cite{AMM}. Moreover, 
they also admit natural compactifications obtained by successive real blow-ups. 

The boundary strata of codimension $1$ of $\mathcal D_{A,B}^+$ are given in the following list:
\begin{enumerate}
\item[$i)$] there is a subset $A_1$ of $A$, such that
\begin{equation}\label{eq-diskb1}
\partial_{A_1,0}\mathcal D_{A,B}^+\cong \mathcal D_{A_1}\times \mathcal D_{A\smallsetminus A_1,B}^+.
\end{equation}
Intuitively, this corresponds to the situation, where points in $D^\times$ labelled by $A_1$ tend together to the origin.
Clearly, we must have $2|A_1|-1\geq 0$ and $2(|A|-|A_1|)+|B|-1\geq 0$.
\item[$ii)$] There is a subset $A_1$ of $A$, such that
\begin{equation}\label{eq-diskb2}
\partial_{A_1}\mathcal D_{A,B}^+\cong \mathcal C_{A_1}\times \mathcal D_{A\smallsetminus A_1\sqcup \{*\},B}^+.
\end{equation}
This corresponds to the situation, where points in $D^\times$ labelled by $A_1$ collapse together to a point in $D^\times$, labelled by $*$.
We must impose $2|A_1|-3\geq 0$ and $2(|A|-|A_1|+1)+|B|-1\geq 0$.
\item[$iii)$] Finally, there is a subset $A_1$ of $A$ and an ordered subset $B_1$ of successive elements of $B$, 
such that
\begin{equation}\label{eq-diskb3}
\partial_{A_1,B_1}\mathcal D_{A,B}^+\cong \mathcal C_{A_1,B_1}^+\times \mathcal D_{A\smallsetminus A_1,B\smallsetminus B_1\sqcup \{*\}}^+,
\end{equation}
which describes the situation, where points in $D^\times$ labelled by $A_1$ and successive points labelled by $B_1$ collapse together to a point in $S^1$, labelled by $*$.
We have to impose $2|A_1|+|B_1|-2\geq 0$ and $2(|A|-|A_1|)+(|B|-|B_1|+1)-1\geq 0$.
\end{enumerate}

\subsubsection{An identification}\label{sss-3-1-3}
Considering the special case $A=[n]$ and $B=[m]$, $m\geq 1$, we may use the action of $S^1$ to construct a 
section of $D_{A,B}^+$ by fixing the first point in $S^1$ to $1$.
This section is diffeomorphic, by means of the M\"obius transformation 
$$
\psi\,:\,\mathcal H\sqcup\mathbb R\,\longrightarrow\,D\sqcup S^1\smallsetminus\{\mathrm 1\}\,;
\,z\,\longmapsto\,\frac{z-\mathrm i}{z+\mathrm i}\,,
$$
where $D$ is the unit disk, to a smooth section of $C_{n+1,m-1}^+$, given by fixing e.g.\ the first point in 
the complex upper half-plane $\mathcal H$ to $\mathrm i$ by means of the action of $G_2$.

Then, the compactified configuration space $\mathcal D_{n,m}^+$ can be identified with $\mathcal C_{n+1,m-1}^+$, 
and we observe that the cyclic order of the $m$ points in $S^1$ translates naturally into an order of the $m-1$ points 
on the real axis $\mathbb R$. 

We point out that in certain situation it is better to use the compactified configuration spaces $\mathcal D_{n,m}^+$
instead of the equivalent $\mathcal C_{n+1,m-1}^+$, because $i)$ a cyclic order is visualized in an easier way on 
$S^1$ and $ii)$ we need two special points (the origin and the first point in $S^1$), which are also better visualized 
in the punctured disk $D$ with boundary. 

We further consider the manifold $D_n$, for $n\geq 1$ and notice the identification $D_n\cong C_{n+1}$: to be more 
precise, by means of complex translation, we may put e.g.\ the first point in $C_{n+1}$ at the origin, and using 
rescalings, one can put the remaining points in the punctured unit disk with boundary. 
Analogously as before, the compactification $\mathcal D_n$ of $D_n$ can be identified with $\mathcal C_{n+1}$.  

\medskip

More generally, this identification remains possible for arbitrary $A,B$, after the choice of distinguished elements $\bul\in A$ and $\circ\in B$. 
We consequently identify the codimension $1$ boundary strata of $\mathcal D_{A,B}^+$ with those of $\mathcal C_{A\sqcup\{\bullet\},B\smallsetminus\{\circ\}}^+$ (higher codimension can be worked out along the same lines very easily):  
\begin{itemize}
\item[$i)$] the situation \eqref{eq-diskb1} where points labelled by $A_1$ collapse 
together to the origin corresponds to the situation \eqref{eq-upbound2}, where points labelled by 
$A_1\sqcup\{\bullet\}$ collapse together to a single point in $\mathcal H$, which takes the r\^ole of the marked point $\bul$. 
\item[$ii)$] The situation \eqref{eq-diskb2}, where points labelled by $A_1$ collapse to a single point 
in $D$, corresponds to the situation \eqref{eq-upbound2}, where points labelled by $A_1$ 
collapse together to a single point in $\mathcal H$, which will not be the new marked point $\bul$. 
\item[$iii)$] The situation \eqref{eq-diskb3}, where points labelled by $A_1\sqcup B_1$, with $\circ\notin B_1$, 
collapse to a single point in $S^1$ corresponds to the situation \eqref{eq-upbound1}, where points labelled 
by $A_1\sqcup B_1$ collapse to a single point in $\mathbb{R}$, which will not be the new marked point $\circ$. 
\item[$iv)$] Finally, the situation \eqref{eq-diskb3}, where points labelled by $A_1\sqcup B_1$, with $\circ\in B_1$, 
collapse to a single point in $S^1$ corresponds to the situation \eqref{eq-upbound1}, where points labelled 
by the set $(A\backslash A_1\sqcup\{\bullet\})\sqcup (B\backslash B_1)$ collapse to a single point in $\mathbb{R}$, which will be the new marked point $\circ$. 
\end{itemize}

\subsection{Two remarkable compactified configuration spaces}\label{ss-3-2}
We describe two remarkable compactified configuration spaces: 
{\bf the eye} and the {\bf I-cube}. 

\subsubsection{The eye}\label{sss-3-2-1}
We now describe explicitly the compactified configuration space $\mathcal C_{2,0}$, known as {\bf Kontsevich's eye}. 
Here is a picture of it, with all boundary strata of codimension $1$, labelled by Greek letter, and codimension $2$, labelled by Latin letters, which we will describe shortly afterwards:
\begin{center}
\includegraphics[scale=0.37]{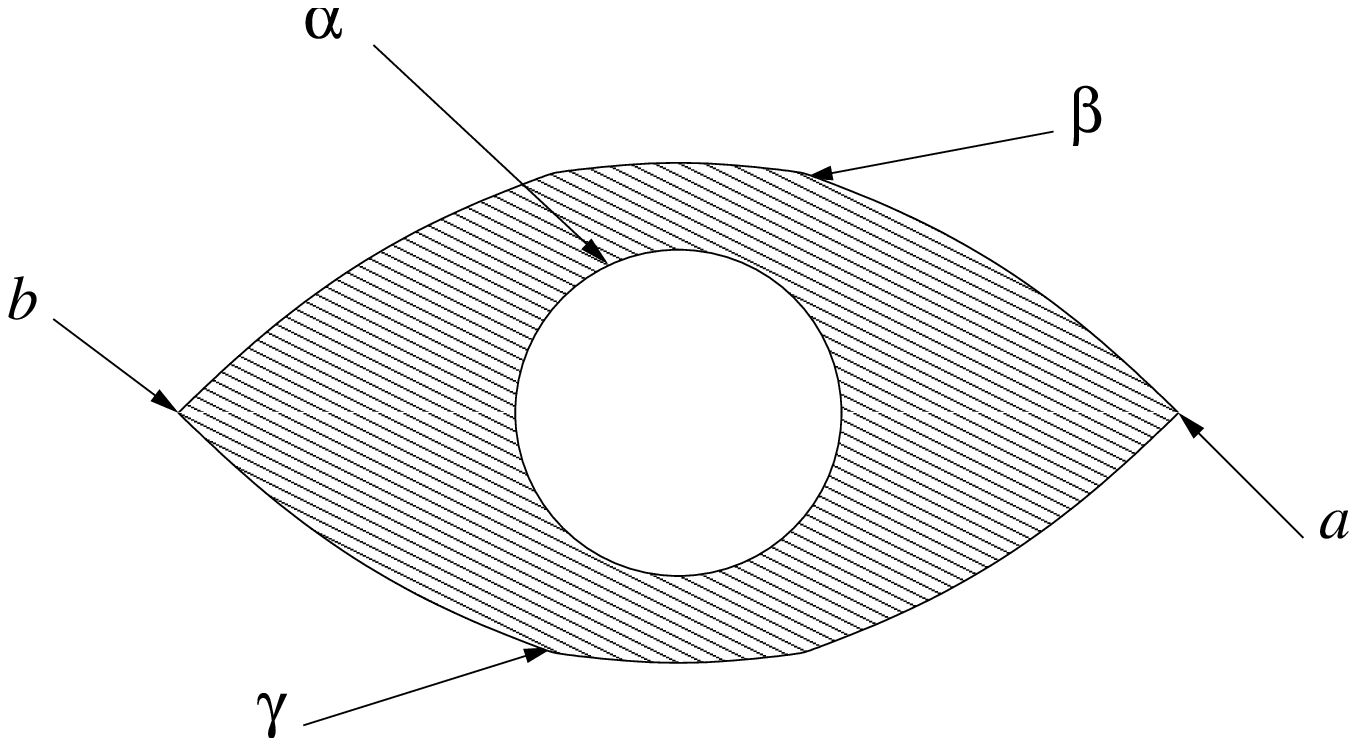} \\
\text{Figure 5 - Kontsevich's eye} \\
\end{center}
We first describe the boundary strata of codimension $1$.
\begin{itemize}
\item[$i)$] The stratum labelled by $\alpha$ corresponds to $\mathcal C_2=S^1$: intuitively, it describes the situation, where the two points collapse to a single point in $\mathcal H$; 
\item[$ii)$] the stratum labelled by $\beta$ corresponds to $\mathcal C_{1,1}\cong[0,1]$: it describes the situation. where the first point goes to the real axis; 
\item[$iii)$] the stratum labelled by $\gamma$ corresponds to $\mathcal C_{1,1}\cong[0,1]$: it describes the situation, where the second point goes to the real axis. 
\end{itemize}
As already observed, we have the identification $\mathcal C_{2,0}\cong\mathcal D_{1,1}$, and we can thus reinterpret its codimension $1$ boundary strata as follows:
\begin{itemize}
\item[$i)$] the stratum labelled by $\alpha$ corresponds to $\mathcal D_1=S^1$, which describes the situation, where the point in $D^\times$ goes to the origin; 
\item[$ii)$] the stratum labelled by $\beta$ corresponds to $\mathcal C_{1,1}\cong[0,1]$, which describes the situation, where the point in $D^\times$ and the point on $S^1$ collapse together in $S^1$; 
\item[$iii)$] the stratum labelled $\gamma$ corresponds to $\mathcal D^+_{0,2}\cong[0,1]$, which describes the situation, the point in $D^\times$ goes to $S^1$, without collapsing to the single point of $S^1$. 
\end{itemize}
We use the following pictorial representation for the latter situations, where we consider Kontsevich's eye as $\mathcal D_{1,1}$: 
\begin{center}
  \includegraphics[scale=0.30]{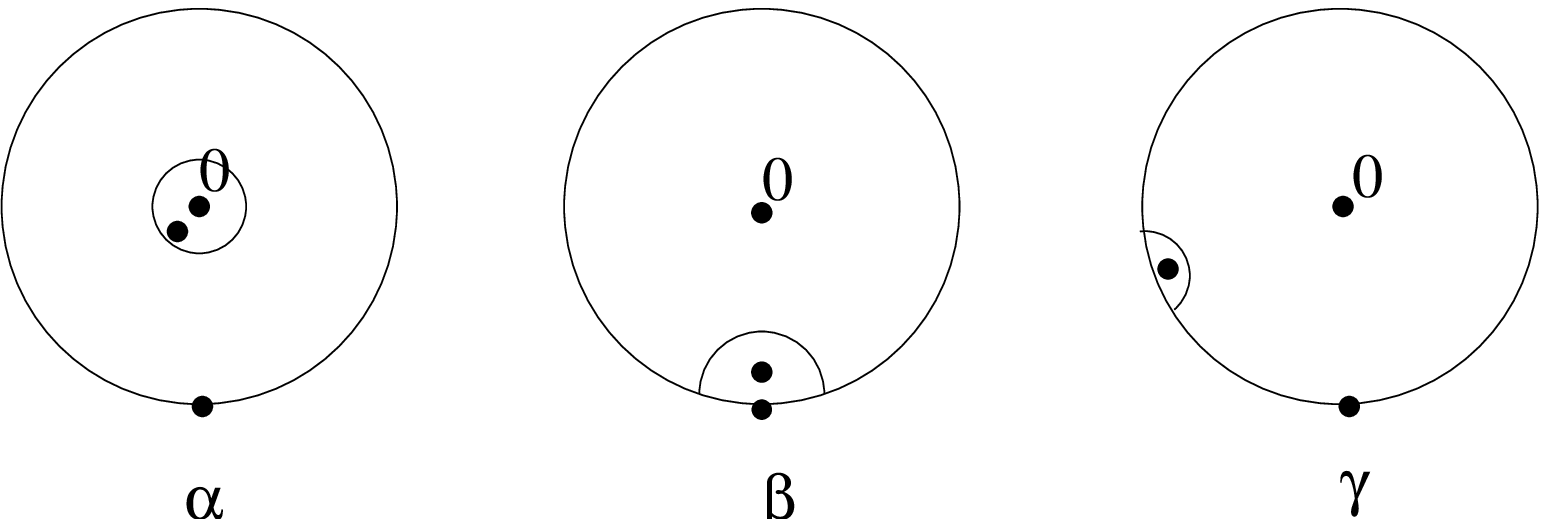} \\
\text{Figure 6 - Boundary strata of codimension $1$ of $\mathcal D_{1,1}$} \\
\end{center}
Finally, the two boundary strata of codimension $2$ are each one a copy of $\mathcal C_{0,2}^+=\{{\rm pt}\}$. 
The situations they describe can be depicted as follows: 
\begin{center}
  \includegraphics[scale=0.30]{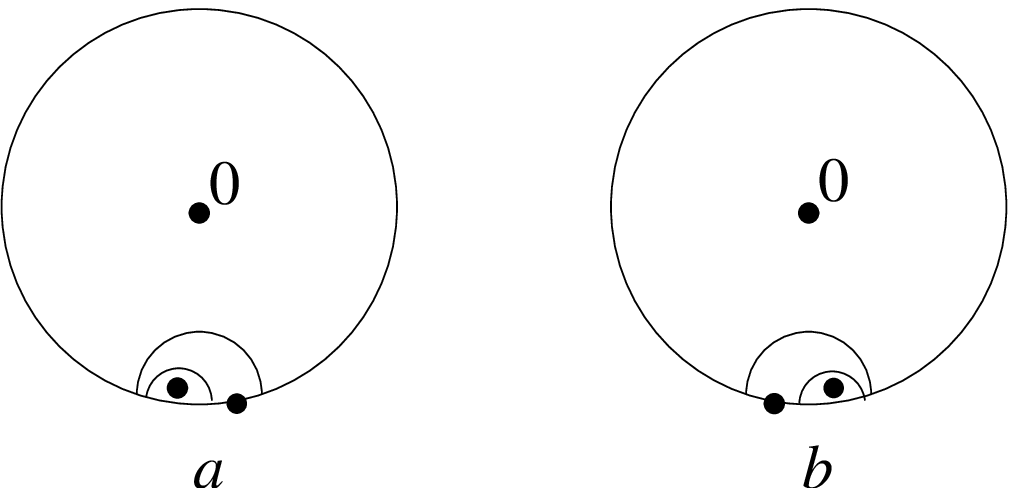} \\
\text{Figure 7 - Boundary strata of codimension $2$ of $\mathcal D_{1,1}$} \\
\end{center}

\subsubsection{The I-cube}\label{sss-3-2-2}
We now describe shortly the compactified configuration space $\mathcal C_{2,1}\cong\mathcal D_{1,2}^+$, 
which will be called the {\bf I-cube}: in particular, we are interested in its boundary strata of codimension $1$ and $2$. 
As in Subsubsection~\ref{sss-3-2-1}, we use Greek letters, resp.\ Latin letters, for 
labelling boundary strata of codimension $1$, resp.\ $2$. 

Pictorially, the I-cube looks like as follows:
\begin{center}
  \includegraphics[scale=0.38]{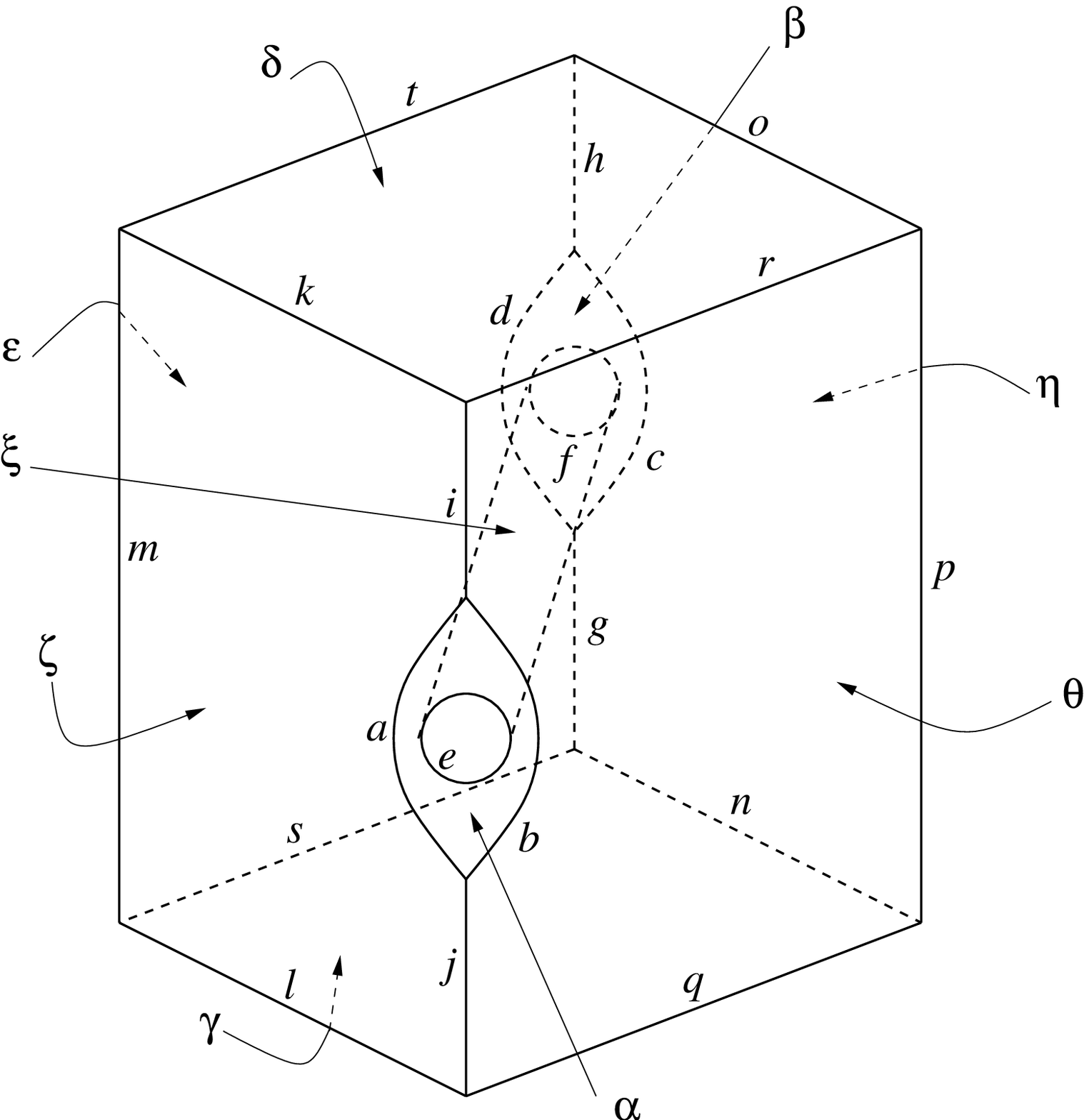} \\
\text{Figure 8 - The I-cube $\mathcal C_{2,1}\cong \mathcal D_{1,2}^+$} \\
\end{center}
Its boundary stratification consists of $9$ strata of codimension $1$, $20$ strata of codimension $2$ and $12$ 
strata of codimension $3$. 

\paragraph{{\bf Boundary strata of codimension $1$}}\label{p-2-2-2-1}
We illustrate explicitly the boundary strata of codimension $1$: again, before describing them mathematically, it is better to depict them:
\bigskip
\begin{center}
  \includegraphics[scale=0.28]{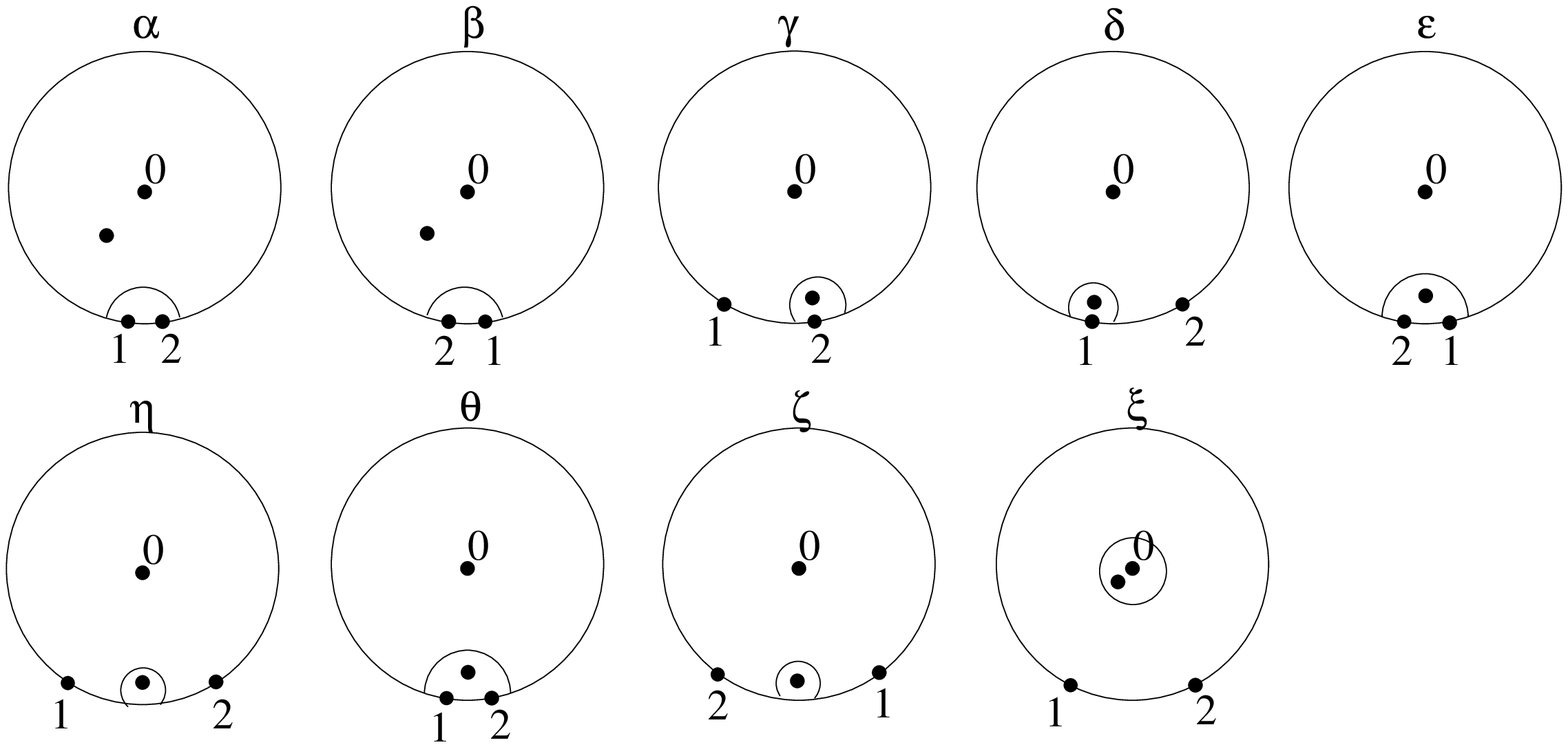} \\
\text{Figure 9 - Boundary strata of the I-cube of codimension $1$ } \\
\end{center}
\bigskip
The strata labelled by $\alpha$ and $\beta$ are both described by $\mathcal C_{0,2}^+\times \mathcal D_{1,1}$, depending on the cyclic order 
of the two points in $S^1$: since $\mathcal C_{0,2}^+$ is $0$-dimensional, the strata $\alpha$ and $\beta$ are two copies of Kontsevich's eye $\mathcal D_{1,1}$.

As for the strata labelled by $\gamma$ and $\delta$, they are both described by $\mathcal C_{1,1}\times\mathcal D_{0,2}^+$: both $\mathcal C_{1,1}$ 
and $\mathcal D_{0,2}^+$ correspond to closed intervals, whence $\gamma$ and $\delta$ are topologically two squares.

The strata labelled by $\varepsilon$ and $\theta$ correspond both to $\mathcal C_{1,2}^+\times\mathcal D_{0,1}$, depending on the cyclic order of 
points in $S^1$: recalling the results of Subsection~\ref{ss-2-1}, $\varepsilon$ and $\theta$ are topologically two copies of the hexagon (this 
will be also clearer after the description of the boundary strata of codimension $2$ of the I-cube).

On the other hand, the strata labelled by $\eta$ and $\zeta$ are both described by $\mathcal C_{1,0}\times\mathcal D_{0,3}^+$, depending on the 
cyclic order on $S^1$: again, since $\mathcal C_{1,0}$ is $0$-dimensional, an inspection of $\mathcal D_{0,3}^+$ shows that $\eta$ and $\zeta$ 
are topologically two copies of the hexagon (again, we deserve a careful explanation, when dealing with boundary strata of codimension $2$ of the I-cube). 

Finally, the stratum labelled by $\xi$ corresponds to $\mathcal D_1\times \mathcal D_{0,2}^+$: since $\mathcal D_1=S^1$ and $\mathcal D_{0,2}^+$ 
is a closed interval, topologically $\xi$ looks like a cylinder.

The above picture describes the boundary strata of codimension $1$ of $\mathcal D_{1,2}^+$: using the prescriptions of Subsubsection~\ref{sss-3-1-3}, it is then easy to identify these boundary strata with the corresponding boundary strata of codimension $1$ of $\mathcal C_{2,1}$.

\paragraph{{\bf Boundary strata of codimension $2$}}\label{p-2-2-2-2}
We discuss now some relevant boundary strata of the I-cube of codimension $2$: we first illustrate all of them pictorially as follows, referring 
to the picture of the I-cube for the notations:
\bigskip
\begin{center}
  \includegraphics[scale=0.28]{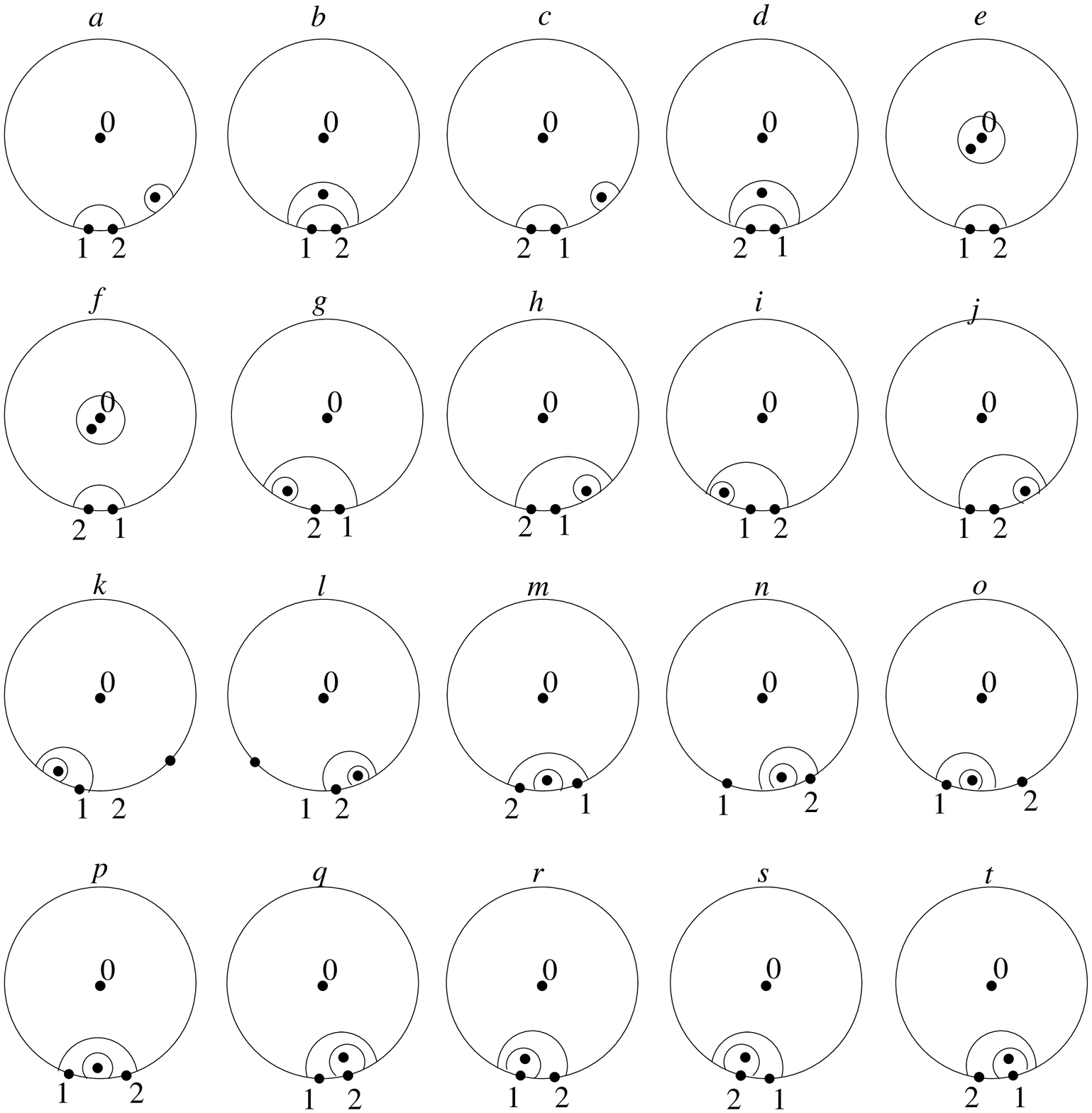} \\
\text{Figure 10 - Boundary strata of the I-cube of codimension $2$ } \\
\end{center}
\bigskip
For our purposes, we need only describe explicitly the boundary strata labelled by $e$, $f$, $h$, $j$, $p$, $q$ and $o$: in fact, these describe 
certain boundary components of a particular imbedding of the plane square into the I-cube, which will be useful later on.

The strata labelled by $e$ and $f$ are described as $\mathcal D_1\times \mathcal C_{0,2}^+\times\mathcal D_{0,1}$: since $ \mathcal C_{0,2}^+$ 
and $\mathcal D_{0,1}$ are both $0$-dimensional, while $\mathcal D_1=S^1$, they can be identified with the pupils of the two brave new eyes. 

The strata labelled by $h$, $j$ and $p$ correspond all to $\mathcal C_{1,0}\times \mathcal C_{0,3}^+\times \mathcal D_{0,1}$, hence they correspond 
all to a closed interval, the only difference depending on the cyclic order on $S^1$ and on the corresponding order on $\mathbb R$.

As for the stratum labelled by $o$, it is described as $\mathcal C_{1,0}\times\mathcal C_{0,2}^+\times\mathcal D_{0,2}^+$, 
which is also topologically, by previous arguments, a closed interval.

Finally, the stratum labelled by $q$ is $\mathcal C_{1,1}\times\mathcal C_{0,2}^+\times\mathcal D_{0,1}$: once again, 
it is topologically a closed interval.

It is left as an exercise to identify the boundary strata of $\mathcal D_{1,2}^+$ with the boundary strata 
of $\mathcal C_{2,1}$, according to the prescriptions of Subsubsection~\ref{sss-3-1-3}.

\subsection{Integral weights associated to graphs}\label{ss-3-3}

In this Subsection we recall Kontsevich's, resp.~Shoikhet's, angle forms and the corresponding weights, 
resp.~modified weights, associated to graphs. 

\subsubsection{Angle forms}\label{sss-3-3-1}

We first need to specify a smooth $1$-form on the configuration space $\mathcal C_{2,0}$.
For any two distinct points $p,q$ in $\mathcal H\sqcup\mathbb R$, we define
\[
\varphi(p,q)=\frac1{2\pi}\mathrm{arg}\!\left(\frac{q-p}{q-\overline p}\right)\,.
\] 
The real number $\varphi(p,q)$ represents the (normalized) angle from the geodesic from $p$ to the point $\infty$ on the 
positive imaginary axis to the geodesic between $p$ and $q$ w.r.t.\ the hyperbolic metric of 
$\mathcal H\sqcup\mathbb R$, measured in counterclockwise direction. It is therefore defined up to a constant, and thus 
$\omega:={\rm d}\varphi$ is a well-defined $1$-form. 
\begin{Lem}\label{l-angle}
The $1$-form $\omega$ extends to a smooth $1$-form on Kontsevich's eye $\mathcal C_{2,0}$, with the 
following properties: 
\begin{itemize}
\item[$a)$] the restriction of $\omega$ to the boundary stratum $\mathcal C_2=S^1$ equals the total derivative of 
the (normalized) angle measured in counterclockwise direction from the positive imaginary axis;
\item[$b)$] the restriction of $\omega$ to $\mathcal C_{1,1}$, where the first point goes to the real axis, 
vanishes. \hfill\qed
\end{itemize}
\end{Lem}

\medskip

We then define a smooth $1$-form on the configuration space $\mathcal C_{3,0}$. 
For any three pairwise distinct points $p,q,r$ in $\mathcal H\sqcup\mathbb R$, we define the 
modified Kontsevich angle function
\begin{equation}\label{eq-Dprop}
\varphi_D(p,q,r)=\varphi(q,r)-\varphi(q,p)\,,
\end{equation}
and set $\omega_D=\mathrm d\varphi_D$.
\begin{Lem}\label{l-angleD}
The $1$-form $\omega_D$ extends to a smooth $1$-form on $\mathcal C_{3,0}$, with the 
following properties: 
\begin{itemize}
\item[$a)$] its restriction to $\mathcal C_{2,1}$, when the second point approaches the real axis, vanishes; 
\item[$b)$] its restriction to $\mathcal C_{2,0}\times\mathcal C_{1,1}$, 
when the first and second point collapse together to the real axis, 
equals $-\pi_1^*\omega$;
\item[$c)$] its restriction to $\mathcal C_2\times\mathcal C_{2,0}$, 
when the first and second point collapse together in the upper half-plane, 
equals $\pi_2^*\omega-\pi_1^*\omega$;
\item[$d)$] its restriction to $\mathcal C_{2,0}\times\mathcal C_{1,1}$ (resp.~$\mathcal C_2\times\mathcal C_{2,0}$), 
when the first and third point collapse together to the real axis (resp.~in the upper half-plane), 
vanishes;
\item[$e)$] its restriction to $\mathcal C_{2,0}\times\mathcal C_{1,1}$, 
when the second and third point collapse together to the real axis, 
equals $\pi_1^*\omega$. 
\item[$f)$] its restriction to $\mathcal C_2\times\mathcal C_{2,0}$, 
when the second and third point collapse together in the upper half-plane, 
equals $\pi_1^*\omega-\pi_2^*\omega$. 
\hfill\qed
\end{itemize}
\end{Lem}
Using the prescriptions of Subsubsection~\ref{sss-3-1-3}, it is easy to identify the boundary strata of $\mathcal C_{2,0}$ in Lemma~\ref{l-angle} and the boundary strata of $\mathcal C_{3,0}$ in Lemma~\ref{l-angleD} with the corresponding boundary strata of $\mathcal D_{1,1}$ and $\mathcal D_{2,1}$.

In particular, here is a useful pictorial description of the angle function~\eqref{eq-Dprop}:
\bigskip
\begin{center}
\resizebox{0.23 \textwidth}{!}{\input{modKangle.pstex_t}} \\
\text{Figure 11 - The modified Kontsevich's angle functions $\varphi_D$} \\
\end{center}
\bigskip

\subsubsection{Integral weights associated to graphs}\label{sss-3-3-2}

We consider, for given positive integers $n$ and $m$, directed graphs $\Gamma$ with $m+n$ vertices 
labelled by the set $\mathcal V(\Gamma)=\{1,\dots,n,\overline{1},\dots,\overline{m}\}$. 
Here, ``directed'' means that each edge of $\Gamma$ carries an orientation. 
Additionally, the graphs we consider are required to have no loop 
(a loop is an edge beginning and ending at the same vertex). 
To any edge $e=(i,j)$ of such a directed graph $\Gamma$, we associate 
the smooth $1$-form $\omega_e:=\pi_e^*\omega$ on $\mathcal C_{n,m}^+$, where 
$\pi_e\,:\,\mathcal C_{n,m}^+\to\mathcal{C}_{2,0}$ is the smooth map given by 
$$
\big[(z_1,\dots,z_n,z_{\overline{1}},\dots,z_{\overline{m}})\big]\,\longmapsto\,\big[(z_i,z_j)\big]\,.
$$

Then, to any directed graph $\Gamma$ without loop, and denoting by $\mathcal E(\Gamma)$ the set of its edges, 
we associate a differential form 
\begin{equation}\label{eq-weight-f}
\omega_\Gamma:=\bigwedge_{e\in \mathcal E(\Gamma)}\omega_e
\end{equation}
on the (compactified) configuration space $\mathcal{C}_{n,m}^+$. 
\begin{Rem}
We observe that, {\em a priori}, it is necessary to choose an ordering of the edges of $\Gamma$ since $\omega_\Gamma$ 
is a product of $1$-forms: two different orderings of the edges of $\Gamma$ simply differ by a sign. 
This sign ambiguity precisely coincide (and thus cancel) with the one appearing in the definition of $B_\Gamma$ 
(see next Section). 
\end{Rem}
We recall that $\mathcal C_{n,m}^+$ is orientable and its orientation specifies an orientation for any boundary stratum thereof. 
\begin{Def}\label{d-weight}
The {\em weight} $W_\Gamma$ of the directed graph $\Gamma$ is 
\begin{equation}\label{eq-weight}
W_\Gamma:=\int_{\mathcal C_{n,m}^+} \omega_\Gamma\,.
\end{equation}
\end{Def} 
Observe that the weight (\ref{eq-weight}) truly exists since it is an integral of a smooth differential form 
over a compact manifold with corners. 

\medskip

In the same way, we define a modified weight associated to a graph without loop $\Gamma$ with 
$m+n+1$ vertices labelled by $\mathcal{V}(\Gamma):=\{0,\dots,n,\overline{1},\dots,\overline{m}\}$. 
To any edge $e=(i,j)\in\mathcal E(\Gamma)$, we associate a smooth $1$-form $\omega_{D,e}$ on $\mathcal D_{n,m}^+$ by the following rules: 
\begin{itemize}
\item if neither $i$ nor $j$ lies in $\{0,\overline{1}\}$, then $\omega_{D,e}:=\pi_{(0,i,j)}^*\omega_D$, where 
\begin{eqnarray*}
\pi_{(0,i,j)}\,:\,\mathcal D_{n,m}^+\cong\mathcal C_{n+1,m-1}^+ & \longrightarrow & \mathcal{C}_{3,0} \\
\big[(z_0,\dots,z_n,z_{\overline{2}},\dots,z_{\overline{m}})\big]  & \longmapsto & \big[(z_0,z_i,z_j)\big]\,;
\end{eqnarray*}
\item if $i=0$ and $j\neq\overline{1}$, then $\omega_{D,e}:=\pi_{(i,j)}^*\omega$, where 
\begin{eqnarray*}
\pi_{(i,j)}\,:\,\mathcal D_{n,m}^+\cong\mathcal C_{n+1,m-1}^+ & \longrightarrow & \mathcal{C}_{2,0} \\
\big[(z_0,\dots,z_n,z_{\overline{2}},\dots,z_{\overline{m}})\big] & \longmapsto & \big[(z_i,z_j)\big]\,;
\end{eqnarray*}
\item if $j=\overline{1}$ and $i\neq0$, then $\omega_{D,e}:={\rm p}_{(i,j)}^*\omega$, where
\begin{eqnarray*}
{\rm p}_{(i,j)}\,:\,\mathcal D_{n,m}^+ & \longrightarrow & \mathcal D_{1,1}\cong\mathcal{C}_{2,0} \\
\big[(z_1,\dots,z_n,z_{\overline{1}},\dots,z_{\overline{m}})\big] & \longmapsto & \big[(z_i,z_j)\big]\,;
\end{eqnarray*}
\item if $i=\overline{1}$ or $j=0$ or $(i,j)=(0,\overline{1})$, then $\omega_{D,e}=0$. 
\end{itemize}
Then, as above, 
\begin{equation}\label{eq-weight-f-mod}
\omega_{D,\Gamma}:=\bigwedge_{e\in \mathcal E(\Gamma)}\omega_{D,e}
\end{equation}
defines a differential form on $\mathcal D_{n,m}^+$. 
\begin{Def}\label{d-weight-mod}
The {\em modified weight} $W_{D,\Gamma}$ of the directed graph $\Gamma$ is 
\begin{equation}\label{eq-weight-mod}
W_{D,\Gamma}:=\int_{\mathcal D_{n,m}^+} \omega_{D,\Gamma}\,.
\end{equation}
\end{Def} 

\subsection{Explicit formul\ae\ for Kontsevich's and Tsygan's formality morphisms}\label{ss-3-4}

We quickly review the construction of $i)$ Kontsevich's $L_\infty$-quasi-isomorphism 
$\mathcal U$, and $ii)$ Shoikhet's $L_\infty$-quasi-isomorphism $\mathcal S$ respectively.

\subsubsection{The $L_\infty$-quasi-isomorphism $\mathcal U$}\label{sss-3-4-1}

For any pair of non-negative integers $(n,m)$, a {\bf K-admissible graph} $\Gamma$ 
of type $(n,m)$ is by definition a directed graph without loops and with labels obeying the following requirements:
\begin{enumerate}
\item[$i)$] the set of vertices $\mathcal{V}(\Gamma)$ is given by $\{1,\dots,n,\overline{1},\dots,\overline{m}\}$; 
vertices labelled by $\{1,\dots,n\}$, resp.\ $\{\overline{1},\dots,\overline{m}\}$, are called vertices 
of the first, resp.\ second, type; 
\item[ii)] every edge in $\mathcal{E}(\Gamma)$ starts at some vertex of the first type and there is at most one 
edge between any two distinct vertices of $\Gamma$. 
\end{enumerate}
For a given vertex $v$ of $\Gamma$, we denote by $\mathrm{star}(v)$ the subset of $\mathcal{E}(\Gamma)$ of 
edges starting at $v$: then, we assume that, for any vertex of the first type $v$ of $\Gamma$, the elements of 
$\mathrm{star}(v)$ are labelled as $(e_v^1,\dots,e_v^{|\mathrm{star}(v)|})$. 
By definition, the {\bf valence} of a vertex $v$ is the cardinality of the star of $v$.
The set of K-admissible graphs of type $(n,m)$ is denoted by $\mathcal G^K_{n,m}$.
\begin{Rem}
In the following we will use integral weights associated to graphs introduced in the previous Section. 
We can restrict ourselves safely to K-admissible graphs such that $|\mathcal{E}(\Gamma)|=2n+m-2$ as one 
can easily see that the weights vanish in any other situation. 
\end{Rem}
Finally, we define the $n$-th structure map $\mathcal U_n$ of Kontsevich's $L_\infty$-quasi-isomorphism by
\begin{equation}\label{eq-taylorU}
\mathcal U_n:=\sum_{m\geq 0}\sum_{\Gamma\in\mathcal G^K_{n,m}}W_\Gamma U_\Gamma\,:\,
\bigwedge^n T_\mathrm{poly}^\bullet(V)\to \mathcal D_\mathrm{poly}^\bullet(V)[1-n]\,,
\end{equation}
where $U_\Gamma(\alpha_1,\dots,\alpha_n)$ is a $m$-polydifferential operator naturally associated to 
the graph $\Gamma$ and polyvector fields $\alpha_1,\dots,\alpha_n$, as defined in \cite{K} (see also \cite[Appendix A.8]{CF}). 
\begin{Thm}[Kontsevich]\label{t-kont}
The Taylor components (\ref{eq-taylorU}) combine to an $L_\infty$-quasi-isomorphism
\[
\mathcal U:T_\mathrm{poly}^\bullet(V)\to \mathcal D_\mathrm{poly}^\bullet(V)
\]
of $L_\infty$-algebras, whose first order Taylor component reduces to the Hochschild--Kostant--Rosenberg quasi-isomorphism in cohomology.
\end{Thm}
The complete proof of Theorem~\ref{t-kont} is given in~\cite{K}: the main argument of the proof relies on a 
clever use of Stokes' Theorem to derive quadratic identities for the weights (\ref{eq-weight-f}) of (\ref{eq-taylorU}), 
which in turn imply the quadratic identities for (\ref{eq-taylorU}), corresponding to the fact that $\mathcal U$ is an $L_\infty$-morphism. 

\subsubsection{The $L_\infty$-quasi-isomorphism $\mathcal S$}\label{sss-3-4-2}

The construction of $\mathcal S$ is similar, in principle, to the construction sketched in the previous Subsection, 
but presents certain subtleties, which we need to discuss also for later purposes. 

An {\bf S-admissible graph} of type $(n,m)$ is a directed labelled graph $\Gamma$ without loops and such that: 
\begin{enumerate}
\item[$i)$] the set of vertices $\mathcal{V}(\Gamma)$ is given by $\{0,\dots,n,\overline{1},\dots,\overline{m}\}$; 
vertices labelled by $\{1,\dots,n\}$, resp.\ $\{\overline{1},\dots,\overline{m}\}$, are called vertices 
of the first, resp.\ second, type; 
\item[ii)] every edge in $\mathcal{E}(\Gamma)$ starts at some vertex of the first type and there is at most one 
edge between any two distinct vertices of $\Gamma$; 
\item[iii)] there is no edge ending at the special vertex $0$. 
\end{enumerate}
The set of S-admissible graphs of type $(n,m)$ is denoted by $\mathcal G^S_{n,m}$.
\begin{Rem}
As above we can restrict ourselves safely to S-admissible graphs such that $|\mathcal{E}(\Gamma)|=2n+m-1$, 
as the modified weights vanish in any other situation. 
\end{Rem}

We now consider an S-admissible graph in $\mathcal G^S_{n,m}$, such that $|\mathrm{star}(0)|=l$. 
To $n$ polyvector fields $\{\gamma_1,\dots,\gamma_n\}$ on $V$, such that $|\mathrm{star}(k)|=|\gamma_k|+1$, 
$k=1,\dots,n$, and to a Hochschild chain $c=(a_0|a_1|\cdots|a_{m-1})$ of degree $-m+1$, we associate an $l$-form 
on $V$ (whose actual degree is $-l$, following the grading in~\cite{Sh}) defined via 
\begin{equation}\label{eq-diffS}
\langle\alpha,S_{\Gamma}(\gamma_1,\dots,\gamma_n;c)\rangle:=U_\Gamma(\alpha,\gamma_1,\dots,\gamma_n)(a_0,\dots,a_{m-1})\,.
\end{equation}
The graded vector space $C_{-\bullet}^{\rm poly}(V)$ with Hochschild differential and DGLA action $\mathrm L$ over 
the DGLA $D_\mathrm{poly}(V)$ has a structure of a DGM: if we compose $\mathrm L$ with the L$_\infty$-quasi-isomorphism 
$\mathcal U$, $C_{-\bullet}^{\rm poly}(V)$ becomes an $L_\infty$-module over the DGLA $T_\mathrm{poly}(V)$.
The $n$-th Taylor component $\mathcal S_n$ of $\mathcal S$ from the $L_\infty$-module $C_{-\bullet}^{\rm poly}(V)$ 
to the $L_\infty$-module $\Omega^{-\bullet}(V)$ (actually, this is a true DGM, with trivial differential and action 
$\mathrm L$ given by the Lie derivative w.r.t.\ polyvector fields) over $T_\mathrm{poly}^\bullet(V)$ is given by
\begin{equation}\label{eq-taylorS}
\mathcal S_n:=\sum_{m\geq 1}\sum_{\Gamma\in\mathcal G^S_{n,m}}W_{D,\Gamma}S_{\Gamma}\,:\,\bigwedge^n T_\mathrm{poly}^\bullet(V)\otimes C_{-\bullet}^{\rm poly}(V)\to \Omega^{-\bullet}(V)[-n]\,.
\end{equation}
\begin{Thm}[Shoikhet]\label{t-shoikhet}
The Taylor components (\ref{eq-taylorS}) combine to an $L_\infty$-quasi-isomorphism 
\[
\mathcal S:C_{-\bullet}^{\rm poly}(V)\to \Omega^{-\bul}(V)
\]
of $L_\infty$-modules over $T_\bul(V)$, whose $0$-th order Taylor component reduces to the Hochschild--Kostant--Rosenberg quasi-isomorphism in homology.
\end{Thm}
We refer to~\cite{Sh} for a complete proof of Theorem~\ref{t-shoikhet}: the proof of the quadratic identities satisfied by the weights \eqref{eq-weight-mod} of 
(\ref{eq-taylorS}) can be found in~\cite{Sh}, and relies again on a clever use of Stokes' Theorem.

\section{The compatibility between cup products}\label{s-4}
We borrow the notation from Sections~\ref{s-1},~\ref{s-2} and~\ref{s-3}: in particular, 
$(\mathfrak m,{\rm d}_{\mathfrak m})$ is as in the introduction, 
and accordingly we consider the twisted DGLAs $T_{\mathrm{poly}}^\mathfrak m(V)$, 
$\mathcal D_{\mathrm{poly}}^\mathfrak m(V)$ with corresponding new gradings, products, differentials etc.

We further consider a general MCE $\gamma$ in $T_{\mathrm{poly}}^\mathfrak n(V)$ 
\begin{equation}\label{eq-MCel}
\gamma=\gamma_{-1}+\gamma_0+\gamma_1+\gamma_2+\cdots,
\end{equation}
where the suffix refers to the polyvector degree, which satisfies the MC equation
\[
\mathrm d_\mathfrak m\gamma+\frac{1}2[\gamma,\gamma]=0.
\]
We denote by $\mathcal U(\gamma)$ its image w.r.t.\ the (extension of the) $L_\infty$-quasi-isomorphism 
$\mathcal U$ of Theorem~\ref{t-kont}, Subsubsection~\ref{sss-3-4-1}, i.e.\
\[
\mathcal U(\gamma)=\sum_{n\geq 1}\frac{1}{n!}\mathcal U_n(\underset{n}{\underbrace{\gamma,\dots,\gamma}}).
\]
It is a MCE in $D_\mathrm{poly}^\mathfrak n(V)$, with infinitely many components of different polydifferential 
operator degree.

\subsection{The homotopy argument for the cup product}\label{ss-4-2}

We consider Kontsevich's eye $\mathcal C_{2,0}$, and a smooth curve $\ell$ therein, with starting point $\ell(0)$ 
on the pupil $\mathcal C_2$ and final point $\ell(1)$ in any one of the boundary strata of codimension $2$, 
$\mathcal C_{0,2}^+$, and such that $\ell(t)$ in $C_{2,0}$, for $t$ in $(0,1)$, e.g.
\bigskip
\begin{center}
\resizebox{0.31\textwidth}{!}{\input{K_curve.pstex_t}}\\
\text{Figure 12 - The curve $\ell$ in $\mathcal C_{2,0}$} \\
\end{center}
\bigskip 
More generally, for any pair of non-negative integers $(n,m)$, such that $n\geq 2$ 
(hence, automatically, $2n+m-2\geq 0$), we consider the subset $\mathcal Z_{n,m}^+$ of $\mathcal C_{n,m}^+$, 
which consists of those configurations, whose projection onto $\mathcal C_{2,0}$ through $\pi_{1,2}$ is in $\ell$.

Subsets of the form $\mathcal Z_{n,m}^+$ were introduced in~\cite{K}, and they were analyzed more carefully in~\cite{MT}: they are actually submanifolds with corners of $\mathcal C_{n,m}^+$ of codimension $1$, and they inherit an orientation from the orientation of $\ell$ and of the spaces $\mathcal C_{n,m}^+$ themselves, as shown in~\cite{MT}.

Another important feature of $\mathcal Z_{n,m}^+$ is the characterization of its boundary: for our purposes, we are interested only in its boundary strata of codimension $1$, which are of the following type:
\begin{enumerate}
\item[$i)$] Configurations of $\mathcal C_{n,m}^+$, such that their projection onto $\mathcal C_{2,0}$ is $\ell(0)$: this boundary stratum is denoted by $\mathcal Z_{n,m,0}^+$: explicitly, a general component $Z$ of $\mathcal Z_{n,m,0}^+$ splits as
\[
Z\cong \mathcal C_A^0\times \mathcal C_{n-|A|+1,m},
\]
where $A$ is a subset of $[n]$ with $2\leq |A|\leq n$, which contains the points labelled by $1$ and $2$, and $\mathcal C_A^0$ denotes the subset of $\mathcal C_A$, such that the projection onto $\mathcal C_2\cong S^1$ (corresponding to the points $1$ and $2$) is a fixed point of $S^1$.
\item[$ii)$] Configurations of $\mathcal C_{n,m}^+$, such that their projection onto $\mathcal C_{2,0}$ is $\ell(1)$: this boundary stratum is denoted by $\mathcal Z_{n,m,1}^+$.
A general component $Z$ of $\mathcal Z_{n,m,1}^+$ splits as
\[
Z\cong \mathcal C_{A_1,B_1}^+\times \mathcal C_{A_2,B_2}^+\times \mathcal C_{A_3,B_3}^+,
\]
where $1\leq |A_1|\leq n$, $0\leq |B_1|\leq m$, $1\leq |A_2|\leq n$, $0\leq |B_2|\leq m$, $0\leq |A_3|\leq n$, $2\leq |B_3|\leq m$, $1\in A_1$ and $2\in A_2$. 
\item[$iii)$] Non-trivial intersections of boundary strata of codimension $1$ of $\mathcal C_{n,m}^+$ with the interior $\overset{\circ}{\mathcal Z}_{n,m}^+$ of $\mathcal Z_{n,m}^+$: this we denote by $Z_{n,m}^+$.
We observe that, in this case, the first and second point of $\mathcal C_{n,m}^+$ are distinct and lie on the curve $\ell$.
The explicit form of a general component $Z$ of $Z_{n,m}^+$ will be described explicitly later on.
\end{enumerate}
Pictorially, configurations of points in the boundary strata of $\mathcal Z_{n,m}^+$ of type $i)$ and $ii)$ look like as follows: 
\bigskip
\begin{center}
\resizebox{0.95 \textwidth}{!}{\input{K-eye-h_i+ii.pstex_t}}\\
\text{Figure 13 - Typical configurations in the boundary strata of $\mathcal Z_{n,m}^+$ of type $i)$ and $ii)$} \\
\end{center}

\bigskip

For a MCE $\gamma$ as in (\ref{eq-MCel}) and any two $\mathfrak m$-valued polyvector fields $\alpha$, $\beta$ on $V$, 
we will construct a bilinear map $\mathcal H_\gamma^{\mathcal U}$ from 
$T_\mathrm{poly}^\mathfrak m(V)\otimes T_\mathrm{poly}^\mathfrak m(V)$ 
to $\mathcal D_\mathrm{poly}^\mathfrak m(V)$ such that the following identity holds true,
\begin{equation}\label{eq-homotopy-U}
\begin{aligned}
\mathcal U_\gamma(\alpha\cup\beta)-\mathcal U_\gamma(\alpha)\cup\mathcal U_\gamma(\beta)&=\mathcal H_\gamma^\mathcal U((\mathrm d_\mathfrak m\alpha+[\gamma,\alpha],\beta)+(-1)^{|\alpha|}\mathcal H_\gamma^\mathcal U(\alpha,\mathrm d_\mathfrak m\beta+[\gamma,\beta])+\\
&\phantom{=}+(\mathrm d_\mathfrak m+\mathrm d_\mathrm H+\mathrm L_{\mathcal U(\gamma)})\mathcal H_\gamma^\mathcal U(\alpha,\beta)\,.
\end{aligned}
\end{equation}

First of all, for a non-negative integer $m$, we define 
\begin{equation*}
\mathcal H_\gamma^{\mathcal U,m}(\alpha,\beta)=\sum_{n\geq 0}\frac{1}{n!}\sum_{\Gamma\in\mathcal G^K_{n+2,m}}\overset{\circ}{W}_\Gamma U_\Gamma(\alpha,\beta,\underset{n}{\underbrace{\gamma,\dots,\gamma}}),
\end{equation*}
where now the weight $\overset{\circ}{W}_\Gamma$ of an admissible graph $\Gamma$ in $\mathcal G^K_{n+2,m}$ is 
\[
\overset{\circ}W_\Gamma=\int_{\overset{\circ}{\mathcal Z}^+_{n+2,m}}\omega_\Gamma,
\]
with the same notations as above.
Finally, we set 
\[
\mathcal H_\gamma^\mathcal U(\alpha,\beta)=\sum_{m\geq 0}\mathcal H_\gamma^{\mathcal U,m}(\alpha,\beta).
\]
We want to reinterpret (\ref{eq-homotopy-U}) in terms of equalities between weights: we observe that the 
$\mathfrak m$-valued polyvector fields $\alpha$ and $\beta$ are put at the first and second vertex of the 
first type of any graph $\Gamma$ appearing in the morphisms $\mathcal U_n$. 

We first observe that both expressions on the left-hand side of (\ref{eq-homotopy-U}) can be re-written as 
\begin{align}
\label{eq-hom1}\mathcal U_\gamma(\alpha\cup\beta)&=\sum_{n\geq 0}\frac{1}{n!}\mathcal U_{n+2}^0(\alpha,\beta,\underset{n}{\underbrace{\gamma,\dots,\gamma}}),\\
\label{eq-hom2}\mathcal U_\gamma(\alpha)\cup\mathcal U_\gamma(\beta)&=\sum_{n\geq 0}\frac{1}{n!}\mathcal U_{n+2}^1(\alpha,\beta,\underset{n}{\underbrace{\gamma,\dots,\gamma}}),
\end{align}
where the morphisms $\mathcal U_{n+2}^i$, $i=0,1$, $n\geq 0$, are defined as in (\ref{eq-taylorU}), the only difference being that the weights (\ref{eq-weight}) have been replaced by 
\[
W_\Gamma^i=\int_{\mathcal Z_{n+2,m,i}^+}\omega_\Gamma,\ i=0,1,
\]
for any graph $\Gamma$ in $\mathcal G^K_{n+2,m}$.
We refer to~\cite{MT,CR} for an explicit proof of (\ref{eq-hom1}) and (\ref{eq-hom2}), which we omit here.

Stokes' Theorem implies the following identity between weights:
\[
W_\Gamma^1-W_\Gamma^0=\widetilde{W}_\Gamma:=\int_{Z_{n+2,m}^+}\omega_\Gamma,
\]
for any graph $\Gamma$ in $\mathcal G^K_{n+2,m}$.
Hence, the proof of (\ref{eq-homotopy-U}) is equivalent to evaluating explicitly the weights on the right-hand 
side of the previous identity.

We inspect more carefully the weights $\widetilde{W}_\Gamma$, for a general graph $\Gamma$ in $\mathcal G^K_{n+2,m}$.
For this purpose, we consider $Z_{n+2,m}^+$: as was observed earlier, $Z_{n+2,m}^+$ is the intersection of boundary strata of codimension $1$ of $\mathcal C_{n+2,m}^+$ with the interior of $\mathcal Z_{n+2,m}^+$.
Hence, the possible strata $Z$ of $Z_{n+2,m}^+$ are of the following form, recalling that in $Z_{n+2,m}^+$, the first two points in $\mathcal H$ remain distinct and lie on the curve $\gamma$:
\begin{enumerate}
\item[$i)$] there is a subset $A$ of $[n+2]$, which contains {\bf at most} one of the first two points of $\mathcal C_{n+2,m}^+$, such that 
\[
Z\cong (\mathcal C_A\times \mathcal C_{n-|A|+3,m}^+)\cap \overset{\circ}{\mathcal Z}_{n+2,m}^+=\mathcal C_A\times \overset{\circ}{\mathcal Z}_{n-|A|+3,m}^+,
\]
where $2\leq |A|\leq n+1$.
\item[$ii)$] There are a subset $A$ of $[n+2]$ and a subset $B$ of $[m]$ of successive integers, such that $A$ does not contain neither the first nor the second point of $\mathcal C_{n+2,m}^+$, such that 
\[
Z\cong (\mathcal C_{A,B}^+\times \mathcal C_{n-|A|+2,m-|B|+1}^+)\cap \overset{\circ}{\mathcal Z}_{n+2,m}^+=\mathcal C_{A,B}^+\times \overset{\circ}{\mathcal Z}_{n-|A|+2,m-|B|+1}^+,
\]
where $0\leq |A|\leq n$ and $0\leq |B|\leq m$.
\item[$iii)$] There is a subset $A$ of $[n+2]$, which contains both the first two points of $\mathcal C_{n+2,m}^+$, such that 
\[
Z\cong (\mathcal C_{A,B}^+\times \mathcal C_{n-|A|+2,m-|B|+1}^+)\cap \overset{\circ}{\mathcal Z}_{n+2,m}^+=\overset{\circ}{\mathcal Z}_{A,B}^+\times \mathcal C_{n-|A|+2,m-|B|+1}^+,
\]
where $2\leq |A|\leq n+2$ and $0\leq |B|\leq m$.
\end{enumerate}
We consider now the restriction of weights to strata of type $i)$, $ii)$ and $iii)$ of $Z_{n+2,m}^+$.

\subsubsection{Strata of type $i)$}\label{sss-4-2-1}

Graphically, a typical configuration of points in a general component of the boundary stratum $Z_{n,m}^+$ of type $i)$ looks like as follows:
\bigskip
\begin{center}
\resizebox{0.375 \textwidth}{!}{\input{K-eye-h_iii_1.pstex_t}}\\
\text{Figure 14 - A typical configuration in a general component of the boundary stratum of $Z_{n,m}^+$ of type $i)$} \\
\end{center}
\bigskip
We have two subcases of $i)$, namely, when $i_1)$ exactly one of the first two points is in $A$, or $i_2)$ neither of them is in $A$.
Any weight (\ref{eq-weight}) splits as
\[
\int_Z \omega_\Gamma=\int_{\mathcal C_A}\omega_{\Gamma_A} \int_{\overset{\circ}{\mathcal Z}_{n-|A|+3,m}^+}\omega_{\Gamma^A},
\]
where $\Gamma_A$, resp.\ $\Gamma^A$, denotes the subgraph of $\Gamma$, whose vertices are labelled by $A$ and whose edges have both endpoints in $A$, resp.\ obtained by contracting the subgraph $\Gamma_A$ to a single vertex.

By Kontsevich's Lemma~\ref{l-K-1}, Appendix~\ref{app}, the first integral on the right-hand side does not vanish, only if $|A|=2$ and $\Gamma_A$ consists of a single edge connecting the two points in $A$, in which case, by Lemma~\ref{l-angle}, it equals $1$.
Therefore, the graph $\Gamma^A$ is in $\mathcal G^K_{n+1,m}$; the weighted sum of polydifferential operators associated to subgraphs $\Gamma^A$ corresponds either to the action of $\alpha$ or $\beta$ on $\gamma$, in case $i_1)$, or, in case $i_2)$, to the action of $\gamma$ on itself w.r.t.\ the Schouten--Nijenhuis brackets.

The sum over all possible admissible graphs $\Gamma$, whose splitting as above is non-trivial, of the corresponding weights and polydifferential operators gives the first two terms of (\ref{eq-homotopy-U}), up to $\mathrm d_\mathfrak n$, and polydifferential operators containing the Schouten--Nijenhuis brackets of $\gamma$ with itself.

\subsubsection{Strata of type $ii)$}\label{sss-4-2-2}

For a better understanding, here is the graphical representation of a typical configuration of points in a general component of the boundary stratum $Z_{n,m}^+$ of type $ii)$:
\bigskip
\begin{center}
\resizebox{0.375 \textwidth}{!}{\input{K-eye-h_iii_2.pstex_t}}\\
\text{Figure 15 - A typical configuration in a general component of $Z_{n,m}^+$ of type $ii)$} \\
\end{center}

\bigskip 

We consider the case $ii)$: for any graph $\Gamma$ in $\mathcal G^K_{n+2,m}$, the weight $\widetilde{W}_\Gamma$ 
restricted to a component $Z$ of $Z_{n+2,m}$ splits as
\[
\int_Z \omega_\Gamma=\underset{\overset{\circ}{\mathcal Z}_{n-|A|+2,m-|B|+1}^+}
{\int}\left(\int_{\mathcal C_{A,B}^+}\omega_{\Gamma_{A,B}}\right)\omega_{\Gamma^{A,B}}\,,
\]
with notations similar to those in Subsubsection~\ref{sss-4-2-1}.

By Lemma~\ref{l-angle}, there can be no outgoing edge from $\Gamma_{A,B}$: thus, the subgraph $\Gamma_{A,B}$ is in 
$\mathcal G^K_{|A|,|B|}$, and so is $\Gamma^{A,B}$, and we have the splitting
\[
\int_Z \omega_\Gamma=\int_{\mathcal C_{A,B}^+}\omega_{\Gamma_{A,B}}
\underset{\overset{\circ}{\mathcal Z}_{n-|A|+2,m-|B|+1}^+}{\int}\omega_{\Gamma^{A,B}}\,.
\]
Therefore, the sum over all possible admissible graphs $\Gamma$ (whose splitting as above is non-trivial) of 
the corresponding weights and polydifferential operators yields 
$\mathcal H_\gamma^{\mathcal U}(\alpha,\beta)\{\widetilde\gamma\}$, 
recalling the brace identities and $\widetilde\gamma:=\mu+\mathcal U(\gamma)$.
We observe that the standard multiplication $\mu$ appears, when $A_3=\emptyset$ and 
$B_3=\{\overline 1,\overline 2\}$.

\subsubsection{Strata of type $iii)$}\label{sss-4-2-3}

Here is a pictorial representation of a typical configuration of points in a general component of the boundary stratum of $Z_{n,m}^+$ of type $iii)$: 
\bigskip
\begin{center}
\resizebox{0.375 \textwidth}{!}{\input{K-eye-h_iii_3.pstex_t}}\\
\text{Figure 16 - A typical configuration in a general component of $Z_{n,m}^+$ of type $iii)$} \\
\end{center}
\bigskip
For an admissible graph in $\mathcal G_{n+2,m}$, the weight $\widetilde{W}_\Gamma$ restricted to $Z$ splits as
\[
\int_Z \omega_\Gamma=\underset{\mathcal C_{n-|A|+2,m-|B|+1}^+}{\int}\left(\int_{\overset{\circ}{\mathcal Z}_{A,B}^+}\omega_{\Gamma_{A,B}}\right)\omega^{\Gamma_{A,B}},
\]
with the same notations as before.

By Lemma~\ref{l-angle}, again, there are no outgoing edges from $\Gamma_{A,B}$, hence both $\Gamma_{A,B}$ and $\Gamma^{A,B}$ are admissible, and in fact we have the splitting
\[
\int_Z \omega_\Gamma=\underset{\mathcal C_{n-|A|+2,m-|B|+1}^+}{\int}\omega_{\Gamma_{A,B}}\int_{\overset{\circ}{\mathcal Z}_{A,B}^+}\omega_{\Gamma^{A,B}}.
\]
Therefore, summing up over all possible admissible graphs $\Gamma$, with non-trivial splitting as above, of 
the corresponding weights and polydifferential operators gives 
$\widetilde\gamma\{\mathcal H^{\mathcal U}_\gamma(\alpha,\beta)\}$; we observe, once again, that the standard 
multiplication $\mu$ appears, when $A=[n+2]$ and $|B|=m-1$. 

Finally, since $\mathrm d_\mathfrak n$ is a differential, the MC equation for $\gamma$ permits to re-insert it in 
(\ref{eq-homotopy-U}), using the Leibniz rule. 

\section{The compatibility between cap products}\label{s-5}

We now come to the proof of the compatibility between cap product in the case of $X=\mathbb R^d$.
Borrowing the notation from Sections~\ref{s-4} and~\ref{s-5}, we construct a linear operator 
$$
\mathcal H_\gamma^\mathcal S\,:\,T^\mathfrak m_{\rm poly}(V)\otimes C^{{\rm poly},\mathfrak m}(V)\,
\longrightarrow\,\Omega^\mathfrak m(V)\,,
$$
(with abuse of notations from Section~\ref{s-4}), which is required to satisfy the following homotopy property: 
\begin{equation}\label{eq-homotopy1}
\begin{aligned}
&\mathcal S_\gamma\Big(\mathcal U_\gamma(\alpha)\cap c\Big)-\alpha\cap\mathcal S_\gamma(c)
=(\mathrm d_\mathfrak m+\mathrm L_\gamma)\Big(\mathcal H_\gamma^\mathcal S(\alpha,c)\Big)
+\mathcal H_\gamma^\mathcal S\Big(\mathrm d_\mathfrak m\alpha+[\gamma,\alpha],c\Big)+\\
&+(-1)^{|\!|\alpha|\!|}\mathcal H_\gamma^\mathcal S\Big(\alpha,\mathrm d_\mathfrak mc+\mathrm b_\mathrm Hc
+\mathrm L_{\mathcal U(\gamma)}c\Big)\,,
\end{aligned}
\end{equation}
for a MCE $\gamma$ in $T_\mathrm{poly}^\mathfrak n(V)$ as in~\eqref{eq-MCel}, Section~\ref{s-4}, a general 
$\mathfrak{m}$-valued polyvector field $\alpha$ and a general Hochschild chain $c$ on $A$ with values in $A$, 
$A$ being the graded algebra of $\mathfrak m$-valued functions on $V$. 
As usual, we define $\mathcal H_\gamma^\mathcal S$ in terms of graphs and weights associated to them. 
Namely, with $\alpha$ and $c=(a_0|\cdots|a_{m})$ as above ($m\geq0$), we have 
\begin{equation}\label{eq-formulaH}
\mathcal H_{\gamma}^\mathcal S(\alpha,c):=\sum_{n\geq0}\frac1{n!}\sum_{\Gamma\in\mathcal G^S_{n+1,m+1}}
\overset{\circ}{W}_{D,\Gamma}S_{\Gamma}(\alpha,\underbrace{\gamma,\dots,\gamma}_{n\textrm{ times}};c)\,,
\end{equation}
where $\overset{\circ}{W}_{D,\Gamma}$ are certain weights defined as suitable integrals.

\subsection{A curve on Kontsevich's eye $\mathcal D_{1,1}$ and related configuration spaces}\label{ss-5-1}
We consider the curve $\ell$ on Kontsevich's eye $\mathcal D_{1,1}$, with initial point $\ell(0)$ on $\alpha$, 
and final point $\ell(1)=b$, which corresponds to the following imbedding of the open unit interval into the open configuration space $D_{1,1}$: 
\bigskip
\begin{center}
\resizebox{0.45 \textwidth}{!}{\input{BNE_curve.pstex_t}}\\
\text{Figure 17 - The curve $\ell$ in $\mathcal D_{1,1}$} \\
\end{center}
\bigskip 
For the labelling of all boundary strata of the brave new eye, we refer to Subsubsection~\ref{sss-3-2-1}.
We also observe that, under the identification $\mathcal D_{1,1}=\mathcal C_{2,0}$, the curve $\ell$ corresponds to the curve on Kontsevich's Eye of Subsection~\ref{ss-4-2}.

We consider the subset $\mathcal Y_{n,m}^+$ of $\mathcal D_{n,m}^+$, for $n\geq 1$ and $m\geq 1$, consisting of those configurations, whose projection onto $\mathcal D_{1,1}$ (which extends the natural projection from $D_{n,m}$ onto $D_{1,1}$, onto the first point in $D$ and the first point in $S^1$) belongs to the curve $\ell$. 
Pictorially,
\bigskip
\begin{center}
  \includegraphics[scale=0.35]{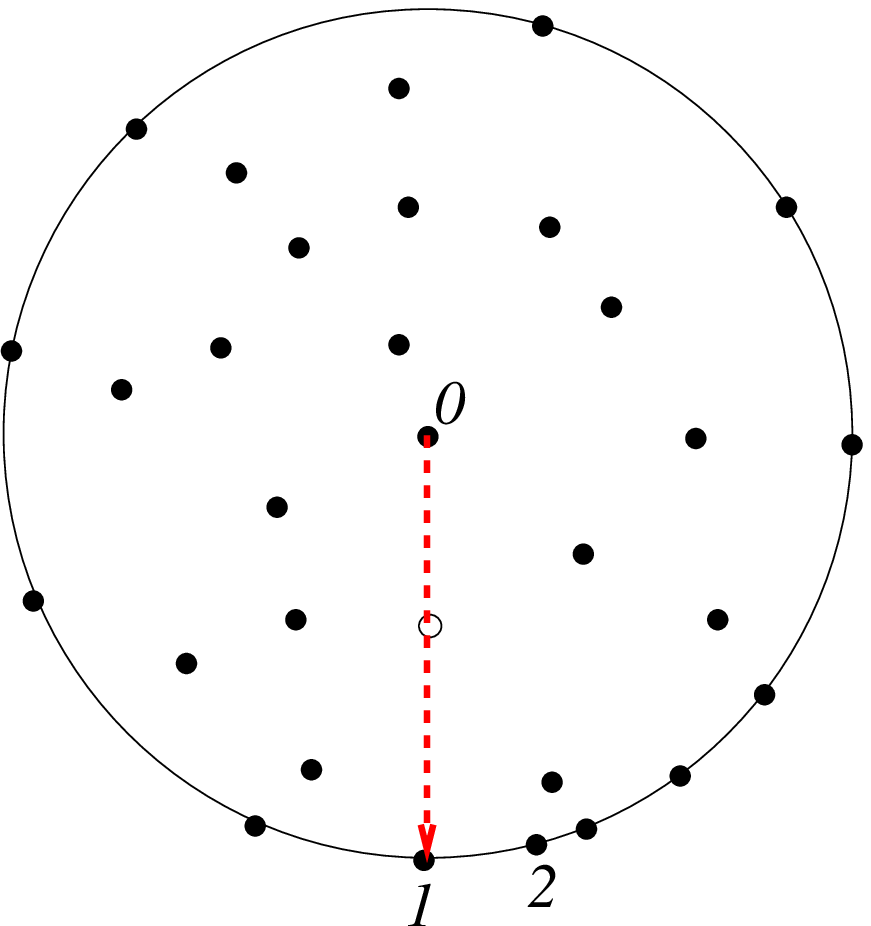} \\
\text{Figure 18 - A typical configuration in $\mathcal Y_{n,m}^+$} \\
\end{center}
\bigskip 
The dashed line represents the curve, along which the first point in $D$ moves, which we have denoted by $\circ$, connecting it with the first point in $S^1$ (w.r.t.\ the cyclic order).

The set $\mathcal Y_{n,m}^+$, for $n\geq 1$, is a smooth submanifold with corners of $\mathcal D_{n,m}^+$, of codimension $1$; it inherits an orientation from the orientation of the curve and of $\mathcal D_{n,m}^+$. 

We need to characterize explicitly the boundary strata of $\mathcal Y_{n,m}^+$ of codimension $1$: this is analogous to what was done in Subsection~\ref{ss-4-2}, whence the strata consist of
\begin{enumerate}
\item[$i)$] configurations in $\mathcal D_{n,m}^+$, whose projection onto $\mathcal D_{1,1}$ is $\ell(0)$ (the corresponding strata are denoted collectively by $\mathcal Y_{n,1,0}^+$);
\item[$ii)$] configurations in $\mathcal D_{n,m}^+$, whose projection onto $\mathcal D_{1,1}$ is $\ell(1)$ (the corresponding strata are denoted collectively by $\mathcal Y_{n,m,1}^+$);
\item[$iii)$] the intersection of boundary strata of codimension $1$ of $\mathcal D_{n,m}^+$ with the interior $\overset{\circ}{\mathcal Y}_{n,m}^+$ of $\mathcal Y_{n,m}^+$ (the corresponding strata are denoted collectively by $Y_{n,m}^+$).
\end{enumerate}
It is clear that all such boundary strata are submanifolds with corners of $\mathcal D_{n,m}^+$ of codimension $2$.

In the forthcoming Subsections, we prove that $\alpha\cap\mathcal S_\gamma(c)$, 
$\mathcal S_\gamma\big(\mathcal U_\gamma(\alpha)\cap c\big)$, and the r.h.s.~of \eqref{eq-homotopy1} can be 
expressed via a formula similar to \eqref{eq-formulaH}, involving new weights $W_{D,\Gamma}^0$, 
$W_{D,\Gamma}^1$ and $\overset{\circ}W_{D,\Gamma}$, where, 
for any S-admissible graph $\Gamma\in\mathcal G_{n,m}^S$, 
$$
W_{D,\Gamma}^i:=\int_{\mathcal Y_{n,m,i}}\omega_{D,\Gamma}\ i=0,1,\ \overset{\circ}W_{D,\Gamma}=\int_{\overset{\circ}{\mathcal Y}_{n,m}^+}\omega_{D,\Gamma}.
$$

\subsection{A formula for $\alpha\cap \mathcal S_\gamma(c)$}\label{ss-5-2}
We first consider the boundary strata $\mathcal Y_{n,m,0}^+$: it is not difficult to verify that a general component $Z$ of $\mathcal Y_{n,m,0}^+$ has the form
\[
Z\cong \mathcal D_A^0\times\mathcal D_{n-|A|,m}^+, 
\]
where $A$, $|A|\geq 1$, contains at least the first point in $D$, and $\mathcal D_A^0$ is a smooth submanifold of $\mathcal D_A$ of codimension $1$, whose elements are configurations in $\mathcal D_A$, whose projection onto the first point is fixed (since $\ell(0)$ represents in fact a point in $D$ which approaches the origin along a fixed direction in $S^1$). 
Graphically,
\bigskip
\begin{center}
  \includegraphics[scale=0.35]{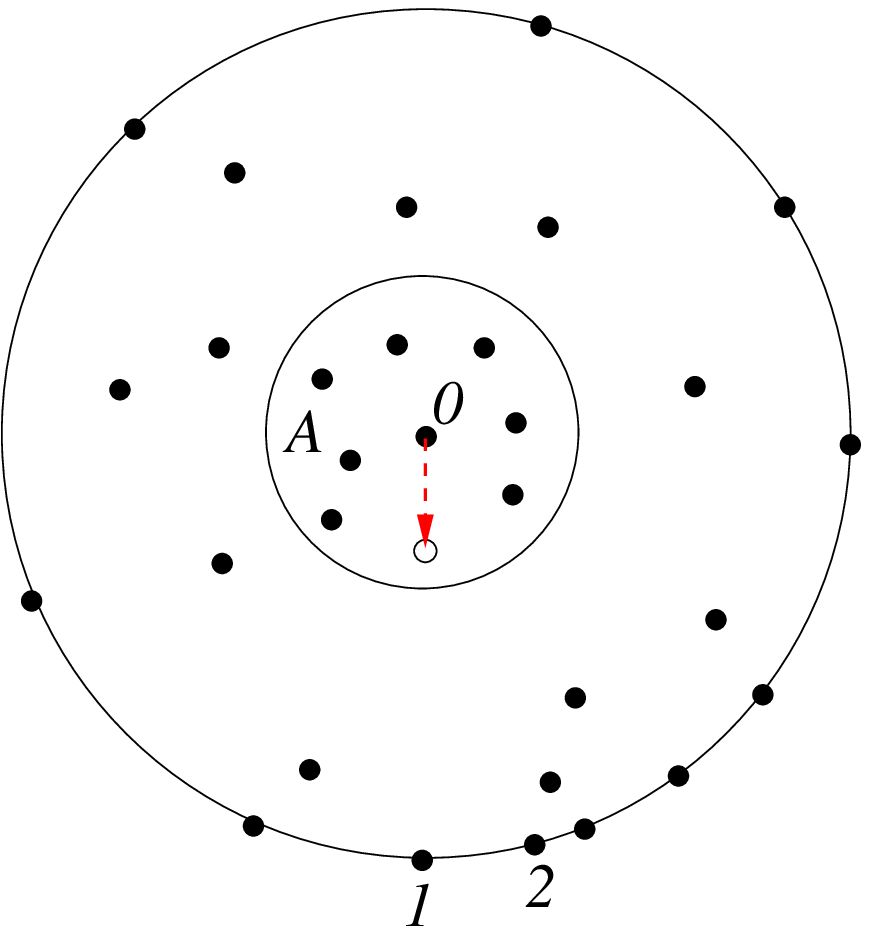} \\
\text{Figure 19 - A typical configuration in $\mathcal Y_{n,m,0}^+$} \\
\end{center}
\bigskip 
\begin{Lem}\label{l-wS-0}
For any admissible graph $\Gamma$ in $\mathcal G^S_{n+1,m+1}$, the weight $W_{D,\Gamma}^0$ vanishes if the vertex $1$ 
has at least one incoming edge. 
Otherwise, the identity
\[
W_{D,\Gamma}^0=W_{D,\Gamma_0}
\]
holds true, where $\Gamma_0\in\mathcal G^S_{n,m+1}$ is obtained from $\Gamma$ by collapsing 
the vertices $0$ and $1$. 
\end{Lem}
\begin{proof}
By the above characterization of the components $Z$ of $\mathcal Y_{n+1,m,0}^+$, for a general admissible graph $\Gamma$ as in the assumptions of Lemma~\ref{l-wS-0}, we have the following factorization 
\[
W_{D,\Gamma}^0\big\vert_Z=\int_{\mathcal D^+_{n-|A|+1,m+1}}\left(\int_{\mathcal D_A^0}\omega_{D,\Gamma_A}\right)\omega_{D,\Gamma^A},
\]
where $Z$ is a general component of $\mathcal Y_{n,m,0}^+$, and where $\Gamma_A$ is the subgraph of $\Gamma$, whose vertices are labelled by $A$ and whose edges have {\bf at least} one endpoint in $A$, and $\Gamma^A$ is obtained from $\Gamma$ by contracting $\Gamma_A$ to the single vertex labelled by $0$.

We focus our attention on the inner integral in the previous factorization.

We first observe that, if there is exactly one edge from $0$ to $1$, then, by means of Lemma~\ref{l-angle}, $a)$, the integrand vanishes, since $\omega|_{\mathcal D_1=\mathcal C_2}$ is the derivative of a constant angle, hence it is trivial. 

Further, using Lemma~\ref{l-angle}, a), and Lemma~\ref{l-angleD}, $c)$ and $f)$ (and the characterization of the restriction of $\omega_D$ to the boundary stratum $\mathcal C_3\times \mathcal C_{1,0}$ of $\mathcal C_{3,0}$) to explicitly evaluate the integrand $\omega_{D,\Gamma_A}$ on $\mathcal D_A^0$, and by dimensional reasons, all factors in $\omega_{D,\Gamma_A}$ which live on $\mathcal D_A^0$ are products of Kontsevich's angle function on $\mathcal C_2$.
We can therefore apply the arguments of the proof of Kontsevich's Lemma~\ref{l-K-1}, Appendix~\ref{app}, whence the only possibly non-trivial integrals appear, when $|A|=1$, i.e.\ $\mathcal D_A^0$ consists of a single point on $S^1$.

Since there are no edges connecting $0$ with $1$, when $|A|=1$, there can be {\bf some edge} from $1$ to some other vertices (of first or second type), or {\bf some edge} with endpoint $1$.
 
We write $\mathrm{star}(1)=\{e_1^1,\dots,e_1^p\}$.
In the first case, by Lemma~\ref{l-angleD}, $c)$, the integrand $\omega_{D,\Gamma_A}$ is a product of the form
\[
\omega_{D,\Gamma_A}|_{\mathcal D_A}=\bigwedge_{k=1}^p(\omega_{D,e_1^k}+\omega|_{\mathcal C_2}).
\]
Again, since Kontsevich's angle function is constant by construction, $\omega|_{\mathcal C_2}$ vanishes, the inner integration over a $0$-dimensional point may be then discarded, and the only non-trivial factor surviving integration is $\bigwedge_{k=1}^p\omega_{D,e_1^k}$.
The previous product can be re-inserted into the remaining integrand $\omega_{D,\Gamma^A}$, and, denoting by $\Gamma_0$ this new graph, the claim follows.

In the second case, Lemma~\ref{l-angleD}, $d)$, immediately implies the claim.
\end{proof}
\begin{Prop}\label{p-cap-0}
For $\gamma$, $\alpha$ and $c$ as above, the following identity holds true:
\begin{equation}\label{eq-cap-0}
\alpha\cap \mathcal S_\gamma(c)=\sum_{n\geq 0}\frac{1}{n!}\sum_{\Gamma\in\mathcal G^S_{n+1,m+1}}
W_{D,\Gamma}^0 S_\Gamma(\alpha,\underset{n}{\underbrace{\gamma,\dots,\gamma}};c)\,.
\end{equation}
\end{Prop}
\begin{proof}[(Sketch of proof)]
On the one hand, the l.h.s.~of (\ref{eq-cap-0}) can be re-written as 
\[
\sum_{n\geq 0}\frac{1}{n!}\sum_{\Gamma_0\in\mathcal G^S_{n,m+1}}W_{D,\Gamma_0}\ 
\iota_\alpha S_{\Gamma_0}(\underset{n}{\underbrace{\gamma,\dots,\gamma}};c)\,.
\]
On the other hand, we consider the r.h.s.~of (\ref{eq-cap-0}): by Lemma~\ref{l-wS-0}, the weights 
$W_{D,\Gamma}^0$ are non-trivial only for those admissible diagrams $\Gamma$, whose vertex labelled 
by $1$ has no incoming edges, in which case the sum simplifies to 
\[
\sum_{n\geq 0}\frac{1}{n!}\sum_{\Gamma_0\in\mathcal G^S_{n,m+1}}
\sum_{\Gamma\in\mathcal G^S_{n+1,m+1}\atop\Gamma_0\prec\Gamma}
W_{D,\Gamma_0} S_\Gamma(\alpha,\underset{n}{\underbrace{\gamma,\dots,\gamma}};c)\,.
\]
In the previous expression, the third summation is exactly over those admissible graphs, whose vertex $1$ has no incoming edges and whose contraction of the vertices $0$ and $1$ is the admissible graph $\Gamma_0$.

Finally, we observe that for any fixed graph $\Gamma_0\in\mathcal G_{n,m+1}^S$, 
$$
\iota_\alpha S_{\Gamma_0}(\underset{n}{\underbrace{\gamma,\dots,\gamma}};c)
=\sum_{\Gamma\in\mathcal G^S_{n+1,m+1}\atop \Gamma_0\prec\Gamma}
S_\Gamma(\alpha,\underset{n}{\underbrace{\gamma,\dots,\gamma}};c)\,,
$$
which ends the proof of the proposition. 
\end{proof}

\subsection{A formula for $\mathcal S_\gamma\big(\mathcal U_\gamma(\alpha)\cap c\big)$}\label{ss-5-3}

A general component $Z$ of the boundary stratum $\mathcal Y_{n,m,1}^+$ has the explicit form
\begin{equation}\label{eq-b-q}
Z\cong \mathcal C_{A_1,B_1}^+\times\mathcal C_{A_2,B_2}^+\times\mathcal D_{A_3,B_3}^+,  
\end{equation}
where $A_i$, $i=1,2,3$, are disjoint subsets of $[n]$, with $1\leq |A_1|\leq n$, $0\leq |A_2|\leq n$, $0\leq |A_3|\leq n$, with $n=|A_1|+|A_2|+|A_3|$, and $B_i$, $i=1,2,3$, are disjoint ordered subsets of $[m]$ of successive elements, such that $1\leq |B_1|\leq m$, $2\leq |B_2|\leq m$, $1\leq |B_3|\leq m$, and $m=|B_1|+|B_2|+|B_3|$.

Here is a pictorial representation of a typical component $Z$ of $\mathcal Y_{n,m,1}^+$:
\bigskip
\begin{center}
  \includegraphics[scale=0.35]{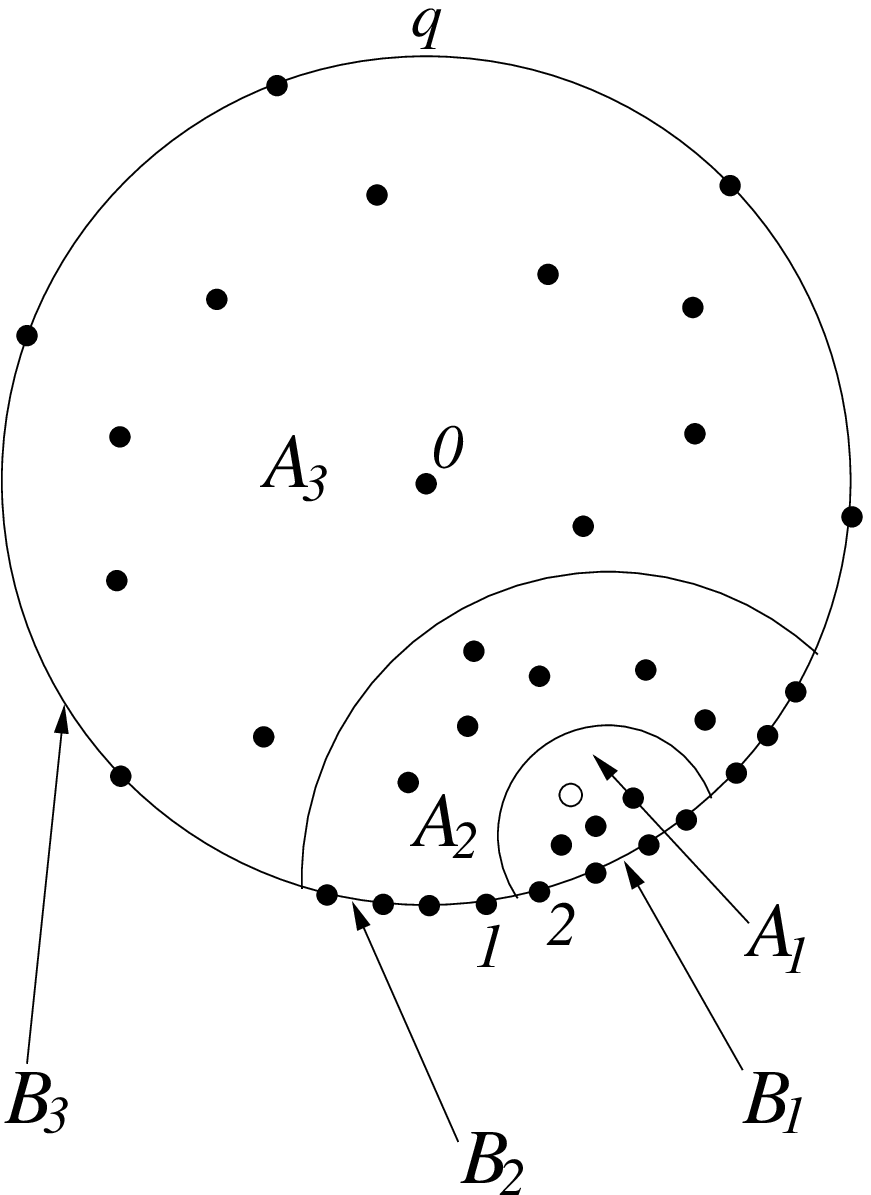} \\
\text{Figure 20 - A typical configuration in $\mathcal Y_{n,m,1}^+$} \\
\end{center}
\bigskip
We consider a component $Z$ of $\mathcal Y_{n+1,m+1}^1$ as in (\ref{eq-b-q}), and, for an S-admissible graph 
$\Gamma$ in $\mathcal G_{n+1,m+1}^S$ as before, we denote by 
\begin{itemize}
\item[$i)$] $\Gamma^Z_1$ the subgraph of $\Gamma$, whose vertices are labelled by $A_1\sqcup B_1$;
\item[$ii)$] $\Gamma^Z_2$ the graph with vertices labelled by $A_2\sqcup (B_2\sqcup\{*\})$, which is the quotient of the subgraph $\widetilde\Gamma$ (with vertices labelled by $(A_1\sqcup A_2)\sqcup (B_1\sqcup B_2)$) by $\Gamma^Z_1$, and $*$ corresponds to the contraction of $\Gamma_1^Z$;  
\item[$iii)$] $\Gamma^Z_3$ the graph with vertices labelled by $(A_3\sqcup\{0\})\sqcup (B_3\sqcup\{\overline{1}\})$, which is the quotient of $\Gamma$ by $\widetilde\Gamma$.  
\end{itemize}
\begin{Lem}\label{l-wS-1}
For a general component $Z$ of $\mathcal Y_{n+1,m+1,1}^+$ as in (\ref{eq-b-q}) and a general $S$-admissible graph $\Gamma$ in $\mathcal G_{n+1,m+1}^S$, the identity 
\begin{equation}\label{eq-wS-1}
\int_Z\omega_{D,\Gamma}=W_{\Gamma_1^Z}W_{\Gamma_2^Z}W_{D,\Gamma_3^Z}
\end{equation}     
holds true. 
\end{Lem}
\begin{proof}
First of all, if there is an edge e.g.\ from $A_1$ to $A_2$, we may apply Lemma~\ref{l-angleD}, $a)$ or $d)$ to show that the corresponding contribution vanishes; the same argument implies the claim in all other cases.
We observe that this also implies that $\Gamma_i$, $i=1,2,3$, is admissible.

Hence, we have the factorization 
\[
W_{D,\Gamma}^1|_{\mathcal C^+_{A_1,B_1}\times\mathcal C_{A_2,B_2}^+\times \mathcal D^+_{A_3,B_3}}=\int_{\mathcal C_{A_1,B_1}^+}\omega_{D,\Gamma_1}\int_{\mathcal C_{A_2,B_2}^+}\omega_{D,\Gamma_2}\int_{\mathcal D_{A_3,B_3}^+}\omega_{D,\Gamma_3}.
\]
Finally, we use Lemma~\ref{l-angleD}, $e)$, to reduce the first two factors in the previous factorization to usual Kontsevich's weights (\ref{eq-weight}): in fact, we have 
\[
\omega_{D,\Gamma_i}|_{\mathcal C_{A_i,B_i}^+}=\omega_{\Gamma_i},\ i=1,2,
\]
and the claim follows directly from the definition of (\ref{eq-weight}).
\end{proof}
Hence, Lemma~\ref{l-wS-1} implies, more generally, the following factorization property:
\begin{equation}\label{eq-factor3}
W_{D,\Gamma}^1=\sum_Z W_{\Gamma_1^Z}W_{\Gamma_2^Z}W_{D,\Gamma_3^Z}\,,
\end{equation}
where $Z$ runs over components of the type (\ref{eq-b-q}) of $\mathcal Y_{n+1,m+1,1}^+$.
\begin{Prop}\label{p-cap-1}
For $\gamma$, $\alpha$ and $c$ as above, the following identity holds true:
\begin{equation}\label{eq-cap-1}
\mathcal S_\gamma(\mathcal U_\gamma(\alpha)\cap c)=\sum_{n\geq 0}\frac{1}{n!}
\sum_{\Gamma\in\mathcal G^S_{n+1,m+1}}W_{D,\Gamma}^1
S_\Gamma(\alpha,\underset{n}{\underbrace{\gamma,\dots,\gamma}};c)\,.
\end{equation}
\end{Prop}
\begin{proof}[Sketch of proof]
We consider the left-hand side of (\ref{eq-cap-1}): it can be re-written as 
\begin{align*}
&\sum_{n_1,n_2,n_3\geq0 \atop m_1,m_2,m_3\geq0}\frac{1}{n_1!n_2!n_3!}
\sum_{\Gamma_1\in\mathcal G^K_{n_1+1,m_1}}
\sum_{\Gamma_2\in\mathcal G^K_{n_2,m_2+1}}
\sum_{\Gamma_3\in\mathcal G^S_{n_3,m_3+1}}
W_{\Gamma_1}W_{\Gamma_2}W_{D,\Gamma_3}
\sum_{0\leq k\leq l\leq p\leq m_1+m_2+m_3 \atop p-l=m_3}\\
& S_{\Gamma_3}\Big(\underset{n_3}{\underbrace{\gamma,\dots,\gamma}};
\big(U_{\Gamma_2}(\underset{n_2}{\underbrace{\gamma,\dots,\gamma}})(a_{p+1},\dots,a_k,
U_{\Gamma_1}(\alpha,\underset{n_1}{\underbrace{\gamma,\dots,\gamma}})(a_{k+1}\dots,a_{k+m_1}),\dots,a_l)\big)
|\dots|a_p\Big).
\end{align*}
As for the right-hand side of (\ref{eq-cap-0}), we apply Lemma~\ref{l-wS-1}:
$$
\sum_{n\geq 0}\frac{1}{n!}\sum_{\Gamma\in\mathcal G^S_{n+1,m+1}}\sum_Z
W_{\Gamma_1^Z}W_{\Gamma_2^Z}W_{D,\Gamma_3^Z}
S_{\Gamma}(\alpha,\underset{n}{\underbrace{\gamma,\dots,\gamma}};c),
$$
where, for an admissible graph $\Gamma$, $Z$ runs over all possible decompositions of $\Gamma$ 
into three admissible graphs as above. 
Finally, for any triple $(\Gamma_1,\Gamma_2,\Gamma_3)$ and any component $Z$ 
of $\mathcal Y_{n+1,m+1}^1$ as above, one can show that\footnote{Here we assume that 
$n_1+n_2+n_3=n$ and $m_1+m_2+m_3=m+1$. }
\begin{align*}
& S_{\Gamma_3}\Big(\underset{n_3}{\underbrace{\gamma,\dots,\gamma}};
\big(U_{\Gamma_2}(\underset{n_2}{\underbrace{\gamma,\dots,\gamma}})(a_{p+1},\dots,a_k,
U_{\Gamma_1}(\alpha,\underset{n_1}{\underbrace{\gamma,\dots,\gamma}})(a_{k+1}\dots,a_{k+m_1}),\dots,a_l)\big)
|\dots|a_p\Big) =\\
& =\sum_{\Gamma\in\mathcal G^S_{n+1,1}\atop \Gamma^Z_i=\Gamma_i}
S_{\Gamma}(\underset{n}{\alpha,\underbrace{\gamma,\dots,\gamma}},c)\,. 
\end{align*}
We observe that the component $Z$ determines the indices $k, l, p$. To finish the 
proof of the Proposition, it remains to compute the number of elements in the sum of 
the r.h.s. of the last identity: we let the reader check that it is precisely $\frac{n!}{n_1!n_2!n_3!}$. 
\end{proof}
\begin{Rem}
In the case $m\geq1$, the projection $\mathcal D_{n+1,m+1}^+\twoheadrightarrow\mathcal D_{1,1}$ 
factors through $\mathcal D_{n+1,m+1}^+\twoheadrightarrow\mathcal D_{1,2}^+\twoheadrightarrow\mathcal D_{1,1}$. 
Moreover, the inverse image of $\ell(1)$ through the last projection consists of the union of the 
following components of the I-cube (see Figure 8): $j$, $q$, $p$, $o$ and $h$. 
The detailed contribution of the inverse image of each of these components 
through $\mathcal D_{n+1,m+1}^+\twoheadrightarrow\mathcal D_{1,2}^+$ is as follows: 
$j)$ $k>0$ and $p\neq m$, 
$q)$ $k=0$, $m_1\neq0$ and $p\neq m$, 
$p)$ $m_1=0$ and $m_2\neq0$, 
$o)$ $m_1=m_2=0$ and $p\neq m$, 
$h)$ $p=m$ and $m_1=m_2=0$. 
\end{Rem}

\subsection{The homotopy formula}\label{ss-5-4}
Summarizing the results of Propositions~\ref{p-cap-0} and~\ref{p-cap-1}, we may write
\[
\mathcal S_\gamma(\mathcal U_\gamma(\alpha)\cap c)-\alpha\cap \mathcal S_\gamma(c)
=\sum_{n\geq 0}\frac{1}{n!}\sum_{\Gamma\in\mathcal G^S_{n+1,m+1}}(W_{D,\Gamma}^1-W_{D,\Gamma}^0)
S_{\Gamma}(\alpha,\underset{n}{\underbrace{\gamma,\dots,\gamma}};c).
\]

Then, for any S-admissible graph $\Gamma\in\mathcal G^S_{n,m}$, we have
\[
0=\int_{\mathcal Y_{n,m}^+}\mathrm d\omega_{D,\Gamma}=\int_{\partial \mathcal Y_{n,m}^+}\omega_{D,\Gamma}
=W^1_{D,\Gamma}-W^0_{D,\Gamma}-\overset{\circ}W_{D,\Gamma}\,,
\]
where the second equality is a consequence of Stokes' Theorem. 
It thus remains to analyze the contributions coming from integration over components of $Y_{n,m}^+$, 
which are of three types: 
\begin{itemize}
\item[$i)$] there is a subset $A$ of $[n]$ which may or may not contain the vertex $1$, such that
\[
Z\cong (\mathcal C_A\times \mathcal D_{n-|A|+1,m}^+)\cap \overset{\circ}{\mathcal Y}_{n,m}^+\cong \mathcal C_A\times \overset{\circ}{\mathcal Y}_{n-|A|+1,m}^+
\]
and $2\leq |A|\leq n$.
\item[$ii)$] There is a subset $A$ of $[n]$, which does not contain the vertex $1$, such that
\[
Z\cong (\mathcal D_A\times \mathcal D_{n-|A|,m}^+)\cap \overset{\circ}{\mathcal Y}_{n,m}^+\cong \mathcal D_A\times \overset{\circ}{\mathcal Y}_{n-|A|,m}^+\,.
\]
\item[$iii)$] There is a subset $A$ of $[n]$, which does not contain the vertex $1$, and an ordered subset $B$ of successive elements in $[m]$, such that
\[
Z\cong (\mathcal C_{A,B}^+\times \mathcal D_{n-|A|,m-|B|+1}^+)\cap \overset{\circ}{\mathcal Y}_{n,m}^+\cong
\mathcal C_{A,B}^+\times \overset{\circ}{\mathcal Y}_{n-|A|,m-|B|+1}^+\,.
\]
\end{itemize}
Namely, such components are intersections of boundary strata of codimension $1$ of $\mathcal D_{n,m}^+$ 
with the interior of $\pi^{-1}\big(\ell(]0,1[)\big)$. 
Hence, a configuration in $Y_{n,m}^+$ is such that the first point in $D^\times$ approaches neither the origin nor $S^1$. 

For any S-admissible graph $\Gamma\in\mathcal G^S_{n,m}$ and any component $Z$ of $Y_{n,m}^+$, 
we write
\begin{equation}\label{eq-wS-i}
W_{D,\Gamma}^{Z}:=\int_Z\omega_{D,\Gamma}\,.
\end{equation}

\subsubsection{Contribution of components of type $i)$}\label{sss-5-4-1}
We consider components of type $i)$ of $Y_{n,m}^+$: to get a better understanding of them, here is a pictorial representation of possible configurations in two distinct general components 
\bigskip
\begin{center}
  \includegraphics[scale=0.35]{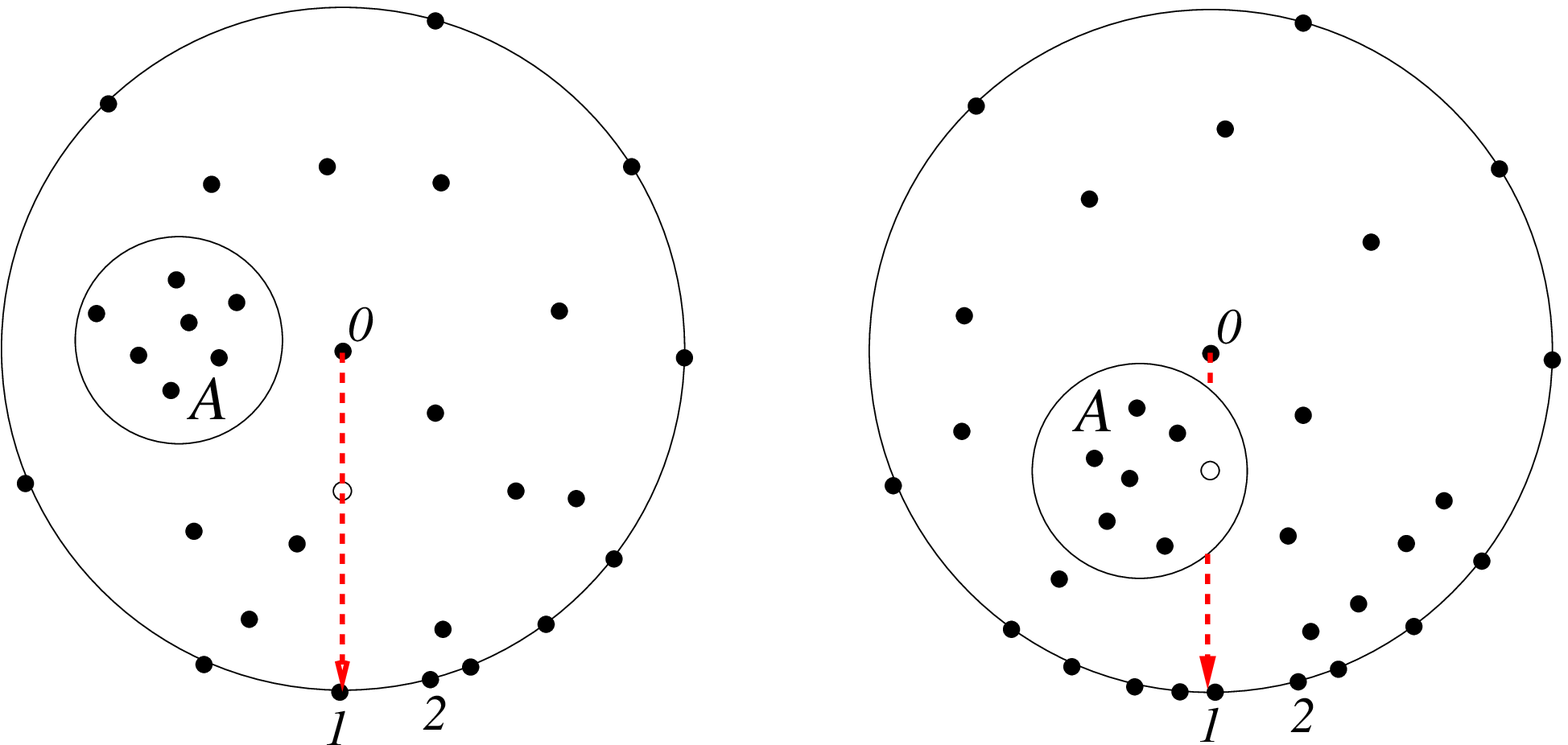} \\
\text{Figure 21 - Two distinct typical configurations of type $i)$ in $Y_{n,m}^+$ } \\
\end{center}
\bigskip
\begin{Lem}\label{l-wS-i}
For an S-admissible graph $\Gamma\in\mathcal G^S_{n+1,m+1}$ and a component $Z$ of $Y_{n+1,m+1}^+$ of type $i)$ 
as before, the weight $W_{D,\Gamma}^{Z}$ is non-trivial only if $|A|=2$ and there is exactly one edge connecting the 
vertices labelled by $A$, in which case we have
\[
W_{D,\Gamma}^{Z}=\overset{\circ}{W}_{\Gamma^A}:=\int_{\overset{\circ}{\mathcal Y}^+_{n,m+1}}\omega_{D,\Gamma^A},
\]
where $\Gamma^A$ is obtained from $\Gamma$ by contracting the vertices labelled by $A$ and eliminating 
the edge between them. 
\end{Lem}
\begin{proof}
First of all, we have the following factorization of (\ref{eq-wS-i}):
\[
W_{D,\Gamma}^{Z,i}=\int_{\overset{\circ}{\mathcal Y}^+_{n-|A|+1,m+1}}\left(\int_{\mathcal C_A}\omega_{D,\Gamma_A}\right)\omega_{D,\Gamma^A},
\]
where $\Gamma_A$, resp.\ $\Gamma^A$, is the subgraph of $\Gamma$, whose vertices are labelled by $A$ and whose edges have at least one endpoint in $A$, resp.\ the graph obtained from $\Gamma$ by contracting $\Gamma_A$ to a single vertex.

Using Lemma~\ref{l-angleD}, $f)$, we re-write the first factor in the previous factorization as 
\[
\int_{\mathcal C_A}\omega_{\Gamma_A},
\]
if $|E_{\Gamma_A}|\geq 2$, whence, by Kontsevich's Lemma~\ref{l-K-1}, we conclude that it vanishes.

Thus, $\Gamma_A$ can have {\bf at most} one edge: dimensional reasons imply that $|A|=2$.
The corresponding integral does not vanish iff $\Gamma_A$ consists of a single edge connecting the two vertices labelled by $A$: by Lemma~\ref{l-angleD}, $f)$, the contribution of such an integral is $1$ (since the integral over $\mathcal C_2=S^1$ of the piece $-\pi_2^*\omega$ of the restriction of $\omega_D$ vanishes, as $\pi_2^*\omega$ is not on $S^1$).
It is also clear that, in this case, $\Gamma^A$ is an admissible graph in $\mathcal G_{n-1,m+1}^S$. 
\end{proof}
\begin{Prop}\label{p-cap-i}
For $\gamma$, $\alpha$ and $c$ as above, the following identity holds true:
\begin{equation}\label{eq-cap-i}
\begin{aligned}
&\sum_{n\geq 0}\frac{1}{n!}\sum_{\Gamma\in\mathcal G^S_{n+1,m+1}}\sum_Z W_{D,\Gamma}^{Z}
S_{\Gamma}(\alpha,\underset{n}{\underbrace{\gamma,\dots,\gamma}};c)=\\
&=\frac12\sum_{n\geq 2}\frac{1}{(n-2)!}\sum_{\Gamma\in\mathcal G^S_{n,m+1}}\overset{\circ}{W}_{D,\Gamma}
S_{\Gamma}(\alpha,[\gamma,\gamma],\underset{n-2}{\underbrace{\gamma,\dots,\gamma}};c)+\\
&\phantom{=}+\sum_{n\geq 1}\frac{1}{(n-1)!}\sum_{\Gamma\in\mathcal G^S_{n,m+1}}\overset{\circ}{W}_{D,\Gamma}
S_{\Gamma}([\gamma,\alpha],\underset{n-1}{\underbrace{\gamma,\dots,\gamma}};c)\,,
\end{aligned}
\end{equation}
where $Z$ runs over components of type $i)$ of $Y_{n+1,m+1}^+$. 
\end{Prop}
\begin{proof}[Sketch of proof]
By Lemma~\ref{l-wS-i}, the only components $Z$ of type $i)$ yielding non-trivial weights are those of 
the form $\mathcal C_A\times \overset{\circ}{\mathcal Y}_{n,m+1}^+$, where $|A|=2$: thus, the left-hand side 
of (\ref{eq-cap-i}) can be re-written as 
\[
\sum_{n\geq 0}\frac{1}{n!}\sum_A \sum_{\Gamma^A\in\mathcal G^S_{n,m+1}}
\sum_{\Gamma\in\mathcal G^S_{n+1,m+1}\atop \Gamma^A\prec \Gamma} \overset{\circ}{W}_{D,\Gamma^A}
S_{\Gamma}(\alpha,\underset{n}{\underbrace{\gamma,\dots,\gamma}};c),
\]
where the notation $\Gamma^A\prec \Gamma$ means that $\Gamma^A$ is obtained from $\Gamma$ by collapsing the vertices 
of $\Gamma$ labelled by $A$ and the only edge between them; the second sum is over all $A\subset[n]$ such that $|A|=2$ 
(i.e.~over the above components $Z$ of type $i)$). 

Terms involving $[\gamma,\alpha]$, resp.~$\frac12[\gamma,\gamma]$, correspond to components $Z$ for which $A$ 
contains $1$, resp.~does not contain $1$. 
\end{proof}
Now, combining the MC equation for $\gamma$ with Leibniz's rule, we get 
$$
(n-1)S_{\Gamma}(\alpha,[\gamma,\gamma],\underset{n-2}{\underbrace{\gamma,\dots,\gamma}};c)
={\rm d}_{\mathfrak m}\big(S_{\Gamma}(\alpha,\underset{n-1}{\underbrace{\gamma,\dots,\gamma}};c)\big)
+S_{\Gamma}\big({\rm d}_{\mathfrak m}(\alpha),\underset{n-1}{\underbrace{\gamma,\dots,\gamma}};c\big)
\pm S_{\Gamma}\big(\alpha,\underset{n-1}{\underbrace{\gamma,\dots,\gamma}};{\rm d}_{\mathfrak m}(c)\big))\,,
$$
which implies, thanks to Proposition \ref{p-cap-i}, that 
\begin{equation}\label{eq-cap-i-bis}
\begin{aligned}
&\sum_{n\geq 0}\frac{1}{n!}\sum_{\Gamma\in\mathcal G^S_{n+1,m+1}}\sum_Z W_{D,\Gamma}^{Z}
S_{\Gamma}(\alpha,\underset{n}{\underbrace{\gamma,\dots,\gamma}};c)=\\
&={\rm d}_{\mathfrak m}\big(\mathcal H_\gamma^\mathcal S(\alpha,c)\big)
+\mathcal H_\gamma^\mathcal S\big({\rm d}_{\gamma}(\alpha),c\big)
\pm \mathcal H_\gamma^\mathcal S\big(\alpha,{\rm d}_{\mathfrak m}(c)\big)\,,
\end{aligned}
\end{equation}

\subsubsection{Contribution of components of type $ii)$}\label{sss-5-4-2}

We begin with a pictorial representation of a general configuration in a possible component $Z$ of type $ii)$ of $Y_{n,m}^+$:
\bigskip
\begin{center}
  \includegraphics[scale=0.35]{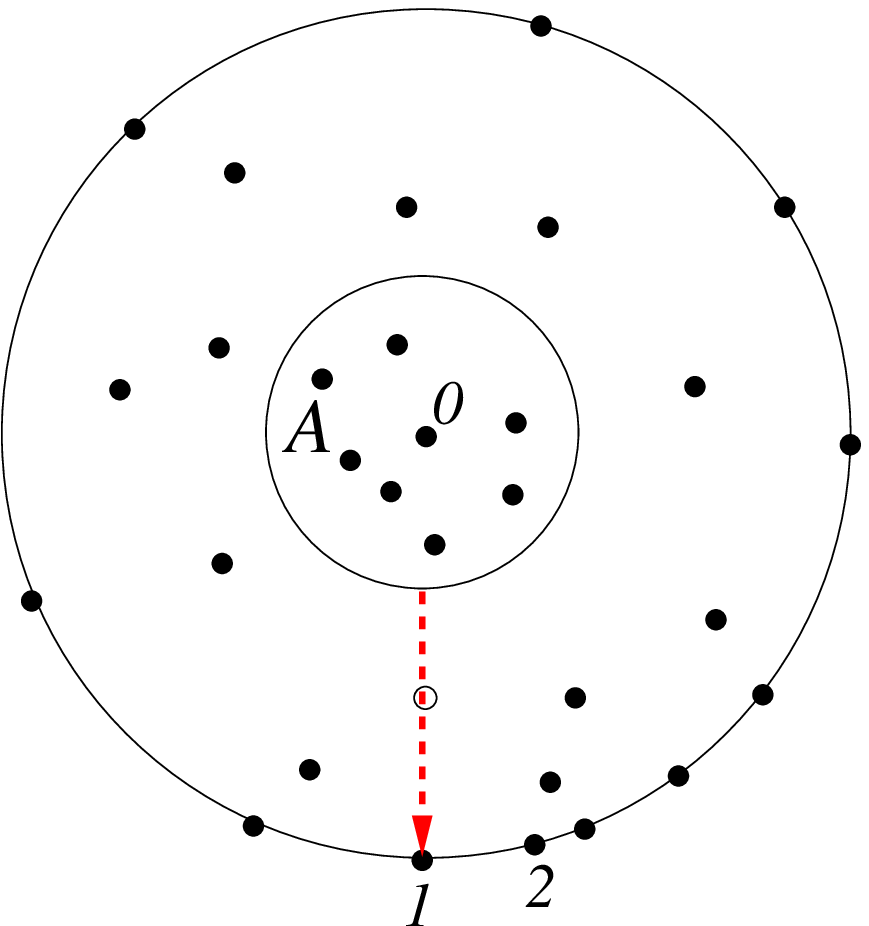} \\
\text{Figure 22 - A typical configurations of type $ii)$ in $Y_{n,m}^+$ } \\
\end{center}
\bigskip
\begin{Lem}\label{l-wS-ii}
For an S-admissible graph $\Gamma\in\mathcal G^S_{n+1,m+1}$ and a component $Z$ of $Y_{n+1,m+1}^+$ 
of type $ii)$ as before, the weight $W_{D,\Gamma}^Z$ is non-trivial only if $|A|=1$ 
and if either $a)$ there is exactly one edge from $0$ to the vertex labelled by $k$, 
or $b)$ there is no edge from $0$ to $k$ and there is at least one edge with starting point labelled by $k$.
We then have
\begin{align*}
W_{D,\Gamma}^{Z}&=\overset{\circ}{W}_{D,\Gamma^A},\ \text{in case $a)$},\\
W_{D,\Gamma}^{Z}&=\sum_{e_{v_A}^l\in\mathrm{star}(v_A)}\pm \overset{\circ}{W}_{D,\Gamma^{A,l}},\ \text{in case $b)$},
\end{align*}
where $\Gamma^A$, resp.\ $\Gamma^{A,l}$, is obtained from $\Gamma$ by contracting $0$ and the vertex 
labelled by $A$ and eliminating the edge from $0$ to this vertex, resp.\ eliminating the edge $e_{v_A}^l$. 
\end{Lem}
\begin{proof}
The decomposition of the component $Z$ implies the factorization of the corresponding weight,
\[
W_{D,\Gamma}^{Z,ii}=\int_{\overset{\circ}{\mathcal Y}^+_{n-|A|+1,m+1}}\left(\int_{\mathcal D_A}\omega_{D,\Gamma_A}\right)\omega_{D,\Gamma^A},
\]
where now $\Gamma_A$ denotes the subgraph of $\Gamma$, whose vertices are labelled by $A$ and whose edges have {\bf at least} one vertex in $A$, and $\Gamma^A$ is obtained from $\Gamma$ by contracting $\Gamma_A$ to a single vertex.

We focus on the first factor: first of all, we use Lemma~\ref{l-angle}, a), and Lemma~\ref{l-angleD}, $c)$ and $d)$ (and again the restriction of $\omega_D$ to $\mathcal C_3\times \mathcal C_{1,0}$) to evaluate the restriction of $\omega_{D,\Gamma_A}$ on $\mathcal D_A$.
All factors of this restriction, which live on $\mathcal D_A$, are pull-backs of Kontsevich's angle form $\omega$ to $\mathcal C_2$: hence, we may apply Kontsevich's Lemma~\ref{l-K-1} to conclude that all such contributions vanish unless $|A|=1$.

If $|A|=1$ (the corresponding vertex is denoted by $v_k$), and there is {\bf at least} one edge, whose endpoint is $v_k$, Lemma~\ref{l-angleD}, $d)$, implies immediately that the corresponding integral vanishes. 

We assume now $\mathrm{star}(v_A)=\{e_{v_A}^1,\dots,e_{v_A}^p\}$.
If there is exactly one edge from $0$ to $v_A$, Lemma~\ref{l-angleD}, $c)$, implies that $\omega_{D,\Gamma_A}$ is a product of the form (forgetting about signs)
\[
\omega_{D,\Gamma_A}|_{\mathcal D_A}=\omega|_{\mathcal C_2}\wedge \bigwedge_{k=1}^p(\omega_{D,e_{v_A}^k}+\omega|_{\mathcal C_2})=\omega|_{\mathcal C_2}\wedge \bigwedge_{k=1}^p\omega_{D,e_{v_A}^k},
\]
by the antisymmetry of the wedge product.
The inner integration over $S^1$ produces hence a factor $1$, and the form $\bigwedge_{k=1}^p\omega_{D,e_{v_A}^k}$ can be inserted into the outer factor $\omega_{\Gamma^A}$.

If there is no edge from $0$ to $v_A$, by Lemma~\ref{l-angleD}, $c)$, $\omega_{D,\Gamma_A}$ is a product (again forgetting about signs)
\[
\omega_{D,\Gamma_A}|_{\mathcal D_A}=\bigwedge_{k=1}^p(\omega_{D,e_{v_A}^k}+\omega|_{\mathcal C_2})=\sum_{k=1}^p(-1)^{k-1}\omega|_{\mathcal C_2}\wedge\bigwedge_{j\neq k}\omega_{D,e_{v_A}^j},
\]
again by the antisymmetry of the wedge product, since we need the factor $\omega|_{\mathcal C_2}$ because of the integral.
The inner integration over $S^1$ produces hence a factor $1$, and the forms $\bigwedge_{j\neq k}\omega_{D,e_{v_A}^j}$, $k=1,\dots,p$ can be inserted into the outer factor $\omega_{\Gamma^A}$, and the claim follows.
\end{proof}
\begin{Prop}\label{p-cap-ii}
For $\gamma$, $\alpha$ and $c$ as above, the following identity holds true:
\begin{equation}\label{eq-cap-ii}
\sum_{n\geq 0}\frac{1}{n!}\sum_{\Gamma\in\mathcal G^S_{n+1,m+1}}\sum_Z W_{D,\Gamma}^{Z}
S_{\Gamma}(\alpha,\underset{n}{\underbrace{\gamma,\dots,\gamma}};c)
=\mathrm{L_\gamma}\big(\mathcal H_\gamma(\alpha,c)\big)\,,
\end{equation}
where $Z$ runs over components of type $ii)$ of $Y_{n+1,m+1}^+$. 
\end{Prop}
\begin{proof}
Using Lemma~\ref{l-wS-ii}, we re-write the left-hand side of (\ref{eq-cap-ii}) as
\begin{align*}
&\sum_{n\geq 0}\frac{1}{n!}\sum_{k=2}^{n+1} \sum_{\Gamma^k\in\mathcal G_{n,m+1}}\sum_{\Gamma\in\mathcal G_{n+1,m+1}\atop \Gamma^k\prec \Gamma} \overset{\circ}{W}_{D,\Gamma^k}
S_{\Gamma}(\alpha,\underset{n}{\underbrace{\gamma,\dots,\gamma}};c)+\\
&+\sum_{n\geq 0}\frac{1}{n!}\sum_{k=2}^{n+1}\sum_{e_k^l\in\mathrm{star}(k)}\sum_{\Gamma^{{k,l}}\in\mathcal G_{n,m+1}}
\sum_{\Gamma\in\mathcal G_{n+1,m+1}\atop \Gamma^{{k,l}}\prec\Gamma}\pm\overset{\circ}{W}_{D,\Gamma^{{k,l}}}
S_{\Gamma}(\alpha,\underset{n}{\underbrace{\gamma,\dots,\gamma}};c),
\end{align*}
borrowing notations from Lemma~\ref{l-wS-ii}.
On the other hand, recalling the homotopy formula
\[
\mathrm L_\gamma=\mathrm d\circ\iota_\gamma\pm\iota_\gamma\circ\mathrm d\,,
\]
and using slight modifications of the arguments of Subsection 4.2 and 4.3 of~\cite{MT} and in the proof of 
Proposition~\ref{p-cap-0}, we have
\begin{align*}
\mathrm L_{\gamma}\!\left(S_{\Gamma_0}(\alpha,\underset{n-1}{\underbrace{\gamma,\dots,\gamma}};c)\right)
&=\sum_{k=2}^{n+1}\sum_{\Gamma\in \mathcal G_{n+1,m+1}^S\atop \Gamma_0=\Gamma^k}
\pm S_{\Gamma}(\alpha,\underset{n}{\underbrace{\gamma,\dots,\gamma}};c)+\\
&\phantom{=}+\sum_{k=2}^{n+1}\sum_{e_k^l\in\mathrm{star}(k)}
\sum_{\Gamma\in\mathcal G_{n+1,m+1}^S\atop \Gamma_0=\Gamma^{{k,l}}}
\pm S_{\Gamma}(\alpha,\underset{n}{\underbrace{\gamma,\dots,\gamma}};c)\,,
\end{align*} 
whence the claim follows.
More precisely, the first, resp.\ second, term on the right-hand side of the previous equality corresponds to 
the composition $\mathrm d\circ\iota_{\gamma}$, resp.\ $\iota_{\gamma}\circ \mathrm d$, by an explicit evaluation of the contraction operation and by means of Leibniz's rule for $\mathrm d$.
\end{proof}

\subsubsection{Contribution of components of type $iii)$}\label{sss-5-4-3}

We discuss weights associated to admissible graphs and to components of type $iii)$ of $Y_{n,m}^+$: before entering into the discussion, a pictorial representation of the two distinct possible configurations in such components could be helpful:
\bigskip
\begin{center}
  \includegraphics[scale=0.35]{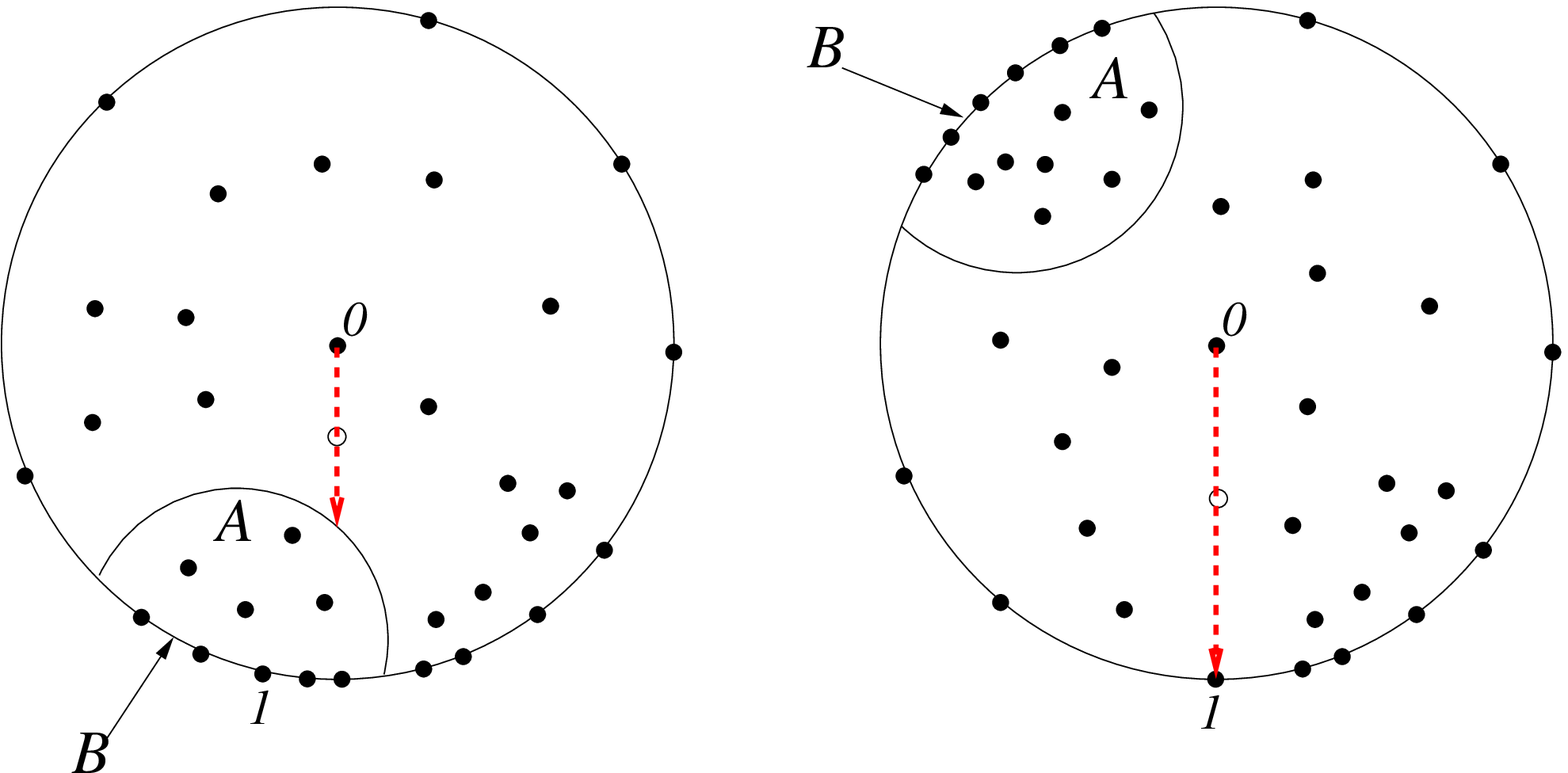} \\
\text{Figure 23 - Two possible configurations in $Y_{n,m}^+$ of type $iii)$ } \\
\end{center}
\bigskip
\begin{Lem}\label{l-wS-iii}
For an S-admissible graph $\Gamma\in\mathcal G^S_{n+1,m+1}$ and a component $Z$ of $Y^+_{n+1,m+1}$ 
of type $iii)$ as before, the weight $W_{D,\Gamma}^{Z}$ vanishes unless there are no edges connecting 
$A$ to its complement. In this case we have
\[
W_{D,\Gamma}^{Z}=W_{\Gamma_{A,B}}\overset{\circ}{W}_{D,\Gamma^{A,B}}\,,
\]
where $\Gamma_{A,B}$ denotes the subgraph of $\Gamma$ whose vertices are labelled by $A\sqcup B$ and the graph $\Gamma^{A,B}$ 
is obtained from $\Gamma$ by contracting $\Gamma_{A,B}$ to a single vertex of the second type. \hfill\qed
\end{Lem}
\begin{proof}
The above weight vanishes, if there is at least one edge connecting a vertex labelled by $A$ to a vertex not labelled by $A$ in virtue of Lemma~\ref{l-angleD}, $a)$, similarly to the first step in the proof of Lemma~\ref{l-wS-1}: this forces, by the way, $\Gamma_A$ and $\Gamma_Z$ to be both admissible, since all stars of $\Gamma_A$ and $\Gamma_Z$ belong to $E_{\Gamma_A}$ and $E_{\Gamma_Z}$ respectively.

Additionally, the following factorization of the above weight holds true:
\[
W_{D,\Gamma}^{Z}=\int_{\mathcal C_{A,B}^+}\omega_{D,\Gamma_{A,B}}\int_{\overset{\circ}{\mathcal Y}^+_{n-|A|+1,m-|B|+2}}\omega_{D,\Gamma^{A,B}}.
\]
Finally, we use Lemma~\ref{l-angleD}, $e)$, to prove that the integrand $\omega_{D,\Gamma_{A,B}}$ equals in fact $\omega_{\Gamma_{A,B}}$, whence the last claim follows from the previous factorization and from the definition of (\ref{eq-weight}).
\end{proof}
\begin{Prop}\label{p-cap-iii}
For $\gamma$, $\alpha$ and $c$ as above, the following identity holds true:
\begin{equation}\label{eq-cap-iii}
\sum_{n\geq 0}\frac{1}{n!}\sum_{\Gamma\in\mathcal G^S_{n+1,m+1}}\sum_Z W_{D,\Gamma}^{Z}
S_{\Gamma}(\alpha,\underset{n}{\underbrace{\gamma,\dots,\gamma}};c)
=\mathcal H_\gamma^\mathcal S(\alpha,\mathrm L_{\mu+\mathcal U(\gamma)}(c))\,,
\end{equation}
where $Z$ runs over components of type $iii)$ of $Y_{n+1,m+1}^+$, and $\mu$ denotes the standard multiplication in $A$. 
\end{Prop}
\begin{proof}
We use Lemma~\ref{l-wS-iii} to re-write the left-hand side of (\ref{eq-cap-iii}):
$$
\sum_{n\geq 0}\frac{1}{n!}\sum_{\Gamma\in\mathcal G^S_{n+1,m+1}\atop A,B}W_{\Gamma_{A,B}}
\overset{\circ}{W}_{D,\Gamma^{A,B}}S_{\Gamma}(\alpha,\underset{n}{\underbrace{\gamma,\dots,\gamma}};c)\,.
$$
We fix $n,m\geq0$, $0\leq k\leq n$, $0\leq l\leq m+1$, and we observe that for pair $(\Gamma_1,\Gamma_2)$ of 
admissible graphs $\Gamma_1\in\mathcal G^K_{k,l}$ and $\Gamma_2\in\mathcal G^S_{n-k+1,m-l+2}$, one finds 
$$
\sum_{A,B \atop |A|=k,|B|=l}\sum_{\Gamma\in\mathcal G^S_{n+1,m+1}\atop\Gamma_{A,B}=\Gamma_1,\Gamma^{A,B}=\Gamma_2}
S_{\Gamma}(\alpha,\underset{n}{\underbrace{\gamma,\dots,\gamma}};c)
=S_{\Gamma_2}\big(\alpha,\underset{k}{\underbrace{\gamma,\dots,\gamma}};
{\rm L}_{\mathcal U_{\Gamma_1}(\gamma,\dots,\gamma)}(c)\big)\,.
$$
We finally observe that the graph $\Gamma_1$, which consists of only two vertices of the second type, yields exactly the multiplication $\mu$.

This ends the proof of the Proposition (we leave to the reader the check of the consistency of the combinatorial coefficients). 
\end{proof}

\subsubsection{End of the proof of the compatibility between cap products in the case $X=\mathbb R^d$}

It follows from identities \eqref{eq-cap-i-bis}, \eqref{eq-cap-ii} and \eqref{eq-cap-iii} that 
$$
\sum_{n\geq0}\frac1{n!}\sum_{\Gamma\in\mathcal G_{n+1,m+1}^S} W_{D,\Gamma}^2S_\Gamma(\alpha,\gamma,\dots,\gamma;c)
$$
is precisely equal to the r.h.s.~of the homotopy equation \eqref{eq-homotopy1}, whence the result follows. 

\section{Some special cases of interest}\label{s-6}

In this Section we discuss some interesting special cases.
\begin{enumerate}
\item[$i)$] We first give the recipe for proving Shoikhet's conjecture~\cite{Sh,CR1} starting from the main result of 
this paper; we want to point out that the proof of Shoikhet's conjecture was the starting point of the investigations 
that have led us to the present paper, thanks to a stimulating question from A.S.~Cattaneo. 
\item[$ii)$] We further discuss the case of a MCE of polyvector degree less or equal than $1$, in view of a forthcoming 
application towards globalisation of the results of~\cite{MT} and of the present paper in the framework of deformation 
quantization (we refer to Section~\ref{s-8} for more details).
\item[$iii)$] We finally consider the special case of a MCE $\gamma$ of polyvector degree $0$: in this case, 
we may compute explicitly both quasi-isomorphisms $\mathcal U_\gamma$ and $\mathcal S_{\gamma}$ using results 
of~\cite{CVdB}: this is an important computational result, which will play a fundamental r\^ole in the proof 
of Caldararu's conjecture~\cite{Cald,CVdB,CRVdB}.
\end{enumerate}

\subsection{The result of~\cite{MT} and Shoikhet's conjecture}\label{ss-6-1}

We consider the special case $\mathfrak m=\mathbb R[\![\hbar]\!]$, viewed as a DGA concentrated in degree $0$ and 
trivial differential ($\mathfrak n=\hbar\mathbb R[\![\hbar]\!]$ being the pronilpotent ideal).

For $V=\mathbb R^d$, we consider the DGLA of $\hbar$-formal polyvector fields 
$T_\mathrm{poly}(V)[\![\hbar]\!]$, resp. of $\hbar$-formal 
polydifferential operators $D_\mathrm{poly}(V)[\![\hbar]\!]$, and the corresponding DGM 
$\Omega(V)^{-\bullet}[\![\hbar]\!]$, resp.\ the Hochschild 
chain complex $C_{-\bullet}^{\mathrm{poly}}[\![\hbar]\!]$. 
All algebraic structures are extended $\hbar$-linearly and tensor products are completed 
w.r.t.\ the $\hbar$-adic topology. 

We consider a MCE $\gamma$ of $\hbar T_\mathrm{poly}(V)[\![\hbar]\!]=T_\mathrm{poly}^{\mathfrak{n}}(V)$ of the form
\begin{equation}\label{eq-MC-def}
\gamma=\pi_\hbar=\hbar\pi_1+\cdots,\ \pi_i\in T_\mathrm{poly}^1(V),\ i\geq 1.
\end{equation}
Its image $\mathcal U(\gamma)$ w.r.t.\ the $L_\infty$-quasi-isomorphism $\mathcal U$ is an element of degree $1$.

Since $\mathrm d_\mathrm H=[\mu,\ ]$, where $\mu$ is the standard product on $A$, the MC equation for 
$\mathcal U(\gamma)$ is equivalent to the fact that $\mu+\mathcal U(\gamma)$ defines an $\hbar$-linear associative 
product, which we denote by $\star$, on $C^\infty(V)[\![\hbar]\!]$. 

We now consider the Gerstenhaber algebras (up to homotopy) $(T_\mathrm{poly}(V)[\![\hbar]\!],[\gamma,\ ],[\ ,\ ],\cup)$ 
and $(D_\mathrm{poly}(V)[\![\hbar]\!],\mathrm d_\mathrm H+[\mathcal U(\gamma),\ ],[\ ,\ ],\cup_\star)$; for the latter, 
we have $\mathrm d_\mathrm H+[\mathcal U(\gamma),\ ]=\mathrm d_{\mathrm H,\star}$, where 
$\mathrm d_{\mathrm H,\star}$ is the Hochschild differential on the Hochschild cochain complex of 
the $\mathbb R[\![\hbar]\!]$-algebra $C^\infty(V)[\![\hbar]\!]$, endowed with the deformed product $\star$, with values 
in itself. Furthermore, the product $\cup$ of degree $1$ takes the form
\[
(D_1\cup_\star D_2)(a_1,\dots,a_n)=D_1(a_1,\dots,a_{|D_1|+1})\star D_2(a_{|D_1|+2},\dots,a_n),\ D_i\in D_\mathrm{poly}(V)[\![\hbar]\!],\ a_i\in A,
\]
and $n=|D_1|+|D_2|+2$ (up to a sign depending on $D_i$).

On the other hand, we also consider the $T$-module, resp.\ $T$-module up to homotopy, 
$(\Omega^{-\bullet}(V)[\![\hbar]\!],\mathrm L_\gamma,\mathrm L,\cap)$, 
resp.\ 
$(C_{-\bullet}^{\mathrm{poly}}(V)[\![\hbar]\!],\mathrm b_\mathrm H+\mathrm L_{\mathcal U(\gamma)},\mathrm L,\cap_\star)$: 
in the latter, we have $\mathrm b_\mathrm H+\mathrm L_{\mathcal U(\gamma)}=\mathrm b_{\mathrm H,\star}$, where 
$\mathrm b_{\mathrm H,\star}$ is the Hochschild differential on the Hochschild chain complex of 
$\Big(C^\infty(V)[\![\hbar]\!],\star\Big)$ with values in itself. 
Furthermore, the degree $1$ pairing $\cap_\star$ takes the form
\[
D\cap_\star c=(a_0\star D(a_1,\dots,a_{|D|+1})|a_{|D|+2}|\cdots),\ D\in D_\mathrm{poly}(V)[\![\hbar]\!],\ 
c=(a_0|a_1|\cdots)\in C_{-\bullet}^{\rm poly}(V)[\![\hbar]\!]
\]
(up to a sign depending on $D$ and $c$).
All these structures have been introduced and discussed in~\cite{Sh,CR1}.

Then, Shoikhet's conjecture follows from our main result, together with the compatibility between cup products~\cite{MT}.
\begin{Thm}\label{t-shoik}
For $V=\mathbb R^d$ and for a MCE $\gamma$ of $T_\mathrm{poly}(V)[\![\hbar]\!]$ as in~\eqref{eq-MC-def}, 
the tangent quasi-isomorphisms $\mathcal U_{\gamma}$ and $\mathcal S_{\gamma}$ fit into the following 
commutative diagram of Gerstenhaber algebras and $T$-modules up to homotopy: 
\[
\xymatrix{(T_\mathrm{poly}^\bullet(V)[\![\hbar]\!],[\gamma,\ ],[\ ,\ ],\cup)\ar[d]\ar[r]^-{\mathcal U_\gamma} & 
(D_\mathrm{poly}^\bullet(V)[\![\hbar]\!],\mathrm d_{\mathrm H,\star},[\ ,\ ],\cup_\star)\ar[d] \\
(\Omega^{-\bullet}(V)[\![\hbar]\!],\mathrm L_\gamma,\mathrm L,\cap) & 
\ar[l]_-{\mathcal S_\gamma} (C_{-\bullet}^{\rm poly}(V)[\![\hbar]\!],\mathrm b_{\mathrm H,\star},\mathrm L,\cap_\star)}\,.
\]
\end{Thm}

\subsection{The case of a MCE of polyvector degree at most $1$}\label{ss-6-2}

We then consider a MCE of $T_\mathrm{poly}^\mathfrak n(V)$ of polyvector degree at most $1$, i.e.\
\begin{equation}\label{eq-MCel-1}
\gamma=\gamma_{-1}+\gamma_0+\gamma_1,
\end{equation}
where $i)$ $\gamma_{-1}$ is a degree $2$, $\mathfrak n$-valued function on $V$, $ii)$ $\gamma_0$ is an 
$\mathfrak n$-valued vector field of degree $1$ on $V$, and $iii)$ $\gamma_1$ is a degree $0$, $\mathfrak n$-valued bivector 
field on $V$. 

The image $\mathcal U(\gamma)$ of a MCE as in~\eqref{eq-MCel-1} satisfies the Maurer--Cartan equation in 
$\mathcal D_{\mathrm{poly}}^\mathfrak n(V)$: since $\gamma$ is the sum of three types of $\mathfrak n$-valued 
polyvector fields on $V$, the degree requirement of the classical morphisms $\mathcal U_n$ and the (graded) 
anticommutativity of the wedge product on $\mathfrak m$-valued polyvector fields implies the  
decomposition of $\mathcal U(\gamma)$
\begin{equation}\label{eq-MCtwist}
\begin{aligned}
\mathcal U(\gamma)&=\sum_{n\geq 1}\frac{1}{n!}\mathcal U_n(\underset{n}{\underbrace{\gamma_1,\dots,\gamma_1}})
+\sum_{n\geq 0}\frac{1}{n!}\mathcal U_{n+1}(\gamma_0,\underset{n}{\underbrace{\gamma_1,\dots,\gamma_1}})+\\
&\phantom{=}+\sum_{n\geq 0}\frac{1}{n!}\mathcal U_{n+1}(\gamma_{-1},\underset{n}{\underbrace{\gamma_1,\dots,\gamma_1}})
+\frac{1}2\left(\sum_{n\geq 0}\frac{1}{n!}\mathcal U_{n+2}(\gamma_0,\gamma_0,\underset{n}{\underbrace{\gamma_1,\dots,\gamma_1}})\right)\,.
\end{aligned}
\end{equation}
For the sake of simplicity, we write from now on $B$, resp.\ $Q$, resp.\ $F$, for the first term, resp.\ second 
term, resp.\ sum of the third and fourth term, in (\ref{eq-MCtwist}).

The graded commutative product on $A=C^{\infty}(V)\otimes\mathfrak m$ defines a $1$-cocycle $\mu$ of 
$\mathcal D_{\mathrm{poly}}^\mathfrak m(V)$. We may then consider the $\mathfrak m$-valued bidifferential 
operator $\mu+B$ of degree $0$ and the linear operator $\widetilde Q=\mathrm d_\mathfrak m+Q$ of degree $1$ on 
$\mathfrak m$-valued functions on $V$. Accordingly, the Maurer--Cartan equation for $\mathcal U(\gamma)$ is 
equivalent to 
\begin{enumerate}
\item[$i)$] $\mu+B$ defines an $\mathfrak m$-linear associative product $\star$ of degree $0$ on $A$.
\item[$ii)$] $\widetilde Q$ is a derivation of degree $1$ of $(A,\star)$; its square equals
\[
\widetilde Q^2=-[F,\ ]_\star,
\]
where $[\ ,\ ]_\star$ denotes the graded commutator w.r.t.\ the product $\star$. 
\item[$iii)$] The $\mathfrak m$-valued function $F$ of degree $2$, which, by the previous equation, can be viewed 
as a sort of ``curvature'' of the ``connection'' $\widetilde Q$, satisfies the Bianchi identity, i.e.\ it is 
annihilated by $\widetilde Q$: $\widetilde{Q}(F)=0$. 
\end{enumerate}
In other words, $A$ equipped with the product $\star$, the derivation $\widetilde{Q}$ and the element $F$, 
is a curved DGA; in the framework of~\cite{CFT}, $\widetilde Q$ is a Weyl connection on $A$ with Weyl curvature $F$, 
see Section~\ref{s-8}. 

For a MCE $\gamma$ as in~\eqref{eq-MCel-1}, we consider the $T$-algebra 
$\left(\big(T_{\mathrm{poly}}^\mathfrak m(V),\mathrm d_\mathfrak m+[\gamma,\ ],[\ ,\ ],\cup\big),
\big(\Omega^\mathfrak m (V),\mathrm L_\gamma,\mathrm L,\cap\big)\right)$.

We also consider the Gerstenhaber algebra up to homotopy 
$\left(\mathcal D_{\mathrm{poly}}^\mathfrak m(V),\mathrm d_\mathfrak m+\mathrm d_\mathrm H+[\mathcal U(\gamma),\ ],
[\ ,\ ],\cup_\star\right)$, where
\begin{equation}\label{eq-cupstar}
\left(D_1\cup_\star D_2\right)\!(a_1,\dots,a_n)
=D_1(a_1,\dots,a_{|D_1|+1})\star D_2(a_{|D_1|+2},\dots,a_n),\ n=|D_1|+|D_2|+2,
\end{equation}
for $D_i$, $i=1,2$, general elements of $D_\mathrm{poly}^\mathfrak m(V)$, and $a_j$, $j=1,\dots,|D_1|+|D_2|+2$, 
general elements of $A$. Additionally, we have 
$\mathrm d_\mathfrak m+\mathrm d_\mathrm H+[\mathcal U(\gamma),\ ]=\mathrm d_{\mathrm H,\star}+[\widetilde Q,\ ]+[F,\ ]$, 
where $\mathrm d_{\mathrm H,\star}$, denotes the Hochschild differential w.r.t.\ the product $\star$.
We further consider the $T$-module up to homotopy 
$\left(C_{-\bullet}^{{\rm poly},\mathfrak m}(V),\mathrm d_\mathfrak m+\mathrm b_\mathrm H+\mathrm L_{\mathcal U(\gamma)},\mathrm L,\cap_\star\right)$, 
where
\begin{equation}\label{eq-capstar}
D\cap_\star (a_0|a_1|\cdots|a_n)=(a_0\star D(a_1,\dots,a_{|D|+1})|a_{|D|+2}|\cdots|a_n),
\end{equation}
for $D$, resp.\ $c$, a general element of $D_\mathrm{poly}^\mathfrak m(V)$, 
resp.\ $C_{-\bullet}^{{\rm poly},\mathfrak m}(V)$. Furthermore, we have 
$\mathrm d_\mathfrak m+\mathrm b_\mathrm H+\mathrm L_{\mathcal U(\gamma)}
=\mathrm b_{\mathrm H,\star}+\mathrm L_{\widetilde Q}+\mathrm L_F$. 
\begin{Thm}\label{t-cupprod}
For any MCE $\gamma$ as in (\ref{eq-MCel}), $\mathcal U_\gamma$ and $\mathcal S_\gamma$ are quasi-isomorphisms 
of Gerstenhaber algebras and $T$-modules up to homotopy respectively, fitting into the following commutative diagram:
\[
\xymatrix{\left(T_{\mathrm{poly}}^\mathfrak m(V),\mathrm d_\mathfrak m+[\gamma,\ ],[\ ,\ ],\cup\right)
\ar[d]\ar[r]^-{\mathcal U_\gamma} & 
\left(\mathcal D_{\mathrm{poly}}^\mathfrak m(V),\mathrm d_{\mathrm H,\star}
+[\widetilde Q,\ ]+[F,\ ],[\ ,\ ],\cup_\star\right)
\ar[d] \\
\left(\Omega^\mathfrak m(V),\mathrm d_\mathfrak m+\mathrm L_\gamma,\mathrm L,\cap \right) & 
\ar[l]_-{\mathcal S_\gamma} 
\left(C_{-\bullet}^{{\rm poly},\mathfrak m}(V),\mathrm b_{\mathrm H,\star}
+\mathrm L_{\widetilde Q}+\mathrm L_F,\mathrm L,\cap_\star\right)}\,.
\] 
\end{Thm}

\subsection{Explicit computation of the tangent quasi-isomorphisms}\label{ss-6-3}

Borrowing notation from Subsection~\ref{ss-6-2}, we consider the case of a MCE $\gamma=\gamma_0$ 
concentrated in polyvector degree $0$. 

We consider now the morphism $\mathcal S_\gamma$: for a general Hochschild chain $c=(a_0|\cdots|a_m)$ of (Hochschild) degree $-m$, $m\geq 0$, we have
\begin{equation}\label{eq-wheeledsum}
\mathcal S_{\gamma}(c)=\sum_{n\geq 0}\frac{1}{n!}\sum_{\Gamma\in\mathcal G_{n,m+1}^S}W_{D,\Gamma}\mathcal S_{\Gamma}(\underset{n}{\underbrace{\gamma,\dots,\gamma}},c),
\end{equation}
with the notations from Subsection~\ref{ss-3-2}.

We first observe that the valence of any vertex of the first type of an admissible graph $\Gamma$ as in (\ref{eq-wheeledsum}) is $1$, since it associated to a copy of the vector field $\gamma$.
Thus, we sum only over those admissible graphs $\Gamma$ with univalent vertices of the first type.

Dimensional reasons imply that the weight (\ref{eq-weight-mod}) of an admissible graph $\Gamma$ in $\mathcal G_{n,m+1}^S$ is non-trivial, only if $2n+m$ equals the degree of the integrand, which is in this case $n+l$, where $l$ is the number of edges starting from the vertex $0$, whence $l=n+m$.
Since such an admissible graph has exactly $n+m+1$ vertices (of the first and second type), and since there are neither multiple edges nor loops by assumption, there is exactly one edge joining the vertex $0$ to all vertices except one, namely the first vertex of the second type w.r.t.\ the cyclic order: this is because the integrand $\omega_{D,\Gamma}$ vanishes, if $\mathrm{star}(0)$ contains an edge $e$ joining $0$ to the first vertex in $S^1$ w.r.t.\ the cyclic order, by the constructions of Subsubsection~\ref{sss-3-3-2}.

Since $m\geq 0$, we use the section of $D_{n,m+1}^+$, which, by means of the M\"obius transformations $\psi$ is diffeomorphic to a section of $C_{n+1,m}^+$, see also Subsubsection~\ref{sss-3-1-3} (we observe that the origin $0$ of the disk is mapped to $\mathrm i$, while the point $1$ is mapped to the half-circle at infinity in the complex upper half-plane $\mathcal H$).

Recalling (\ref{eq-Dprop}), Subsubsection~\ref{sss-3-2-1}, the weight (\ref{eq-weight-mod}) of an admissible graph $\Gamma$ in $\mathcal G_{n,m+1}^S$ is mapped to a weight of type (\ref{eq-weight}), the only difference is that the factors of $\omega_{D,\Gamma}$ are mapped to $i)$ usual forms $\omega_e$, whenever $e$ is an edge from $0$ to some vertex (of the first and of the second type), and the ``new'' $e$ is now an edge from $\mathrm i$ to the image w.r.t.\ $\psi$ of the endpoint, and $ii)$ differences between $\omega_e$ and $\omega_{e(\mathrm i)}$, if $e$ is an edge between two vertices (of the first and second type, neither of which is $0$), and the new edge $e$ connects the images of the endpoints w.r.t.\ $\psi$, while $e(\mathrm i)$ is an edge, whose starting point is the starting point of $e$ and whose endpoint is $\mathrm i$.    
Graphically, we have the correspondence
\bigskip
\begin{center}
  \includegraphics[scale=0.28]{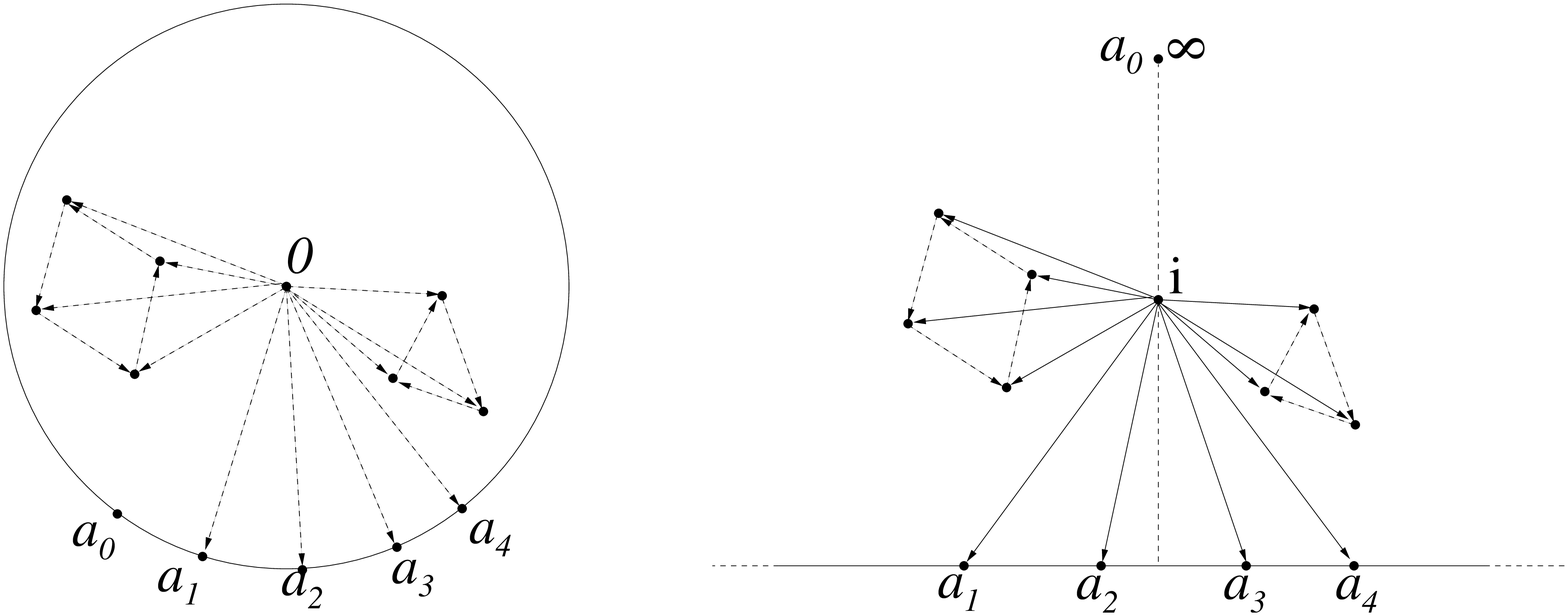} \\
\text{Figure 24 - Correspondence between weights on $\mathcal D_{n,m+1}^+$ and on $\mathcal C_{n+1,m}^+$ } \\
\end{center}
\bigskip
We have used dashed arrows to denote forms on $\mathcal D_{n,m+1}^+$ as in (\ref{eq-weight-mod}), in the graph on the left-hand side, while we have used black, resp.\ dashed, arrows to denote forms on $\mathcal C_{n+1,m}^+$ as in (\ref{eq-weight}), resp.\ differences of such forms. 

From now on, when considering weights (\ref{eq-weight-mod}) of admissible graphs in $\mathcal G_{n,m+1}^S$ as above, we implicitly assume that we are considering them on $\mathcal C_{n+1,m}^+$ by means of the previous correspondence.
\begin{Lem}\label{l-K-var}
If an admissible graph $\Gamma$ in $\mathcal G_{n,m+1}^S$ as above has a vertex of the first type of valence $1$, which is the endpoint of exactly one edge, then its weight vanishes.
\end{Lem}
\begin{proof}
It is more instructive to give a graphical proof (which will be also useful later on)
\bigskip
\begin{center}
  \includegraphics[scale=0.34]{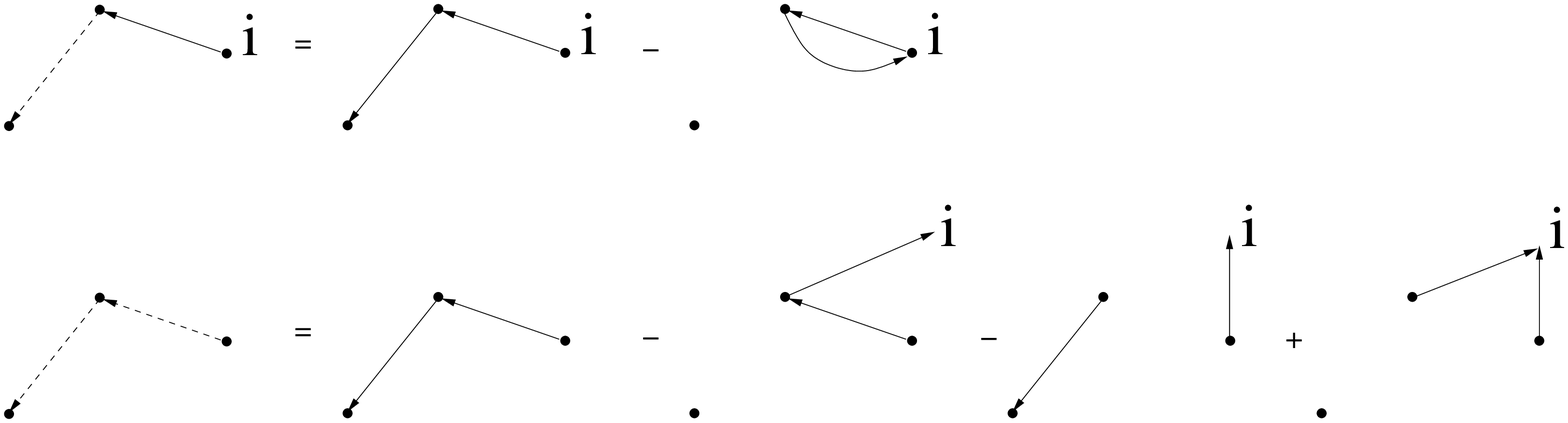} \\
\text{Figure 25 - The two cases in the above situation } \\
\end{center}
\bigskip
On the right-hand side of both equalities above, we may apply Lemmata~\ref{l-K-2} and~\ref{l-K-3}, Appendix; we also notice that, by dimensional reasons, whenever there is a $1$-valent vertex, the corresponding weight vanishes (since we integrate a $1$-form over a $2$-dimensional space, a subset of the complex upper half-plane $\mathcal H$).
\end{proof}
We consider now an admissible graph $\Gamma$ in $\mathcal G_{n,m+1}^S$, satisfying the above dimensional non-triviality condition.
We consider a vertex $v_1$ of the first type: it is the endpoint of an edge starting at $\mathfrak i$, and exactly one edge departs from it.
Moreover, the edge $e_1$ starting at $v$ must connect it to a different vertex of the first type: if not, $e_1$ connects $v$ to a vertex of the second type.
By Lemma~\ref{l-K-var}, there must be an edge $e_2$ from a vertex $v_2$ of the second type with endpoint $v_1$; again by Lemma~\ref{l-K-var}, there must be an edge $e_3$ from a vertex $v_3$ of the second type with endpoint $v_2$, and so on, until we arrive at the vertex $v_n$, which is necessarily as in Lemma~\ref{l-K-var} by dimensional reasons, whence the weight vanishes.

By the very same procedure, we find that all admissible graphs appearing in (\ref{eq-wheeledsum}) and having possibly non-trivial weights must be as in Figure 24 on the left-hand side, i.e.\ they must be wheeled trees (with dashed and black directed edges).
Using the first graphical identity in the proof of Lemma~\ref{l-K-var}, we may replace all dashed arrows by black ones: in fact, 
\bigskip
\begin{center}
  \includegraphics[scale=0.34]{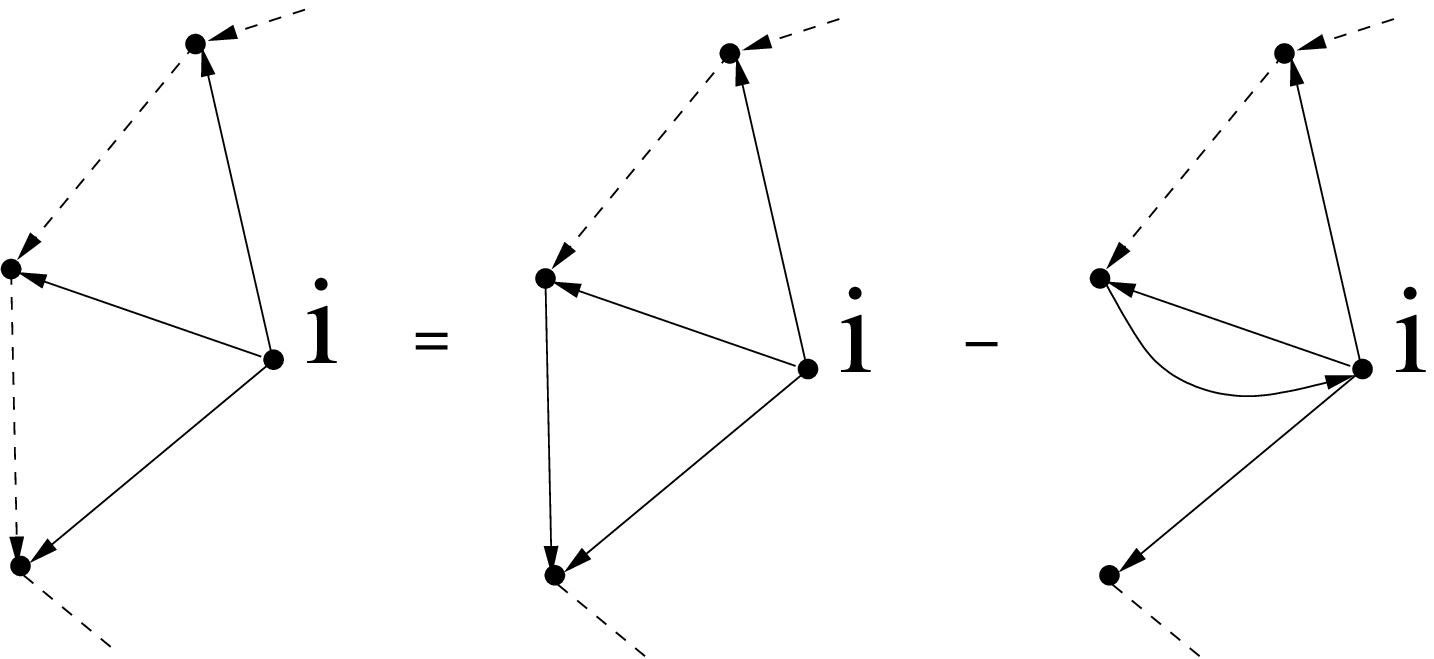} \\
\text{Figure 26 - Replacing dashed edges by black ones } \\
\end{center}
\bigskip 
Any wheeled tree with dashed edges, by repeatedly applying this computation, can be written as the sum of the same wheeled tree and other graphs with only black edges; except the wheeled tree with only black edges, each of the remaining graphs has {\bf at least} one vertex of the first type as in Lemmata~\ref{l-K-2} or~\ref{l-K-3}, whence they vanish.
 
We denote by $\mathcal T_{n,m}^K\subset \mathcal G_{n,m}^K$ the set of wheeled trees as above, which we view as admissible graphs as in Subsubsection~\ref{ss-3-1}: then, by the previous computations, we may re-write (\ref{eq-wheeledsum}) as
\begin{equation}\label{eq-only-wh}
\mathcal S_{\gamma}(c)=\sum_{n\geq 0}\frac{1}{n!}\sum_{\Gamma\in\mathcal T_{n,m+1}^K}W_\Gamma \mathcal S_{\Gamma}(\underset{n}{\underbrace{\gamma,\dots,\gamma}},c).
\end{equation}
We consider a wheeled tree $\Gamma$ in $\mathcal T_{n,m}^K$: by the same arguments as in~\cite{K}, Paragraph 8.3.3.1, if $\Gamma$ contains at least a wheel with an odd number of vertices, then its weight vanishes.
Thus, in (\ref{eq-only-wh}), we may sum only w.r.t.\ even integers $n$.

From now on, we may follow the arguments of Subsection 10.1 of~\cite{CVdB} to evaluate (\ref{eq-only-wh}): of course, there are some sign modifications to keep into account, but the end result turns out to be the same (since $n$ is even and since only wheels with an even number of vertices appear on the above sum).
First, we write 
\[
\gamma=\gamma_{\alpha}\otimes m_\alpha,\ \gamma_{\alpha}\in T_{\mathrm{poly}}^0(V),\ m_\alpha\in\mathfrak m_1,
\]
and, following the notations of~\cite{CVdB}, we introduce the $\mathfrak m$-valued, matrix-valued $1$-form $\Xi$ via
\[
\Xi=\Xi_\alpha\otimes m_\alpha=(\partial_i\partial_k \gamma_{\alpha}^j \mathrm d x_k)\otimes m_\alpha.
\]
Then, following almost verbatim the computations of Subsection 10.1 of~\cite{CVdB}, we find 
\begin{equation}\label{eq-todd}
\mathcal S_{\gamma}(c)=\det\!\sqrt{\frac{\Xi}{\mathrm e^{\frac{\Xi}2}-\mathrm e^{-\frac{\Xi}2}}}\wedge \mathrm{HKR}(c)=j(\gamma)\wedge \mathrm{HKR}(c),
\end{equation}
The right-hand side of (\ref{eq-todd}) needs some explanations.
First of all, we have improperly written a determinant: in fact, it should be denoted by the more appropriate notation $\mathrm{Ber}$, which represents the super-determinant, or Berezinian.
In fact, $\Xi$ represents a $1$-form on $V$ with values in $\mathfrak m$-valued matrices: $\mathfrak m$-valued matrices form a GA, hence usual trace and determinant have to be replaced by their super-analoga.

Further, the square root of the quotient in the Berezinian has to be interpreted as a power series.
More precisely, we have
\[
\frac{1}2 \log\!\left(\frac{e^{\frac{x}2}-e^{-\frac{x}2}}x\right)=\sum_{l\geq 0} \beta_l x^l,
\] 
and the coefficients $\beta_l$ are called {\bf modified Bernoulli numbers}.
 
Aside from some sign differences, which, as already remarked, do not cause changes in the main arguments, the only point we want to stress is that, by the above arguments, the weights in (\ref{eq-only-wh}) are the same weights examined in~\cite{CVdB}, whose computation has been performed, using different approaches, in~\cite{VdB,W}.

Summarizing all the computations so far, we have the following 
\begin{Thm}\label{t-wheel}
For a MCE of $T_\mathrm{poly}^\mathfrak m(V)$ of polyvector degree $0$, the following identity holds true:
\begin{equation}\label{eq-wheel}
\mathcal S_{\gamma}(c)=j(\gamma)\wedge\mathrm{HKR}(c),\ c\in C^{\mathrm{poly},\mathfrak m}_{-\bullet}(V),
\end{equation}
where $\mathrm{HKR}$ is the Hochschild--Kostant--Rosenberg quasi-isomorphism in homology, and $j(\gamma)$ 
is the rooted Todd class analogon appearing in the main result of~\cite{CVdB}.
\end{Thm}

\section{Application : (co)homological Duflo isomorphism}\label{s-7}

We consider a finite dimensional Lie algebra $\mathfrak g$ over a field $k$ of characteristic zero. 

\subsection{Statement of the result}\label{ss-7-1}

We recall the definition of the {\bf (modified) Duflo element} 
$$
J:=\det\Big(\frac{e^{\mathrm{ad}/2}-e^{-\mathrm{ad}/2}}{\mathrm{ad}}\Big)\in\widehat{\S}(\mg^*)^{\mg}\,.
$$
We also remind the reader that the completed algebra $\widehat{\S}(\mg^*)$ naturally 
acts on $\S(\mg)$: 
$$
\xi^k\cdot x^n:=\frac{n!}{(n-k)!}\xi(x)^kx^{n-k}\qquad(x\in\mg\,,\,\xi\in\mg^*\,,\,k>0\,,\,n>0)\,.
$$

The following result (proved in \cite{K,PT}) is a cohomological extension of the original 
Duflo isomorphism \cite{Du}. 
\begin{Thm}[Cohomological Duflo isomorphism]\label{thm-cdi}
The morphism of $\mg$-modules 
$$
\mathcal D\,:=\,{\rm sym}\circ(J^{1/2}\cdot):\S(\mg)\,\longrightarrow\,\U(\mg)
$$
induces an algebra isomorphism 
$$
H^\bul(\mg,\S(\mg))\,\tilde\longrightarrow\,H^\bul(\mg,\U(\mg))
$$
at the level of Chevalley-Eilenberg cohomology. 
\end{Thm}
We now observe that, if $A$ is an algebra on which $\mg$ acts by derivations, the Chevalley--Eilenberg Lie algebra homology $H_{-\bul}(\mg,A)$ is equipped with 
an $H^\bul(\mg,A)$-module structure in the following way: on the level of the complexes, for any Chevalley-Eilenberg 
cochain $\alpha=\xi\otimes a$, resp.~chain $c=x\otimes a'$, one defines
$$
\alpha(c)=\iota_\xi(x)\otimes aa'\,,
$$
where $\iota$ denotes the usual contraction operation\footnote{In fact, this defines an actual DG-module structure 
on the level of the complexes. }. 
In what follows we will prove the following homological version of the Duflo isomorphism. 
\begin{Thm}[Homological Duflo isomorphism]\label{thm-hdi}
The morphism $\mathcal D$ induces an isomorphism of $H^\bul\big(\mg,\S(\mg)\big)$-modules 
$$
H_{-\bul}\big(\mg,\S(\mg)\big)\,\tilde\longrightarrow\,H_{-\bul}\big(\mg,\U(\mg)\big)
$$
at the level of Chevalley-Eilenberg homology. 
\end{Thm}

Considering the degree zero (co)homology, one obtains
\begin{Cor}\label{cor-chdi}
$\mathcal D$ restricts to an isomorphism of algebras 
$\S(\mg)^\mg\,\tilde\longrightarrow\,\U(\mg)^\mg=\mathcal Z\big(\U(\mg)\big)$ on invariants, 
and induces an isomorphism of $\S(\mg)^\mg$-modules 
$\S(\mg)_\mg\,\tilde\longrightarrow\,\U(\mg)_\mg=\mathcal{A}\big(\U(\mg)\big)$ on coinvariants.  
\end{Cor}
Here $\mathcal{Z}(B)$ denotes the center of an algebra $B$, and $\mathcal A(B)=B/[B,B]$ its abelianization. 

\subsection{Proof of the results}\label{ss-7-2}

\subsubsection{Proof of Theorem \ref{thm-cdi}}\label{sss-7-2-1}

In this Subsubsection we follow closely \cite[Subsection 5.2]{CR}. 

\medskip

Let us consider the superspace $V=\Pi\mg$: then, the GA $A$ of superfunctions on $V$ is 
$A=\wedge^\bul(\mg^*)$. 
Therefore the Chevalley-Eilenberg differential $\d_\C$ defines a cohomological 
vector field $\gamma$ on $V$ (i.e.~a degree one derivation squaring to $0$)\footnote{This is abstract nonsense. }. 
In other words, $\gamma$ is a MCE in $T_{\rm poly}^\bul(V)$. 

On the one hand, $T_{\rm poly}^\bul(V)$ is naturally isomorphic to $\wedge^\bul(\mg^*)\otimes\S(\mg)$ and, 
under this identification, $[\gamma,\ ]$ precisely gives the coboundary operator $\d_\C$ of the 
Chevalley-Eilenberg cochain complex of $\mg$ with values in $\S(\mg)$. 
On the other hand, $(D_{\rm poly}^\bul(V),\d_\H+[\gamma,\ ])$ identifies with the complex $CC^\bullet(A,\mathrm d_C)$ of 
Hochschild cochains of the DGA $(A,\d_\C)$ with values in itself. 

Now we observe that we have a quasi-isomorphism of DGAs 
$$
\ell\,:\,(D_{\rm poly}^\bul(V),\d_\H+[\gamma,\ ])\,\longrightarrow\,C^\bul\big(\mg,\U(\mg)\big)\,,
$$
where $C^\bul\big(\mg,\U(\mg)\big)$ denotes the Chevalley-Eilenberg cochain complex of $\mg$ with values in $\U(\mg)$, 
given by the following composition of maps 
$$
D_{\rm poly}^\bul(V)=\wedge^\bul(\mg^*)\otimes{\rm T}\big(\wedge^\bul(\mg)\big)
\twoheadrightarrow\wedge^\bul(\mg^*)\otimes{\rm T}(\mg)\twoheadrightarrow
\wedge^\bul(\mg^*)\otimes\U(\mg)=C^\bul\big(\mg,\U(\mg)\big)\,.
$$
This is a manifestation of the fact that the quadratic DGA 
$\big(\wedge^\bul(\mg^*),\d_\C\big)$ and the quadratic-linear algebra $\U(\mg)$ are 
related by a Koszul-type duality (see e.g.~\cite{P}). 
Moreover, the following diagram of quasi-isomorphisms of complexes commutes : 
$$
\xymatrix{
\Big(T_{{\rm poly}}^\bul(V),[\gamma,\ ]\Big) \ar@{=}[d]\ar[rr]^{{\rm HKR}} & &
\Big(D_{\mathrm{poly}}^\bul(V),\d_\H+[\gamma,\ ]\Big) \ar[d]^{\ell} \\
C^\bul(\mg,\S(\mg))\ar[rr]^{{\rm sym}} & &
C^\bul(\mg,\U\mathfrak(g))\,,
}
$$

\medskip

Finally, we recall (see e.g.~\cite{CR}) that $\mathcal U_{\gamma}={\rm HKR}\circ\iota_{j(\gamma)}$, and that 
one has the following 
\begin{Lem}[\cite{CR}, Lemma 5.6]
Under the obvious identification $\Pi V=\mg$, the supermatrix valued $1$-form $\Xi$, 
restricted to $\mg$, which we implicitly identify with the space of vector fields on $V$ with 
constant coefficients, satisfies $\Xi=\mathrm{ad}$.
\end{Lem}
Therefore, since 
$$
j(\gamma):=\det\sqrt{\frac{e^{\Xi/2}-e^{-\Xi/2}}{\Xi}}\,,
$$
we have the following commutative diagram of algebra isomorphisms 
$$
\xymatrix{
H^\bul\Big(T_{{\rm poly}}^\bul(V),[\gamma,\ ]\Big) \ar@{=}[d]\ar[rr]^{\mathcal U_{\gamma}} & &
H^\bul\Big(D_{\mathrm{poly}}^\bul(V),\d_\H+[\gamma,\ ]\Big) \ar[d]^{\ell} \\
H^\bul(\mg,\S(\mg))\ar[rr]^{\mathcal D} & &
H^\bul(\mg,\U\mathfrak(g))\,.
}
$$
Hence Theorem~\ref{thm-cdi} follows. \hfill \qed

\subsubsection{Proof of Theorem \ref{thm-hdi}}

We keep the same notations as in the previous paragraph. 

\medskip

First of all, we observe that $\Omega^\bul(V)$ is naturally isomorphic to $\wedge^\bul(\mg^*)\otimes\S(\mg^*)$ 
and that, under this identification, $\mathrm{L}_{\gamma}$ precisely gives the coboundary operator of the 
Chevalley-Eilenberg cochain complex of $\mg$ with values in $\S(\mg^*)$. 

Then, $\big(C^{\rm poly}_\bullet(V),{\rm b}_{\rm H}+\mathrm{L}_{\gamma}\big)$ identifies with the complex 
$CC_{-\bullet}(A,\rm d_C)$ 
of Hochschild chains (with reversed grading) of the DGA $(A,\d_\C)$ with values in itself. 

Applying Theorem~\ref{t-wheel} to the present situation, we have $\mathcal S_{\gamma}=j(\gamma)\wedge{\rm HKR}$. 
We observe that, in order for this map to be well-defined, we need to consider completed versions 
$\widehat{\Omega}^\bul(V)=\wedge^\bul (\mg^*)\otimes\widehat{\S}(\mg^*)$ and 
$\widehat{C}_{-\bul}^{\rm poly}(V)=\wedge^\bul(\mg^*)\otimes\widehat{\rm T}\big(\wedge^\bul(\mg^*)\big)$ 
of the spaces involved in the formality for chains. 

Now, we recall that we have a quasi-isomorphism of complexes 
$$
\lambda\,:\,C^\bul(\mg,\U(\mg)^*)\longrightarrow\Big(\widehat{C}_{-\bul}^{\rm poly}(V),{\rm b}+{\rm L}_{\gamma}\Big)
$$
given by the following composition of maps 
$$
C^\bul\big(\mg,\U(\mg)^*\big)=\wedge^\bul(\mg^*)\otimes\U(\mg)^*\hookrightarrow
\wedge^\bul(\mg^*)\otimes{\rm T}(\mg)^*=\wedge^\bul(\mg^*)\otimes\widehat{{\rm T}}(\mg^*)
\hookrightarrow\wedge^\bul(\mg^*)\otimes\widehat{{\rm T}}\big(\wedge^\bul(\mg^*)\big)\,,
$$
which induces an isomorphism of $H^\bul\big(\mg,\U(\mg)\big)$-modules on cohomology. 

Moreover, the following diagram of quasi-isomorphisms of complexes commutes : 
$$
\xymatrix{
\left(\widehat\Omega(V),{\rm L}_{\gamma}\right) \ar@{=}[d] & &
\left(\widehat{C}_{-\bul}^{\rm poly}(V),{\rm b}+{\rm L}_{\gamma}\right) \ar[ll]_{{\rm HKR}}  \\
C^\bul(\mg,\widehat{\S}(\mg^*)) & &
C^\bul(\mg,\U\mathfrak(\mg)^*) \ar[u]_{\lambda} \ar[ll]_{{\rm sym}^*}\,,
}
$$
We observe that, for any $\mg$-module $M$, the Chevalley-Eilenberg cochain complex $C^\bul(\mg,M^*)$ 
is naturally isomorphic to the dual of the Chevalley-Eilenberg chain complex $C_{-\bul}(\mg,M)$ with reversed grading. 
Moreover, a direct computation shows that 
\begin{Lem} For any $\omega\in\widehat{\Omega}^\bul(V)=C^\bul\big(\mg,\widehat{\S}\mathfrak(\mg)^*\big)$ and any 
$c\in C_{-\bul}(\mg,S(\mg))$, \\
\indent i) $\langle j(\gamma)\wedge\omega,c\rangle=\langle\omega,{j(\gamma)}\cdot c\rangle$; \\
\indent ii) for any $\alpha\in T_{\rm poly}^\bul(V)=C^\bul(\mg,S(\mg))$, 
$\langle\iota_\alpha\omega,c\rangle=\langle\omega,\alpha(c)\rangle$. 
\end{Lem}
Therefore the transpose of $\mathcal{S}_{\gamma}$ induces an isomorphism of $H^\bul(\mg,S(\mg))$-modules 
$$
H_{-\bul}(\mg,\S(\mg))\longrightarrow H_{-\bul}(\mg,\U\mathfrak(g))
$$
which is precisely $\mathcal{D}$, whence the proof of Theorem \ref{thm-hdi}. \hfill \qed
\begin{Rem}
As we already mentioned, there is a duality between the DGA $(A,\d_\C)$ and the quadratic-linear algebra 
$\U(\mg)$: in~\cite{CR1}, we give a more direct proof of Corollary~\ref{cor-chdi} in the same spirit of Kontsevich's approach to the original Duflo isomorphism~\cite{K}, which does not make use of the aforementioned duality. 
\end{Rem}

\subsection{Why we can work over $\mathbb{Q}$}\label{ss-7-4}

In this Subsection we explain why Theorem \ref{t-wheel}, Theorem \ref{thm-cdi}, Theorem \ref{thm-hdi} 
and Corollary \ref{cor-chdi} are valid over any field of zero characteristic. 

First of all, we observe that we have been able to compute explicitly $\mathcal{U}_{\gamma}$ (Section 9 of \cite{CR}, 
see also \cite{CVdB}) and $\mathcal{S}_{\gamma}$ (Section \ref{s-6} of the present paper), and both have rational coefficients. 

Then to prove that the mentioned results remain true over $\mathbb{Q}$ (and thus over any field of zero characteristic) 
we have to find homotopies with rational coefficients. 

Finally, we observe that both homotopy equations~\eqref{eq-homotopy-U} and~\eqref{eq-homotopy1} are linear 
w.r.t.\ the weights of graphs appearing in the homotopy operator $\mathcal{H}_\gamma$, respectively. 

To conclude, we have a real solution of a system of linear equations with rational coefficients. 
Therefore a rational solution exists. \hfill\qed

\section{Proof of the main result}\label{s-8}

The main goal of this final Section is to globalize to a general manifold $X$ the local results obtained 
above in the paper. The globalisation procedure is based on~\cite{Dol} (see also \cite{CFT}). 

\subsection{Fedosov resolutions and the globalisation procedure}\label{ss-8-1}

We consider a general $d$-dimensional manifold $X$.
According to Section 3,~\cite{Dol}, we consider the algebra $B=C^\infty(X)$ of smooth functions on $X$; we 
associate to $X$ the DGLAs $(T_\mathrm{poly}(X),0,[\ ,\ ])$ of polyvector fields on $X$ with trivial differential 
and Schouten--Nijenhuis bracket, and $\left(D_\mathrm{poly}(X),\mathrm d_\mathrm H,[\ ,\ ]\right)$ of polydifferential operators on $X$, which is viewed as the subcomplex of the Hochschild cochain complex of $B$, consisting of cochains, which are smooth differential operators w.r.t.\ any of their arguments, and with induced Hochschild differential and 
Gerstenhaber bracket. Both DGLAs are graded w.r.t.\ the shifted degree. 

Furthermore, we have corresponding DGMs, $\left(\Omega(X),0,\mathrm L\right)$, with trivial differential and action 
$\mathrm L$ by polyvector fields given by Lie derivative, and 
$\left(C^\mathrm{poly}(X),\mathrm b_\mathrm H,\mathrm L\right)$, where $C^\mathrm{poly}(X)$ has been defined 
in~\cite{Dol}: it is defined as a suitable completion of the Hochschild chain complex of $C^\infty(X)$ with negative grading, with the Hochschild differential $\mathrm b_\mathrm H$ and the action $\mathrm L$.
We observe that both DGMs are negatively graded.

We quote (without proof) from~\cite{Dol} the main result towards the globalisation of Kontsevich's and Tsygan's 
formality $L_\infty$-quasi-isomorphism.
\begin{Thm}\label{t-glob}
For a given $d$-dimensional manifold $X$, there exist DGLAs $\mathfrak g_i^X$, resp.\ DGMs 
$\mathfrak M_i^X$, $i=1,2$, and $L_\infty$-quasi-isomorphisms $\mathfrak U_X$ and $\mathfrak S_X$, which fit into the following commutative diagram of DGLAs and DGMs:
\begin{equation}\label{eq-glob-CD}
\xymatrix{T_\mathrm{poly}(X)\ar@{^{(}->}[r]\ar[d] & 
\mathfrak g_1^X\ar[d]\ar[r]^-{\mathfrak U_X} & 
\mathfrak g_2^X \ar[d] & 
D_\mathrm{poly}(X)\ar@{_{(}->}[l]\ar[d] \\
\Omega(X)\ar@{^{(}->}[r] & 
\mathfrak M_1^X & 
\ar[l]_-{\mathfrak S_X} \mathfrak M_2^X 
& C^\mathrm{poly}(X)\ar@{_{(}->}[l]}\,,
\end{equation} 
where the vertical arrows denote DGLA-actions, and the hooked arrows denote quasi-isomorphisms.
\end{Thm}
Here we briefly discuss the construction of the DGLAs $\mathfrak g_i^X$ and corresponding DGMs 
$\mathfrak M_i^X$, $i=1,2$, which are called Fedosov resolutions of the corresponding DGLAs and DGMs, 
which have been defined above. Using the notations of~\cite{Dol}, we have
\[
\begin{aligned}
\mathfrak g_1^X&=\Omega(X,\mathcal T_\mathrm{poly}),\ &\ \mathfrak g_2^X&=\Omega(X,\mathcal D_{\mathrm{poly}}),\\
\mathfrak M_1^X&=\Omega(X,\mathcal E),\ & \ \mathfrak M_2^X&=\Omega(X,\mathcal C^\mathrm{poly}),
\end{aligned}
\]
where $\mathcal T_\mathrm{poly}$, resp.\ $\mathcal D_\mathrm{poly}$, denotes the bundle over $X$, whose global sections 
are formal polyvector fields, resp.\ formal polydifferential operators, w.r.t.\ fiber coordinates of $TX$, and 
$\mathcal E$, resp.\ $\mathcal C^\mathrm{poly}$, denotes the bundle, whose global sections are formal differential 
forms, resp.\ formal Hochschild chains, w.r.t.\ fiber coordinates of $TX$.
$\mathcal T_\mathrm{poly}$ and $\mathcal D_\mathrm{poly}$ are bundles of DGLAs, the former with trivial differential and 
$B$-linear Schouten--Nijenhuis bracket w.r.t.\ the formal coordinates of $TX$, the latter with Hochschild differential 
and Gerstenhaber bracket induced from the DGLA-structure on the Hochschild cochain complex of the algebra 
$F=\mathbb R[\![y_1,\dots,y_d]\!]$; similarly, $\mathcal E$ and $\mathcal C^\mathrm{poly}$ are bundles of DGMs over 
$\mathcal T_\mathrm{poly}$ and $\mathcal D_\mathrm{poly}$ respectively, the former with trivial differential and Lie 
derivative w.r.t.\ the fiber coordinates of $TX$, the latter with Hochschild differential and action $\mathrm L$ 
induced from the DGM-structure on the Hochschild chain complex of $F$. 

Obviously, the De Rham differential $\mathrm d$ on $X$ defines, by linear extension, a differential (which we denote by 
the same symbol) on the DGLAs $\mathfrak g_i^X$ and on the DGMs $\mathfrak M_i^X$, $i=1,2$, which is compatible with 
all aforementioned algebraic structures; it is then clear that all DGLAs and DGMs considered so far are naturally 
bi-graded.

A very important tool in the proof of Theorem~\ref{t-glob} is the Fedosov connection on $\mathfrak g_i^X$ and 
$\mathfrak M_i^X$, $i=1,2$: it is customary to denote it as $D=\mathrm d+\omega$, where $\omega$ decomposes as $\omega=A+\widetilde\omega$, 
where $A$ is the connection $1$-form of a torsion-free connection $\nabla=\mathrm d+A$ on $TX$ and $\widetilde \omega$ is 
an element of $\Omega^1(X,\mathcal T_\mathrm{poly}^0)$. 
The Fedosov connection is flat, i.e.\ $D^2=0$, or, equivalently, 
${\rm d}\omega+\frac{1}2[\omega,\omega]=0$; we observe that, since $\nabla$ is torsion-free, then 
$\nabla$ extends to a derivation of degree $1$ on $\mathfrak g_i^X$ and $\mathfrak M_i^X$, $i=1,2$. 
Since $D$ is flat and is compatible with all algebraic structures, we may consider the cohomology of $\mathfrak g_i^X$ 
and $\mathfrak M_i^X$, $i=1,2$ w.r.t.\ $D$: it turns out that all cohomologies are concentrated in degree $0$, and that 
we have isomorphisms of DGLAs and DGMs 
\begin{equation}\label{eq-FedoReso}
\begin{aligned}
\mathrm H^0(\Omega(X,\mathcal T_\mathrm{poly}),D) & \cong T_\mathrm{poly}(X)\,,\ & 
\ \mathrm H^0(\Omega(X,\mathcal D_\mathrm{poly}),D) & \cong D_\mathrm{poly}(X)\,, \\
\mathrm H^0(\Omega(X,\mathcal E),D)&\cong \Omega(X)\,,\ & 
\ \mathrm H^0(\Omega(X,\mathcal C^\mathrm{poly}),D) & \cong C^\mathrm{poly}(X)\,.
\end{aligned}
\end{equation}

\medskip

We finally briefly present the construction of the $L_\infty$-quasi-isomorphisms $\mathfrak U_X$ and $\mathfrak S_X$.

We pick a local chart $U$ of $X$, and we set $(\mathfrak m,\mathrm d_\mathfrak m)=(\Omega(U),\mathrm d)$, 
$\mathrm d$ being the De Rham differential: we observe that $\Omega(U)$ is commutative (in the graded sense), 
and has a decomposition into a nilpotent part $\mathfrak n=\Omega^{\geq1}(U)$ and a commutative algebra 
(concentrated in degree $0$), which contains a unit annihilated by $\mathrm d$. 

On the other hand, Kontsevich's and Tsygan's $L_\infty$-quasi-isomorphisms $\mathcal U$ and $\mathcal S$ extend to 
$L_\infty$-quasi-isomorphisms on formal polyvector fields, polydifferential operators, differential forms and 
Hochschild chains on $V=\mathbb R^d$; we denote by $V_\mathrm{formal}$ the formal linear manifold, whose algebra 
of functions is $\mathcal O_{V_\mathrm{formal}}=\mathbb R[\![y_1,\dots,y_d]\!]$.
This way, we have the following identifications of DGLAs and DGMs
\[
\begin{aligned}
\mathfrak g_1^U&\cong T_\mathrm{poly}^\mathfrak m(V_\mathrm{formal}),\ &\ 
\mathfrak g_2^U&\cong D_\mathrm{poly}^\mathfrak m(V_\mathrm{formal})\,, \\
\mathfrak M_1^U&\cong \Omega^\mathfrak m(V_\mathrm{formal}),\ &\ 
\mathfrak M_2^U&\cong C^{{\rm poly},\mathfrak m}(V_\mathrm{formal})\,,
\end{aligned}
\]
and a corresponding commutative diagram of DGLAs and DGMs. 

The connection $1$-form $\omega|_U$ of the Fedosov connection $D$, when restricted on $U$, is a MCE of polyvector 
degree $0$ in $T_\mathrm{poly}^\mathfrak n(V_\mathrm{formal})$.
We observe that, picking a distinct local chart $V$ of $X$ intersecting $U$ non-trivially, the difference between 
$\omega|_U$ and $\omega|_V$ is a $1$-form on the intersection $U\cap V$ with values in the linear vector fields on $F$. 

Then, the $L_\infty$-quasi-isomorphisms $\mathfrak U_X$ and $\mathfrak S_X$ are constructed by gluing the local 
$L_\infty$-morphisms $\mathfrak U_U$ and $\mathfrak S_U$, obtained by twisting the natural extensions of $\mathcal U$ 
and $\mathcal S$ w.r.t.\ the MCE $B|_U$; we only observe that the additional properties of $\mathcal U$ and 
$\mathcal S$, proved in~\cite{K,Dol}, imply that a change of local charts does not affect $\mathfrak U_U$ and 
$\mathfrak S_U$, whence it follows that they glue together.

\subsection{End of the proof of Theorem~\ref{t-glob}}\label{ss-8-1bis}

Let $(\mathfrak m,{\rm d}_{\mathfrak m})$ be a commutative DGA as in the introduction, 
and consider a MCE $\gamma$ in $T_\mathrm{poly}^{\mathfrak n}(X)$. 
First of all, we observe that the commutative diagram \eqref{eq-glob-CD} can be extended 
by $\mathfrak m$-linearity, and we denote by $\mathcal U$ (resp.~$\mathcal S$) the $L_\infty$-quasi-isomorphism 
obtained by means of a quasi-inverse $\mathcal I$, resp.~$\mathcal J$ (both not unique), of the right-most, resp.~left-most, hooked arrow of the first, resp.~second, line of \eqref{eq-glob-CD}. 

Then, the quasi-isomorphism of DGLAs $T_{\rm poly}^{\mathfrak m}(X)\hookrightarrow\mathfrak{g}_1^{X,\mathfrak m}$ 
produces, out of $\gamma$, a MCE $\gamma_1$ in $\mathfrak{g}_1^{X,\mathfrak n}$, which satisfies by construction, 
and by means of \eqref{eq-FedoReso}, $D(\gamma_1)=0={\rm d}_{\mathfrak m}(\gamma_1)+\frac12[\gamma_1,\gamma_1]$. 

Similarly, the $L_\infty$-quasi-isomorphism $\mathfrak U_X$ produces, out of $\gamma_1$, a MCE 
$\gamma_2$ in $\mathfrak{g}_2^{X,\mathfrak n}$, which itself gives rise to a MCE $\gamma'$ obtained 
by means of the quasi-inverse $\mathcal I$ of $D_{\rm poly}^{\mathfrak m}(X)\hookrightarrow\mathfrak{g}_2^{X,\mathfrak m}$. 

\medskip

By construction, $\gamma'=\mathcal U(\gamma)$, and $\mathcal U_\gamma$ can be computed as the composed 
$L_\infty$-morphism 
\begin{equation}\label{eq-UUU}
\xymatrix{T_\mathrm{poly}^{\mathfrak m}(X)_\gamma \ar@{^{(}->}[r] & 
\mathfrak g_{1,\gamma_1}^{X,\mathfrak m} \ar[r]^-{\mathfrak U_{X,\gamma_1}} & 
\mathfrak g_{2,\gamma_2}^{X,\mathfrak m} \ar[r]^-{\mathcal I_{\gamma_2}} & 
D_\mathrm{poly}^{\mathfrak m}(X)_{\gamma'}\,.}
\end{equation}
It follows from \cite{Dol,CDH} that hooked arrows in \eqref{eq-glob-CD} preserve all algebraic structures 
(namely, $B_\infty$-structures). 
Moreover, the compatibility between cup products from Section 
\ref{s-4} leads us to the following 
\begin{Lem}\label{lem-UUU}
The first Taylor component $\mathfrak U_{X,\gamma_1}^1$ of $\mathfrak U_{X,\gamma_1}$ induces an isomorphism 
of Gerstenhaber algebras 
$$
{\rm H}^\bullet\Big(\mathfrak g_{1,\gamma_1}^{X,\mathfrak m}\Big)\,\longrightarrow\,
{\rm H}^\bullet\Big(\mathfrak g_{2,\gamma_2}^{X,\mathfrak m}\Big)\,.
$$
\end{Lem}
\begin{proof}
We pick up a local chart $U$ of $X$ and write $D=\mathrm d+\omega|_U$. Then, we set $\widetilde{\gamma_1}=(\omega+\gamma_1)|_U$: 
$\widetilde{\gamma_1}$ is a MCE of $\mathfrak g_1^{U,\mathfrak m}$. 
We now introduce the DG algebra $\widetilde{\mathfrak{m}}:=\Omega(U)\otimes\mathfrak m$ and observe that 
\begin{itemize}
\item it is of the form $\widetilde{\mathfrak{m}}=C^\infty(U)\oplus\widetilde{\mathfrak{n}}$, with 
$\widetilde{\mathfrak{n}}$ (pro)nilpotent; 
\item $\mathfrak g_1^{U,\mathfrak m}\cong T_\mathrm{poly}^{\widetilde{\mathfrak m}}(V_\mathrm{formal})$ 
and $\mathfrak g_2^{U,\mathfrak m}\cong D_\mathrm{poly}^{\widetilde{\mathfrak m}}(V_\mathrm{formal})$; 
\item $\widetilde{\gamma_1}$ lies in the (pro)nilpotent part. 
\end{itemize}
We thus obtain a MCE $\mathcal U(\widetilde{\gamma_1})$ in $\mathfrak g_2^{U,\mathfrak m}$, 
which decomposes as $\mathcal U(\widetilde{\gamma_1})=(\omega+\gamma_2)|_U$ by the properties of $\mathfrak U_U$. 
We are therefore in the framework to which the compatibility between cup-products applies: the homotopy identity 
\eqref{eq-homotopy-U} holds true (for $\widetilde{\gamma_1}$ and $\widetilde{\mathfrak m}$). 

Finally, by construction we have $\mathcal U_{\widetilde{\gamma_1}}=\mathfrak U_{U,\gamma_1}$. It thus remains 
to prove that the homotopy operators, which are well-defined locally, glue together to a globally well-defined 
operator. This follows directly from the arguments of \cite[Lemma 10.1.1]{CVdB}. 
\end{proof}
Let us now prove that the first Taylor component of $\mathcal I_{\gamma_2}$ induces an isomorphism of Gerstenhaber 
algebras between ${\rm H}^\bullet\big(\mathfrak g_{2,\gamma_2}^X\big)$ and 
${\rm H}^\bullet\big(D_{\rm poly}^{\mathfrak m}(X)_{\gamma'}\big)$. Observe that 
we can view $\gamma'$ as a MCE in $\mathfrak g_2^X$ since we have an inclusion of $B_\infty$-algebras 
$D_{\rm poly}^{\mathfrak m}(X)\hookrightarrow\mathfrak g_2^{X,\mathfrak m}$. Moreover, since this inclusion is a 
quasi-isomorphism having $\mathcal I$ as a quasi-inverse then $\gamma_2$ and $\gamma'$ are gauge equivalent in 
$\mathfrak g_2^{X,\mathfrak m}$. 
As a consequence the $B_\infty$-algebras $\mathfrak g_{2,\gamma_2}^{X,\mathfrak m}$ and 
$\mathfrak g_{2,\gamma'}^{X,\mathfrak m}$ are isomorphic. 
In the end, we have the following commutative diagram of isomorphisms 
$$
\xymatrix{
{\rm H}^\bullet\Big(\mathfrak g_{2,\gamma_2}^{X,\mathfrak m}\Big) \ar[d]^-{\mathcal I_{\gamma_2,1}} \ar[r] & 
{\rm H}^\bullet\Big(\mathfrak g_{2,\gamma'}^{X,\mathfrak m}\Big) \\
{\rm H}^\bullet\Big(D_{\rm poly}^{\mathfrak m}(X)_{\gamma'}\Big) \ar[ru]\,, & }
$$
two of them being isomorphisms of Gerstenhaber algebras; so the third (i.e.~$\mathcal I_{\gamma_2}$) is. 

We have therefore proved that the first Taylor component $\mathcal U_{\gamma,1}$ of $\mathcal U_\gamma$ induces 
an isomorphism of (Gerstenhaber) algebras 
$$
{\rm H}^\bullet\Big(T_{\rm poly}^{\mathfrak m}(X)_\gamma\Big)\,\longrightarrow\,
{\rm H}^\bullet\Big(D_{\rm poly}^{\mathfrak m}(X)_{\mathcal U(\gamma)}\Big)
$$

\medskip

As for $\mathcal U_\gamma$, the quasi-isomorphism $\mathcal S_{\gamma,0}$ can be decomposed as follows: 
\begin{equation}\label{eq-SSS}
\xymatrix{C^{{\rm poly},\mathfrak m}(X)_{\gamma'} \ar@{^{(}->}[r] & 
\mathfrak M^{X,\mathfrak m}_{2,\gamma'} \ar[r] & 
\mathfrak M^{X,\mathfrak m}_{2,\gamma_2} \ar[r]^-{\mathfrak S_{X,\gamma_1}^0} & 
\mathfrak M_{1,\gamma_1}^{X,\mathfrak m} \ar[r]^-{\mathcal J_{\gamma,0}} & 
\Omega^{\mathfrak m}(X)_{\gamma}\,.}
\end{equation}
Finally, Theorem A is a consequence of the following 
\begin{Prop}
Sequences \eqref{eq-UUU} and \eqref{eq-SSS} of quasi-isomorphisms fit into a commutative diagram of 
Gerstenhaber algebras and $T$-modules
$$
\xymatrix{
{\rm H}^\bullet\Big(T_\mathrm{poly}^{\mathfrak m}(X)_\gamma\Big) \ar[r]\ar[d] & 
{\rm H}^\bullet\Big(\mathfrak g_{1,\gamma_1}^{X,\mathfrak m}\Big) \ar[r]\ar[d] & 
{\rm H}^\bullet\Big(\mathfrak g_{2,\gamma_2}^{X,\mathfrak m}\Big) \ar[r]\ar[d] & 
{\rm H}^\bullet\Big(\mathfrak g_{2,\gamma'}^{X,\mathfrak m}\Big) \ar[r]\ar[d] & 
{\rm H}^\bullet\Big(D_\mathrm{poly}^{\mathfrak m}(X)_{\gamma'}\Big) \ar[d] \\
{\rm H}^\bullet\Big(\Omega^{\mathfrak m}(X)_{\gamma}\Big) & 
{\rm H}^\bullet\Big(\mathfrak M^{X,\mathfrak m}_{1,\gamma_1}\Big) \ar[l] & 
{\rm H}^\bullet\Big(\mathfrak M_{2,\gamma_2}^{X,\mathfrak m}\Big) \ar[l] & 
{\rm H}^\bullet\Big(\mathfrak M^{X,\mathfrak m}_{2,\gamma'}\Big) \ar[l] & 
{\rm H}^\bullet\Big(C^{{\rm poly},\mathfrak m}(X)_{\gamma'}\Big) \ar[l]\,,}
$$
where the vertical arrows denote $T$-module structures and horizontal arrows 
isomorphisms of Gerstenhaber algebras and $T$-modules. 
\end{Prop}
\begin{proof}
We only have to prove that $\mathfrak S_{X,\gamma_1}^0$ induces an isomorphism of 
$T$-modules on cohomology. 

The proof goes along the same lines of arguments as in the proof of Lemma \ref{lem-UUU}. 
It therefore remains to prove that the homotopy operator for the compatibility between cap products 
satisfies the appropriate property so that gluing is allowed. Namely, we have the following 
\begin{Lem}
Let $\Gamma$ be an $S$-admissible graph in $\mathcal G_{n+1,m+1}^S$, $n\geq 1$, $\gamma_i$ general 
elements of $T_{\rm poly}(V_{\rm formal})$, such that, for some $j=2,\dots,n$, $\gamma_j$ is a linear vector field, $c$ a general element of $C^{\rm poly}$ of Hochschild degree $-m$. Then
\[
\overset{\circ}W_{D,\Gamma}\mathcal S_\Gamma(\gamma_1,\gamma_2,\dots,\gamma_n,c)=0.
\] 
\end{Lem}
\begin{proof}
We may assume without loss of generality that $\gamma_2$ is linear.

The first point of the first type in $\overset{\circ}{\mathcal Y}^+_{n+1,m+1}$, which moves along a fixed 
trajectory in the punctured unit disk, corresponds to $\gamma_1$; we observe that any other point in the punctured unit disk $D^\times$ can move freely. 

It is clear that $\mathcal S_\Gamma(\gamma_1,\gamma_2,\dots,c)$ does not 
trivially vanish, only if the vertex of $\Gamma$ corresponding to the linear vector field $\gamma_2$ has at most an incoming 
edge and exactly one outgoing edge: then, Lemma~\ref{l-K-var} yields the vanishing of 
the corresponding integral weight $\overset{\circ}W_{D,\Gamma}$. 
\end{proof}
Now, the homotopy operator for the compatibility between cap products is globally well-defined, 
because a change of local charts of $X$ changes the local expression for the Fedosov connection by a local 
$1$-form with values in linear vector fields: the previous identity implies that such a change does not in 
fact contribute, whence the Proposition follows.  
\end{proof}

\subsection{Relation with deformation quantization and the work of Cattaneo-Felder-Tomassini}

We consider again the special case when $\mathfrak m=\mathbb{R}[\![\hbar]\!]$, and $\gamma$ is a MCE of the form 
$$
\gamma=\pi_\hbar=\hbar\pi_1+\mathcal O(\hbar^2)\in\hbar T_{\rm poly}^1(X)[\![\hbar]\!]\,.
$$
In particular $\pi_1$ defines a Poisson structure on $X$. Below we assume the reader is familiar with 
the subject of deformation quantization. 

Borrowing the notation from the proof of Lemma \ref{lem-UUU} in Subsection \ref{ss-8-1bis}, we see that w.r.t.\ a 
local chart $U$ of $X$, $\widetilde{\gamma_1}=(\omega+\gamma_1)|_U$ has polyvector degree less or equal to $1$, and 
$D(\gamma_1)=0$. 

According to the computations of Subsection~\ref{ss-6-2}, $\mathcal U(\widetilde{\gamma_1})$ induces the following 
structure on 
$A^U:=\mathcal O_{V_{\rm formal}}\otimes\widetilde{\mathfrak m}^0=C^\infty(U,\mathcal O_{V_{\rm formal}})[\![\hbar]\!]$: 
$i)$ a $C^\infty(U)[\![\hbar]\!]$-linear associative product $\star$ on $A^U$, 
$ii)$ a Fedosov connection $\widetilde Q$ on $(A^U,\star)$ and 
$iii)$ the Weyl curvature $F$ of $\widetilde Q$, in the terminology of Subsection 4.2,~\cite{CFT}. 

By inspecting~\eqref{eq-MCtwist} and recalling that $\pi_\hbar=\mathcal O(\hbar)$, it 
follows that $i)$ the Fedosov connection $\widetilde Q$ satisfies $\widetilde Q=D+\mathcal O(\hbar)$, where 
$D$ is the previously introduced Fedosov connection on $\mathfrak g_1^X$, and $ii)$ the Weyl curvature $F$ 
of $\widetilde Q$ satisfies $F=\mathcal O(\hbar)$.

We recall now that $D$ is flat and that the corresponding cohomology on $\mathfrak g_2^X[\![\hbar]\!]$ is concentrated 
in degree $0$: then, the assumptions of Lemmata 4.5 and 4.6 in \cite{CFT} are satisfied, whence we may 
recover\footnote{As a particular case of our claim, in Subsection \ref{ss-8-1bis}, that $\gamma'$ and 
$\gamma_2$ are gauge equivalent in $\mathfrak g_2^{X,\mathfrak m}$. } 
the fact that the MCE in $\mathfrak g_2^X[\![\hbar]\!]$ corresponding to the triple $(\star,\widetilde{Q},F)$ 
is gauge-equivalent to the MCE corresponding to the product $\star$ and the underformed connection $D$, 
with zero curvature. 

This last fact is of crucial use in \cite{CFT} to prove that one has an isomorphism of algebras between 
the algebra of Casimir functions for the formal Poisson structure $\pi_\hbar$ and the center of the corresponding 
quantized algebra, which is precisely the degree zero part of the compatibility between cap products on 
cohomology. 

Let us simply recall that the quantized algebra is constructed as the subspace of $D$-flat sections in $A^X$, 
which is isomorphic to $C^\infty(X)[\![\hbar]\!]$, equipped with the associative product $\star$.

\section{Appendix}\label{app}

In this Appendix, we quote two main technical Lemmata from~\cite{K}, which are used in many computations throughout the paper.  

First of all, we consider the compactified configuration space $\mathcal C_n$, with $n\geq 3$; further, for any two distinct indices $1\leq i\neq j\leq n$, there is natural projection $\pi_{ij}$ from $\mathcal C_n$ onto $\mathcal C_2$, and we denote by $\omega_{ij}$ the pull-back of $\omega|_{\mathcal C_2}$ w.r.t.\ the projection $\pi_{ij}$ (see Lemma~\ref{l-angle}, Subsection~\ref{ss-3-1}). 
\begin{Lem}\label{l-K-1}
For a positive integer $n\geq 3$, the integral
\[
\int_{\mathcal C_n} \bigwedge_{\alpha=1}^{2n-3} \omega_{i_\alpha j_\alpha}
\]
vanishes, for any distinct indices $1\leq i_\alpha\neq j_\alpha\leq n$, $\alpha=1,\dots,2n-3$.
\end{Lem}
For a proof of Lemma~\ref{l-K-1}, we refer to~\cite{K}, Subsection 6.6; in~\cite{BCKT} one can find an alternative proof to the original one of Kontsevich.
Lemma~\ref{l-K-1} is often used in Subsection~\ref{ss-4-2}, Subsection~\ref{ss-5-2} and related Subsubsections, and in Subsection~\ref{ss-5-4} and related Subsubsections. 

We now consider the $1$-form $\omega$ on Kontsevich's eye $\mathcal C_{2,0}$ as in Subsection~\ref{ss-3-1}.
\begin{Lem}\label{l-K-2}
The integral 
\[
\int_{\mathcal C_{2,0}} \omega_{12}\wedge\omega_{21}=\int_{\mathcal H\smallsetminus \{z_0\}}\omega(z_0,z)\wedge \omega(z,z_0)
\]
vanishes, where $\omega_{12}$, resp.\ $\omega_{21}$, denotes the usual form $\omega$ on $\mathcal C_{2,0}$, resp.\ the usual form $\omega$, but with the arguments permuted (we observe that Kontsevich's angle function is not symmetric w.r.t.\ its arguments, hence the vanishing of the integral is non-trivial). 
In the second expression, $z_0$ is some fixed point in the complex upper half-plane $\mathcal H$.
\end{Lem}
\begin{Lem}\label{l-K-3}
If $z_1$, $z_2$ are either two distinct points in the complex upper half-plane $\mathcal H$ or if $z_1$ is in the complex upper half-plane $\mathcal H$ and $z_2$ on the real axis $\mathbb R$, the integral
\[
\int_{\mathcal H\smallsetminus\{z_1,z_2\}} \omega(z_1,z)\wedge \omega(z,z_2)
\]
vanishes, where we integrate w.r.t.\ $z$.
\end{Lem}
For a proof of both Lemmata~\ref{l-K-2} and~\ref{l-K-3}, we refer to~\cite{K}, Lemmata 7.3, 7.4 and 7.5, or to~\cite{CR}.

\end{document}